\documentclass{amsart}
\usepackage{graphicx}
\usepackage{amssymb}
\usepackage{epstopdf}
\usepackage{nicefrac}
\usepackage{enumerate}
\usepackage{hyperref}
\usepackage{float}
\usepackage{color}
\usepackage{pst-all}
\usepackage{tabularx}

\newtheorem{thm}{THEOREM}[section]

\newtheorem{cor}[thm]{COROLLARY}

\newtheorem{defn}[thm]{DEFINITION}

\newtheorem{lemma}[thm]{LEMMA}
\newtheorem{prob}[thm]{PROBLEM}
\newtheorem{prop}[thm]{PROPOSITION}

\newtheorem{remark}[thm]{REMARK}

\newcommand{\ds}{\displaystyle}

\newcommand{\cB}{{\mathcal B}}

\newcommand{\cD}{{\mathcal D}}
\newcommand{\cG}{{\mathcal G}}
\newcommand{\cGF}{\cG_{\F}} 
\newcommand{\cH}{{\mathcal H}}

\newcommand{\cK}{{\mathcal K}}

\newcommand{\cO}{{\mathcal O}}

\newcommand{\cR}{{\mathcal R}}
\newcommand{\cZ}{{\mathcal Z}}
\newcommand{\CO}{{\rm CO}} 
\newcommand{\cP}{{\mathcal P}}

\newcommand{\cU}{{\mathcal U}}
\newcommand{\cV}{{\mathcal V}}
\newcommand{\cW}{{\mathcal W}}
\newcommand{\diam}{{\rm diam}} 
\newcommand{\dX}{d_{\fX}} 
\newcommand{\e}{{\varepsilon}} 
\newcommand{\F}{{\mathcal F}}
\newcommand{\fG}{{\mathfrak{G}}}
\newcommand{\Fix}{{\rm Fix}}

\newcommand{\fX}{{\mathfrak{X}}}

\newcommand{\fZ}{{\mathfrak{Z}}}
\newcommand{\Iso}{{\rm Iso}} 

\newcommand{\mN}{{\mathbb N}}
\newcommand{\mS}{{\mathbb S}}
\newcommand{\mZ}{{\mathbb Z}}
\newcommand{\vp}{{\varphi}}
\newcommand{\whalpha}{\widehat{\alpha}}

\newcommand{\wtG}{\widetilde{G}}

 \newcommand{\Homeo}{{\rm Homeo}} 
 \newcommand{\Cen}{{\rm Cen}} 
  
\newcommand{\Aut}{{\rm Aut}}

\newcommand{\Ad}{{\bf Ad}} 

\newcommand{\whPhi}{\widehat{\Phi}}

\newcommand{\whe}{\widehat{e}}
\newcommand{\whg}{\widehat{g}}
\newcommand{\whh}{\widehat{h}}
\newcommand{\whk}{\widehat{k}}
\newcommand{\whq}{\widehat{q}}

\newcommand{\why}{\widehat{y}}

\newcommand{\whC}{\widehat{C}}
\newcommand{\whG}{\widehat{G}}
\newcommand{\whH}{\widehat{H}}
\newcommand{\whK}{\widehat{K}}

\newcommand{\whU}{\widehat{U}}
\newcommand{\whV}{\widehat{V}}

\newcommand{\oG}{{\mathfrak{G}(\Phi)}}

\newcommand{\ocD}{\overline{\cD}}
\newcommand{\whtheta}{{\widehat{\theta}}}
\newcommand{\whTheta}{{\widehat{\Theta}}}

\newcommand{\whtau}{{\widehat{\tau}}}

 \newcommand{\oK}{{\overline{K}}}
 \newcommand{\ocK}{{\overline{\cK}}}
 \newcommand{\orho}{{\overline{\rho}}}
 \newcommand{\oPhi}{{\overline{\Phi}}}
 \newcommand{\oPsi}{{\overline{\Psi}}}

\newcommand{\ovh}{{\overline{h}}}

\newcommand{\G}{\Gamma}

\input xy
\xyoption {all}

\begin{document}

\title{Limit group invariants  for non-free Cantor actions}

\author{Steven Hurder}
\address{Steven Hurder, Department of Mathematics, University of Illinois at Chicago, 322 SEO (m/c 249), 851 S. Morgan Street, Chicago, IL 60607-7045}
\email{hurder@uic.edu}

\author{Olga Lukina}
\address{Olga Lukina, Faculty of Mathematics, University of Vienna, Oskar-Morgenstern-Platz 1, 1090 Vienna, Austria}
\email{olga.lukina@univie.ac.at}

 \thanks{OL is supported by FWF Project P31950-N35}

 
\thanks{Version date: April 25, 2019; revised January 9, 2020}

  \begin{abstract}
  A Cantor action is a minimal equicontinuous action of a countably generated group G on a Cantor space X. Such actions are also called generalized odometers in the literature. In this work, we introduce two new conjugacy invariants for Cantor actions, the stabilizer limit group and the centralizer limit group.   An action is wild if the stabilizer limit group is an increasing sequence of stabilizer groups without bound, and otherwise is said to be stable if this group chain is bounded. For Cantor actions by a finitely generated group G, we prove that stable actions satisfy a rigidity principle, and furthermore show that the wild property is an invariant of the continuous orbit equivalence class of the action.  
  
 A Cantor action is said to be dynamically wild if it is wild, and the  centralizer limit group is a proper subgroup of the stabilizer limit group. This property is also a conjugacy invariant, and we show    that a Cantor action with a non-Hausdorff element   must be dynamically wild. We then give     examples of wild Cantor actions with non-Hausdorff elements,   using     recursive methods   from Geometric Group Theory to define  actions on the boundaries of trees.
 
\end{abstract}

\maketitle

\section{Introduction}\label{sec-intro}

In this paper, we  investigate the structure of non-free Cantor group actions, and the relationship between the dynamics of the action and the algebraic properties of the group. One of our main results is the definition of direct limit groups which are conjugacy invariants of the action, and we investigate the relation between these new invariants and the dynamics of the action.  Our results are illustrated by   examples of Cantor actions in the literature, and also those constructed  in Section~\ref{sec-arboreal}.

We assume as   standing hypotheses that 
 $G$ is a countably infinite group, $\fX$ is a Cantor space, and $\Phi \colon G \to {\Homeo}(\fX)$   is  an action of $G$.  
  We sometimes assume in addition that $G$ is finitely generated, and when required this hypothesis will be indicated. 
We   denote the action by $(\fX,G,\Phi)$, and write $g \cdot x$ for   $\Phi(g)(x)$  when the action is clear.      

The action   is \emph{minimal} if  for all $x \in \fX$, its   orbit $\cO(x)   = \{g \cdot x \mid g \in G\}$ is dense in $\fX$.  
Let  $\dX$ be a metric on $\fX$ compatible with the topology.   The action  $(\fX,G,\Phi)$ is \emph{equicontinuous} with respect to  $\dX$  if for all $\e >0$ there exists $\delta > 0$, such that for all $x , y \in \fX$ and $g \in G$, $ d_X(x,y) < \delta$ implies that  $d_X(\Phi(g)  \cdot x, \Phi(g) \cdot y) < \e$.
This property   is independent of the choice of the metric $\dX$ on $\fX$. 

 In this paper, a \emph{Cantor action} $(\fX,G,\Phi)$ is     assumed to be minimal and equicontinuous. 
 We also occasionally discuss actions which are not assumed to be minimal or equicontinuous, and these will be called   \emph{general  Cantor actions}.
 
 Cantor actions can be divided into two types, \emph{stable} and \emph{wild}.  For example, a free Cantor action of an abelian group $G$ is stable.
 The stable property  is a weaker form of the more well-known ``topologically free'' property assumed in many works on the study of Cantor actions.  Topologically free general Cantor actions have been   extensively analyzed in the literature \cite{BoyleTomiyama1998,CortezMedynets2016,Li2018,Renault2008}, and the authors' work in \cite{DHL2016c,HL2018b} extends some of these results to stable Cantor actions. The distinction between stable and topologically free actions is discussed in Section~\ref{sec-rigidity}. 
 On the other hand, the structure and properties of wild Cantor actions are less well-known,  and in this work we study this class of actions in depth, obtaining  new invariants and new classification results.

The first example of a wild Cantor action was constructed by Schori in \cite{Schori1966}, as the monodromy action for a weak solenoid which is not homogeneous.  The ``Schori solenoid''  gave a counter-example to a question raised by McCord in \cite{McCord1965}, and became the focus of further study. The   analysis of this example in the authors' work with Jessica Dyer \cite{DHL2016c} introduced the notion of the wild property for a minimal equicontinuous action of a   group $G$ on a Cantor space $\fX$. It was     observed by the authors in \cite{DHL2016c} that the Molino Theory for foliated spaces developed in  \cite{ALM2016} applies for weak solenoids whose monodromy actions are stable,   but fails for weak solenoids whose monodromy actions are  wild. 

The   work \cite{HL2018a} also gave a general method for constructing wild Cantor actions, and used the wild property to show that these examples yield uncountably many classes of non-homeomorphic, non-homogeneous weak solenoids. 
This work suggested that a classification of weak solenoids up to homeomorphism requires a better understanding of the class of wild Cantor actions.
 
In the subsequent work \cite{HL2018b}, the authors studied the classification problem for stable Cantor actions, and gave a generalization of the rigidity theorems for Cantor actions of Cortez and Medynets \cite{CortezMedynets2016} and Li \cite{Li2018}. We showed that for stable Cantor actions, continuous orbit equivalence implies return equivalence of the actions. Our work also showed that  if $G$ is a finitely generated nilpotent group, then every Cantor action by $G$ is stable. Thus, for Cantor actions by finitely generated nilpotent groups, one has a direct approach to their classification via the rigidity principle. 

The study of stable Cantor actions also highlighted the role of non-Hausdorff elements for the action, as discussed in Section~\ref{subsec-nHclosure} and     Definition~\ref{def-nHelement}.
 An action with a non-Hausdorff element is at the opposite extreme of a stable action. Another theme of our work is to investigate     wild actions with non-Hausdorff elements.

 We now discuss the main results of this work. Sections~\ref{sec-models} and \ref{sec-odometer} give the basic background material on Cantor actions that we require. This material overlaps with discussions and results in the authors' previous papers \cite{DHL2016a,HL2018a,HL2018b}, and somewhat also with the works \cite{CortezMedynets2016,Li2018}, but is discussed here in sufficient detail as necessary for the remainder of the paper. 
 
 Section~\ref{sec-models} describes a model for a Cantor action  $(\fX,G,\Phi)$ as a group action on a Cantor homogeneous space, 
 $\fX \cong \oG/\oG_x$, where 
 $\oG \subset \Homeo(\fX)$ is a profinite group   acting transitively on $\fX$, and   $\oG_x \subset \oG$ is the isotropy subgroup of a point $x \in \fX$. The  philosophy of our approach in this paper,   is that the study of the ``adjoint action''  of the isotropy group   $\oG_x$ on $\oG$  yields  conjugacy  invariants for the action. 
 
   Section~\ref{sec-odometer} discusses the odometer model for a Cantor action, and introduces a number of concepts about these actions which are key to subsequent sections. In particular, the stabilizer group chain of an odometer is defined in \eqref{eq-isotropykernel}, and   in Section~\ref{subsec-adjoints}   the     the centralizer group chain is defined.

Section~\ref{sec-limgroups} recalls the formal construction of direct limit groups,   their equivalence, and properties of these groups. Then  Definition~\ref{def-progroupinvariants} and   Theorem~\ref{thm-centraldirectlimgroup}   combined yield:

   \begin{thm}\label{thm-main1}
 Let    $(\fX,G,\Phi)$ be a Cantor action. There are well-defined   direct limit groups 
  \begin{enumerate}
\item the \emph{stabilizer limit group}     $\ds \Upsilon_s(\Phi)$,   
\item the \emph{centralizer limit group}      $\ds \Upsilon_c(\Phi)$,
\end{enumerate}
  defined as the equivalence classes of     the stabilizer and centralizer group chains associated to an odometer model for the action. 
Both    $\Upsilon_s(\Phi)$ and     $\Upsilon_c(\Phi)$  are    conjugacy class invariants of the action, and there is an inclusion of direct limit groups $\ds \Upsilon_c(\Phi) \subset   \Upsilon_s(\Phi)$.
 \end{thm}

  For a Cantor action of an abelian group $G = \mZ^d$, both of the limit group invariants in Theorem~\ref{thm-main1} are   trivial, and the action is stable. However, for Cantor actions of more general groups $G$, the stabilizer and centralizer limit group invariants can be highly non-trivial. One theme of this work is to use these limit group invariants  to obtain a finer classification of non-free Cantor actions,    into the following subtypes of actions. 
  \eject
  
  \begin{defn} [Definition~\ref{def-wildtypes}] \label{def-main1}
A Cantor action     $(\fX,G,\Phi)$ is said to be: 
\begin{enumerate}
\item \emph{stable} if    the stabilizer  group $\Upsilon_s(\Phi)$   is bounded, and   \emph{wild} otherwise;
\item  \emph{algebraically stable} if the  its centralizer  group $\Upsilon_c(\Phi)$   is bounded, and   \emph{algebraically wild} otherwise;
\item \emph{wild of finite type}  if    the stabilizer  group $\Upsilon_s(\Phi)$   is unbounded, and represented by a chain of finite groups;
\item \emph{wild of flat type}  if    the stabilizer  group $\Upsilon_s(\Phi)$   is unbounded, and $\Upsilon_c(\Phi) = \Upsilon_s(\Phi)$; 
\item \emph{dynamically wild}  if    the stabilizer  group $\Upsilon_s(\Phi)$   is unbounded, and is not of flat type.
\end{enumerate}
  \end{defn}
   All of the possibilities in this definition can be realized by examples of Cantor actions.

Section~\ref{sec-dynamics} discusses the notion of a \emph{locally quasi-analytic} Cantor action $(\fX,G,\Phi)$, and the relation between this notion and the more usual notion of a topologically free action. The main result gives an interpretation of the stable property  in terms of the ``analytic properties'' of the adjoint action of the isotropy group $\oG_x \subset \oG$ for a homogeneous model for $\fX \cong \oG/\oG_x$.  
\begin{thm}[Theorem~\ref{thm-boundedLCQA}]\label{thm-main2}
Let        $(\fX,G,\Phi)$  be a Cantor action, then $\Phi$ is stable if and only if the action of $\oG$ on $\fX$ is   locally   quasi-analytic.
\end{thm}

Section~\ref{sec-rigidity} recalls the notion of a continuous orbit equivalence (COE) between Cantor actions, then discusses three notions of rigidity for Cantor actions.  The work of Li  in \cite{Li2018} gives cohomological criteria for when two COE Cantor actions, both of which are topologically free actions,  are necessarily  $\theta$-conjugate, as defined in Definition~\ref{def-thetaconjugacy}.  The work of Cortez and Medynets in  \cite{CortezMedynets2016} shows that two COE Cantor actions, both of which are free, are   
structurally conjugate (or virtually rigid as defined in Definition~\ref{def-Vrigidity}). The notion of return equivalence for Cantor actions was introduced in  \cite{CHL2018a}, and the authors showed in \cite{HL2018b} that  two COE Cantor actions, both of which are stable, are necessarily return equivalent  as   in Definition~\ref{def-return}.

 It is thus natural to ask how the invariants  $\Upsilon_c(\Phi)$ and $\Upsilon_s(\Phi)$ behave under orbit equivalence. Our first result, and technically most involved, shows the following:
\begin{thm}[Theorem~\ref{thm-coestable}]\label{thm-main3}
Let $(\fX,G,\Phi)$  and $(\fX', G', \Psi)$ be continuously orbit equivalent Cantor actions.  If $G$  is finitely generated, and  $(\fX', G', \Psi)$ is stable, then  $(\fX,G,\Phi)$ is stable.
  \end{thm}
  We deduce from this a strong form of invariance for the wild property. 
\begin{cor}\label{cor-main3}
The wild property is an invariant of continuous orbit equivalence for the class of Cantor actions by finitely generated groups.
\end{cor}
Combining  Corollary~\ref{cor-main3} with the results in    \cite{HL2018b}  yields:
  \begin{cor}[Theorem~\ref{thm-coe=re}]\label{thm-coe=re2}
   Let $G$ and $G'$ be a finitely generated groups, and suppose that the Cantor action   $(\fX', G', \Psi)$ is stable. Let $(\fX, G, \Phi)$ be a general Cantor action   which is  continuously orbit equivalent to $(\fX', G', \Psi)$,  then the actions are return equivalent. 
   \end{cor}

The   above results give effective approaches to   the classification of stable Cantor actions, up to conjugacy or continuous orbit equivalence. On the other hand, for wild Cantor actions, much less is known. The remainder of this work then  investigates the properties of wild Cantor actions.

Section~\ref{sec-hausdorff} introduces the notion of a non-Hausdorff element in $\oG$ for a Cantor action $(\fX, G, \Phi)$. It was remarked by Renault in \cite{Renault2008} that if the germinal groupoid associated to a group action is non-Hausdorff, then the action cannot be topologically free.
In the authors' work \cite{HL2018b} this result was extended to the observation that if $\oG$ contains a non-Hausdorff element, then the action is wild. In this work, we give a more precise consequence. The stabilizer limit group $\Upsilon_s(\Phi)$ is said to be \emph{of finite type} if each group in a representative group chain is a finite group.  Then we have the following result, which relates the direct limit invariants and other ideas introduced in this paper with the dynamics of a Cantor action:
 \begin{thm}\label{thm-main4}
  Let $(\fX,G,\Phi)$  be a  Cantor action. Suppose that $\oG$ contains a non-Hausdorff element, then the  action is dynamically wild and not of finite type.  That is, there is a \emph{proper} inclusion of direct limit groups $\Upsilon_c(\Phi) \subset \Upsilon_s(\Phi)$, and $\Upsilon_s(\Phi)$ is represented by an increasing chain of Cantor groups, which in particular are uncountable.
    \end{thm}
The claims of Theorem~\ref{thm-main4} follow from   Corollary~\ref{cor-nH=unbounded} and Theorem~\ref{thm-distinct=nH}.
\begin{cor}\label{cor-Hgroupoid}
 Let $(\fX, G, \Phi)$ be a Cantor action   for which the stabilizer direct limit group $\Upsilon_s(\Phi)$ has finite type, then the germinal groupoid $\cG(\fX, G, \Phi)$ associated to the action is Hausdorff.
\end{cor}
 The examples of wild Cantor actions  given in \cite[Section~8]{HL2018a} have finite type, so Corollary~\ref{cor-Hgroupoid} implies there exists wild Cantor actions with Hausdorff germinal groupoids.
 
It seems to be a difficult problem to give criteria for a dynamically wild Cantor action which suffice to imply the existence of a non-Hausdorff element.    
Section~\ref{sec-arboreal} constructs  examples of wild Cantor actions  using the ``automata''  approach, which is a well-known method in Geometric Group Theory. This method  defines   a homeomorphism of the   boundary of a tree, using a  recursive definition along the branches of the tree. The examples constructed are   inspired by the work of Nekrashevych in   \cite{Nekrashevych2005,Nekrashevych2018}, and Pink in \cite{Pink2013}.
It would be very interesting to know if other methods of construction are possible, and perhaps that  the presence of a non-Hausdorff element for a Cantor action implies some form of underlying recursiveness for the  action of its generators.

We conclude with the two general problems most relevant to the results of this work. 
     \begin{prob}\label{prob-nh}
 Give sufficient conditions for a pair $(\fG, \cD)$ to admit a non-Hausdorff element for the action of $\fG$ on $\fX \equiv \fG/\cD$, where $\fG$ is a profinite group which is finitely generated in the topological sense, and $\cD \subset \fG$ is a ``totally not normal'' closed subgroup.  
  \end{prob}
  The notion of a \emph{totally not normal} subgroup  is given in Definition~\ref{def-tnn}.

 Note that this problem is not just about the algebraic properties of the groups, as the wild property only emerges when considering the transitive action of $\fG$ on the quotient space $\fX$. 
 
 For the examples of wild Cantor actions of flat type given in \cite[Section~8]{HL2018a},  the discriminant groups of these actions are  an infinite product of finite groups. The second question asks whether there is a   general result,   that    the discriminant groups of flat actions have a restricted algebraic structure.  
     \begin{prob}\label{prob-flat}
Let $(\fX,G,\Phi)$ be a wild Cantor action of flat type, so   the inclusion of direct limit groups $\Upsilon_c(\Phi) \subset \Upsilon_s(\Phi)$ is an equality. What restrictions are imposed on  its discriminant group?
  \end{prob}

\section{The profinite model}\label{sec-models}

Given a Cantor action  $(\fX,G,\Phi)$, let  $\Phi(G) \subset Homeo(X)$ denote the image subgroup. Introduce the  closure $\oG \equiv \overline{\Phi(G)} \subset Homeo(X)$ for the \emph{uniform topology of maps}. (This corresponds to the \emph{Ellis group} for the action, as defined in \cite{Auslander1988,Ellis1960,EllisGottschalk1960}; see also \cite[Section~2]{DHL2016a}.)
 That is, each  element   $\whg \in \fG(\Phi)$ is the uniform limit of a sequence of maps   $\{\Phi(g_i) \mid i \geq 1\} \subset \Phi(G)$. By abuse of notation, we    sometimes also denote   the limiting element $\whg$ by $(g_i)$.
It was observed by Ellis in \cite{Ellis1960},   by Auslander in \cite[Chapter 3]{Auslander1988} and again by Glasner   in \cite[Proposition~2.5]{Glasner2007}, that   the assumption   the action is equicontinuous implies    $\fG(\Phi) $  is a separable profinite group. 

For example, if $G$ is an abelian group, then $\fG(\Phi)$ is a compact totally disconnected abelian group, which can be thought of as the group of asymptotic motions of the system. When $G$ is non-abelian,    the action closure $\oG$ can have   more subtle algebraic properties. A main theme of this work is to examine the interplay between the algebraic structure of $\oG$ and the dynamics of the action $\Phi$.

     The philosophy behind our     study of Cantor actions, is to consider $\fX$ as a homogeneous space for $\oG$, in analogy to the case of homogeneous spaces of Lie type. Recall that for $G$ a connected Lie group and $H \subset G$ a closed subgroup,   the quotient space $X = G/H$ is a homogeneous $G$-space, and one studies the geometry of $X$ using the structure of the Lie algebra $\mathfrak{g}$ of $G$, and the adjoint action of $H$ on  $\mathfrak{g}$. See for example   Chapters X and XI of \cite{KN1969}.     
     For a  Cantor action $(\fX,G,\Phi)$, there is no obvious analog of   a Lie algebra associated to $\oG$. None-the-less, one can investigate the properties of the adjoint action of the isotropy group of a point, localized to neighborhoods of $x$.

Let $\whPhi \colon \oG \times \fX \to \fX$ denote the induced action of $\oG$ on $\fX$. For $\whg \in \oG$ we   write its action on $\fX$ by $\whg \cdot x = \whPhi(\whg)(x)$.
If the action $\Phi \colon G \times \fX \to \fX$ is   minimal, then the group $\oG$ acts transitively on $\fX$.  The action $(\fX,G,\Phi)$ is said to be \emph{faithful} if $\Phi(g) \cdot x = x$ for all $x \in \fX$ implies that $g$ is the identity element. Equivalently, the action is faithful if the action map $\Phi \colon G \to \Homeo(\fX)$ is injective.
The action   is   \emph{free}   if for any $x \in \fX$,    $\Phi(g) \cdot x = x$ implies that $g$ is the identity element of $G$.

\subsection{The adjoint action}
Given $x \in \fX$,   introduce the isotropy group  at $x$,  
\begin{align}\label{iso-defn2}
 \oG_x = \{ \whg  \in \overline{\Phi(G)} \mid \whg \cdot x = x\} \subset \Homeo(\fX) \ ,
\end{align}
which is a closed subgroup of $\oG$, and is  thus either finite, or is a Cantor group.    

There is a   natural identification $\fX \cong \oG/\oG_x$ of left $G$-spaces, and thus      
the conjugacy class of $\oG_x$ in $\oG$ is independent of the choice of   $x$. Moreover, if $\oG_x$ is the trivial group, then $\fX$ is identified with a profinite group, and the   action is  free.  However, there exists examples of  free Cantor actions for which the group $\oG_x$ is non-trivial; the first such examples were constructed by  Fokkink and Oversteegen   in \cite[Section~8]{FO2002}, and further examples were constructed in   \cite[Section~10]{DHL2016c}.

The action of $\oG_x$ on $\fX$  can be considered as induced by the adjoint action of $\oG_x$ on $\oG$. For $\whh \in \oG_x$ and $\whg \in \oG$, let $\Ad(\whh)(\whg)  = \whh \ \whg \ \whh^{-1}$. For $y \in \fX$ choose $\whg \in \oG$ such that $y = \whg \cdot \oG_x$. Then for $\whh \in \oG_x$ we have 
\begin{equation}\label{eq-adjoint}
 \whh  \cdot y = \whh \  \whg \cdot \oG_x = \whh \  \whg \ \whh^{-1}   \cdot \oG_x = \Ad(\whh)(\whg) \cdot \oG_x ~.
 \end{equation}
That is, the action of $\oG_x$ on $\fX$ can be considered as the factor of the adjoint action which is  ``normal'' to $\oG_x$.
 In the case of a homogeneous space  of Lie type, this normal action is induced by the adjoint action on the Lie algebra of $\oG$, and is studied  in terms of representation theory of the compact group $\oG_x$.  For a Cantor action, we instead   consider the restriction of the adjoint action to arbitrarily small neighborhoods of $x \in \fX$. If this localized action of $\oG_x$ stabilizes for arbitrarily small neighborhoods of $x$, then the action is said to be \emph{stable}, and otherwise it is \emph{wild}. 
That is, for a stable Cantor action, there is a well-defined local geometric model for $\fX$ near $x$, while for a wild action there is no stable local model. We  next make these statements precise.

\subsection{The finite model}\label{subsec-models}
A profinite group, by definition, is the inverse limit  of finite quotient groups. For an equicontinuous action $(\fX,G,\Phi)$, an analogous statement holds, that the action $\Phi$ is defined by an inverse limit of finite actions. The key concept to show this is that of  clopen subsets of $\fX$ which are adapted to the action. We briefly recall some basic concepts.

Let $\CO(\fX)$ denote the collection  of all clopen (closed and open) subsets of the Cantor space $\fX$, which forms a basis for the topology of $\fX$. 
For $\phi \in \Homeo(\fX)$ and    $U \in \CO(\fX)$, the image $\phi(U) \in \CO(\fX)$.  
The following   result is folklore, and a proof is given in \cite[Proposition~3.1]{HL2018b}.
 \begin{prop}\label{prop-CO}
A minimal Cantor action $(\fX,G,\Phi)$ is  equicontinuous  if and only if, for the induced action $\Phi_* \colon G \times \CO(\fX) \to \CO(\fX)$, the $G$-orbit of every $U \in \CO(\fX)$ is finite.
\end{prop}
 The   proof that each $U \in \CO(\fX)$ has finite orbit is essentially the same as what was called the ``coding method'' for equicontinuous pseudogroup actions on Cantor transversals in \cite{ClarkHurder2013}, and discussed for group actions    
in \cite[Appendix A]{DHL2016a}, and for free actions in \cite[Section~2]{CortezMedynets2016}.

We say that $U \subset \fX$  is \emph{adapted} to the action $(\fX,G,\Phi)$ if $U$ is a   \emph{non-empty clopen} subset, and for any $g \in G$, 
if $\Phi(g)(U) \cap U \ne \emptyset$ implies that  $\Phi(g)(U) = U$.   The proof of  Proposition~3.1 in \cite{HL2018b} shows that given  $x \in \fX$ and clopen set $x \in W$, there is an adapted clopen set $U$ with $x \in U \subset W$. 

The key property of adapted sets, is that   for  $U$   adapted,   the set of ``return times'' to $U$, 
 \begin{equation}\label{eq-adapted}
G_U = \left\{g \in G \mid \vp(g)(U) \cap U \ne \emptyset  \right\}  
\end{equation}
is a subgroup of   $G$, called the \emph{stabilizer} of $U$.      
  Then for $g, g' \in G$ with $g \cdot U \cap g' \cdot U \ne \emptyset$ we have $g^{-1} \, g' \cdot U = U$, hence $g^{-1} \, g' \in G_U$. Thus,  the  translates $\{ g \cdot U \mid g \in G\}$ form a finite clopen partition of $\fX$, and are in 1-1 correspondence with the quotient space $X_U = G/G_U$. Then $G$ acts by permutations of the finite set $X_U$ and so the stabilizer group $G_U \subset G$ has finite index.  
The action of $g \in G$ on $X_U$ is trivial precisely when $g$ is a stabilizer of each coset $h \cdot G_U$, so   $g \in C_U$ where
\begin{equation}\label{eq-CUdef}
 C_U \ = \  \bigcap_{h \in G} \ h  \ G_U \ h^{-1} \ \subset \ G_U  
\end{equation}
  is the largest normal subgroup of $G$ contained in $G_U$. 
The action of    the finite group $Q_U \equiv G/C_U$ on $X_U$ by permutations is a finite approximation of the action of $G$ on $\fX$,  and the isotropy group of the identity coset   $e \cdot G_U$ is  $D_U \equiv G_U/C_U \subset Q_U$. 

 Thus, the finite model of a Cantor action $(\fX,G,\Phi)$ obtained from $U \in \CO(\fX)$ yields the following data:  The   group $Q_U$ which acts on $X_U$, and the isotropy group $D_U \subset Q_U$ for the basepoint defined by the identity coset $x_U = e \cdot G_U \in X_U$. In addition to the group structures of $Q_U$ and $D_U$, one can study how the group $D_U$ is ``algebraically embedded'' in $Q_U$. 
 
Consider the adjoint action of $D_U$ on $Q_U$, which for  $h \in D_U$ is given by 
  \begin{equation}\label{eq-Urestrictedadjoint}
\Ad_{Q_U}(h) \colon Q_U \to Q_U \quad , \quad \Ad_{Q_U}(h)(k) = h \ k \ h^{-1} ~ ,   ~ k \in Q_U ~ .
\end{equation}

The subgroup $D_U$ is not  a normal subgroup in $Q_U$, and in fact a stronger condition holds. Let $h \in D_U$ and suppose that  $k \ h \ k^{-1} \in D_U$ for all $k \in Q_U$.  Then $k \ h \ k^{-1} = h'$ for some $h' \in D_U$, so
$h k^{-1} = k^{-1} \ h' = k^{-1} \ h'' \ h$ for $h'' = h' \ h^{-1}$. Hence $h \ k^{-1} \ h^{-1} = k^{-1} \ h''$. As this holds for all $k \in Q_U$  the adjoint  $\Ad_{Q_U}(h)$ acts trivially on $X_U$ and hence $h$ must be the identity in $Q_U$.  In other words, for every non-trivial $h \in D_U$, there is a $k \in Q_U$ such that  $k \ h \ k^{-1} \not\in D_U$.  This shows that the subgroup $D_U$ satisfies the following condition, which was introduced in the works \cite{DHL2016a,HL2018a}.
\begin{defn}\label{def-tnn}
A subgroup $H \subset G$ is said to be \emph{totally not normal} if for each non-trivial $h \in H$, there exists $g \in G$ such that $g \ h \ g^{-1} \not\in H$.
\end{defn}

 \subsection{Induced actions}\label{subsec-induced}
 We next consider the relations between the finite model constructed in Section~\ref{subsec-models} and the actions of $G$ and $\oG$ on $\fX$. We first establish some technical  results about restriction maps for adapted sets, then  consider the behavior of the adjoint map for $\whh \in \oG_x$ under restriction to adapted sets for the action. Let $U$ be adapted to the Cantor action  $(\fX,G,\Phi)$.

The first observation is that each $\whg \in \oG$ is the uniform limit of a sequence of maps, $\whg = (g_i)$, and the translates of $U$ form a finite partition of $\fX$ by closed sets, hence   the action of $\whg$ on the \emph{partition} $X_U$ equals the action by $g_i$ for $i$ sufficiently large. Thus, there is a well-defined epimorphism $\whq_U \colon \oG \to Q_U$, which maps $\oG_x$ onto $D_U$ for $x \in U$.

  Let   $(U, G_U, \Phi_U)$ denote the action of $G_U$ restricted to $U$, which is again a Cantor action. 
Define
 \begin{equation}\label{eq-HUV}
H_U \equiv \Phi_U(G_U) \subset \Homeo(U) \quad , \quad K_U = \ker \ \{\Phi_U \colon G_U \to H_U\} \subset G_U \subset G ~ .
\end{equation}
\begin{lemma}\label{lem-Utransitive}
The action of $H_U$ on $U$ is minimal.
\end{lemma}
\proof
As $U$ is adapted,   for each $x \in U$, the translates $H_U \cdot x = \{g \cdot x \in U \mid g \in G_U\}$ are dense.
\endproof

The group $H_U$ is called the \emph{holonomy}   of the action $\Phi$ restricted to $U$, in analogy with the case of actions which arise from the holonomy of weak solenoids \cite{DHL2016b,DHL2016c}. For each $g \in K_U$ the action $\Phi(g)$ fixes every point in $U$. Note that $K_{\fX}$ is the trivial group exactly when $\Phi$ is a faithful action, and $K_U$ can be non-trivial even when the action $\Phi$ on $\fX$ is faithful.

Next,    consider the closures of the restricted groups defined in \eqref{eq-HUV}.  There is a subtle point to consider, that the closure can either be taken in the uniform topology on $\Homeo(\fX)$, 
 or after restriction,  in the uniform topology on $\Homeo(U)$. 
 Let     $\whH_U = \overline{\Phi_U(G_U)} \subset \Homeo(U)$ denote the closure of $H_U$ in the uniform topology on $\Homeo(U)$.
Then   Lemma~\ref{lem-Utransitive} implies:
 \begin{cor}\label{cor-transitive}
  Let $U \subset \fX$ be adapted, then $\whH_U$ acts transitively on $U$.
 \end{cor}

The following is a key technical observation, and is often used implicitly in the following.
   \begin{lemma}\label{lem-closures}
 Let $U \subset \fX$ be adapted for the Cantor action $(\fX,G,\Phi)$. Then
 \begin{equation}
\whH_U = \{ \whg |U  \  \mid \whg \in \oG ~ , ~   \whg \cdot U = U   \} ~ .
\end{equation}
 \end{lemma}
 \proof
 Let $\whh = (g_i) \in \whH_U$ for $g_i \in G_U$, where the sequence $\{\Phi_U(g_i) \}$ converges in the uniform topology in $\Homeo(U)$.
Since $\oG$ is compact, there exists a subsequence $\{g_{i_{\ell}} \mid \ell \geq 0\}$ such that $\{\Phi(g_{i_{\ell}})\}$ converges uniformly in $\Homeo(\fX)$, 
hence yields $\whg = (g_{i_{\ell}}) \in \oG$. Moreover, as $U$ is adapted, we can assume that each $\Phi(g_{i_{\ell}})(U) = U$, 
so   $\whg \cdot U = U$. Then  $\whg |U  = \whh$, as uniform convergence on $\fX$ implies uniform convergence on $U$, hence the actions agree on $U$. This shows that $\whH_U \subset \{\whg | U   \mid \whg \cdot U = U   \}$. The converse is immediate, since $\whg \in \oG$ with $  \whg \cdot U = U $ implies, as noted above, that we can assume $\whg =(g_i)$ where each $\Phi(g_i)(U) = U$, hence $\whg | U   \in \whH_U$.
 \endproof

 For adapted sets $V \subset U$, introduce the   subgroups
  \begin{equation}\label{eq-KUV}
\whH_{U,V} \equiv \{ \whh \in \whH_U \mid \whh \cdot V = V \}  \subset \Homeo(U)  \quad , \quad \whK_{U,V} \equiv \{\whh \in  \whH_{U,V} \mid \whh | V = id\} \subset \Homeo(U) ~ .
\end{equation}
Note that $\whH_{\fX,U} \subset \Homeo(\fX)$ while $\whH_U \subset \Homeo(U)$, and Lemma~\ref{lem-closures} implies there is a surjective restriction map $\rho_U \colon \whH_{\fX,U} \to \whH_U$. 
 
 \begin{lemma}\label{lem-inclusions}
 Let $U \subset \fX$ be adapted   for the Cantor action $(\fX,G,\Phi)$ with $x \in U$. Then
 \begin{equation}\label{eq-inclusions}
  \whK_{\fX,U} \subset \oG_x \subset \whH_{\fX,U}   \subset \oG ~ .
\end{equation}
 \end{lemma}
 \proof
For   $\whh \in \whK_{\fX,U}$ then $\whh \cdot x = x$, hence $\whh \in \oG_x$. 
For $\whh \in \oG_x$ we have $\whh \cdot U \cap U \ne \emptyset$, hence $\whh \cdot U = U$ and so $\whh \in \whH_{\fX,U}$. 
 \endproof
 
  \subsection{Induced adjoint actions}\label{subsec-inducedadj}
For   $U \subset \fX$   adapted  with $x \in U$,  we consider the restriction of the adjoint action of $\oG_x$ to  the subgroup  $\whH_{\fX,U} \subset \oG$.
Let $\whh \in \oG_x$ then $\whh \in \whH_{\fX,U}$ by Lemma~\ref{lem-inclusions}, so the adjoint action of $\whh$   restricts to an action    $\Ad_U(\whh) \colon \whH_{\fX,U} \to \whH_{\fX,U}$. 
Consider the induced   action on the quotient space, $\Ad_U(\whh) \colon U \to U$. 
\begin{lemma}\label{lem-adjointtrivial}
Let $U \in \CO(\fX)$,  $x \in U$. Then   $\whh \in \oG_x$  acts as the identity on $U$ if  and only if 
\begin{equation}\label{eq-adjoint44}
 \Ad_U(\whh)(\whg) \cdot \oG_x  = \whg \cdot \oG_x  ~ , ~ {\rm for ~ all} ~ \whg \in  \whH_{\fX,U} ~ .
  \end{equation}
\end{lemma}
\proof
 Suppose that $\whh$ acts as the identity on $U$, so $\whh \in \oG_x \cap \whK_{\fX,U}$.   
 For $y \in U$, by  Corollary~\ref{cor-transitive} we can choose $\whg \in \oG$ such that $y = \whg \cdot \oG_x$.  
 Then $\whg \cdot U \cap U \ne \emptyset$, hence $\whg \in \whH_{\fX,U}$ and since $\whh \cdot y = y$, we have 
\begin{equation}\label{eq-adjoint33}
\Ad_U(\whh)(\whg) \cdot \oG_x =  \whh \  \whg \ \whh^{-1}   \cdot \oG_x  =  \whh \  \whg \cdot \oG_x = \whh  \cdot y =   y   =   \whg \cdot \oG_x ~ ,
 \end{equation}
 which yields the identity \eqref{eq-adjoint44}.
  Conversely, if  $\whh \in \oG_x$ satisfies \eqref{eq-adjoint44}, then     \eqref{eq-adjoint33} implies that $\whh \cdot y = y$ and that   $\whh \in \oG_x \cap \whK_{\fX, U}$.
 \endproof
 
Let us also introduce the adjoint automorphisms   on the isotropy subgroup $\oG_x$. For  $\whh \in \oG_x$ let $\Ad_x(\whh) \colon \oG_x \to \oG_x$ denote the adjoint action, given by 
$ \Ad_x(\whh)(\whk) = \whh \ \whk \ \whh^{-1}$ for $\whk \in \oG_x$.
 
 \begin{remark}\label{rmk-fiberaction}
 {\rm  Note that   the condition \eqref{eq-adjoint44} does not imply that $\Ad_x(\whh)$ acts as the identity on the isotropy subgroup $\oG_x$ as the action of $\Ad_x(\whh)$ is absorbed when passing to cosets of $\oG_x$.   Examples show that the action of $\Ad_x(\whh)$ can be non-trivial, while the induced quotient   action on $U$ is trivial.  This observation lies behind the distinction between the definitions of the   stabilizer and centralizer limit groups in Section~\ref{subsec-centralizerlimgroup}, and the notions  ``wild of flat type'' and ``dynamically wild'' actions in Definition~\ref{def-wildtypes}.
   }
  \end{remark}

\section{Odometer models}  \label{sec-odometer}
The   odometer model for a Cantor action $(\fX,G,\Phi)$  is obtained from a  sequence of    finite approximations of the action, as in Section~\ref{sec-models}.
For example, a Vietoris solenoid is defined as the inverse limit of a sequence of finite coverings of $\mS^1$, which is equivalent to considering a descending chain of finite index subgroups of the fundamental group $\mZ$ of $\mS^1$. For a Cantor action, one considers a descending chain of finite index subgroups of a fixed group $G$, and then forms the inverse limit action associated to this chain.  The odometer model for an action can be thought of as an ``algebraic model" for the action, which enables   relating the algebraic properties of  $G$ with the dynamics of its action.

Downarowicz gives a survey of classical odometers in \cite{Do2005}. Some properties of odometer models for actions of non-abelian groups $G$     are discussed    by Cortez and Petite in the works \cite{CortezPetite2008,CortezPetite2018}, and by the authors in \cite{DHL2016a,DHL2016c,HL2018a}.  Section 2 of  \cite{CortezMedynets2016} discusses     the literature for   non-abelian odometers.

We first describe the construction of the odometer model as  the inverse limit action obtained from finite approximations associated to adapted subsets of $\fX$. We then  consider the equivalence of these odometer models, and finally consider the adjoint action  for these   inverse limit actions.

 \subsection{Group chains}\label{subsec-groupchains}
For a choice of basepoint $x \in \fX$ and   $\e > 0$, Proposition~\ref{prop-CO} implies there exists an adapted clopen set $U \in \CO(\fX)$ with $x \in U$ and $\diam(U) < \e$.   Thus,   given a basepoint $x$, by iterating this process one can always construct the following:
 
\begin{defn}\label{def-adaptednbhds}
Let  $(\fX,G,\Phi)$   be a Cantor  action.
A properly descending chain of clopen sets $\cU = \{U_{\ell} \subset \fX  \mid \ell \geq 0\}$ is said to be an \emph{adapted neighborhood basis} at $x \in \fX$ for the action $\Phi$  if
    $x \in U_{\ell +1} \subset U_{\ell}$ for all $ \ell \geq 1$ with     $\ds \bigcap_{\ell \geq 1}  \ U_{\ell} = \{x\}$, and  each $U_{\ell}$ is adapted to the action $\Phi$.
\end{defn}

Given an adapted neighborhood basis $\cU$ at $x$, for each $\ell \geq 0$ we can repeat the constructions  of Section~\ref{sec-models} for the adapted set $U_{\ell}$. As a matter of notation, after fixing the collection $\cU$, we use the subscript $\ell$ in place of the subscript $U_{\ell}$ when convenient.
For example,      $G_{\ell} \equiv G_{U_{\ell}}$ denotes the stabilizer group of $U_{\ell}$.
Then  we obtain a descending chain of finite index subgroups 
 \begin{equation}\label{eq-groupchain}
 \cG^x_{\cU} = \{G = G_0 \supset G_1 \supset G_2 \supset \cdots \} \ .
\end{equation}
Note that each $G_{\ell}$ has finite index in $G$, and is not assumed to be a normal subgroup.  Also note that while the intersection of the chain $\cU$ is a single point $\{x\}$, the intersection of the stabilizer groups   in  $\cG^x_{\cU}$ need not be the trivial group.
 
Next, set $X_{\ell} = G/G_{\ell}$ with basepoint $x_{\ell} = e G_{\ell} \in X_{\ell}$,  where $e \in G$ is the identity element. Note that  $G$ acts transitively on the left on   $X_{\ell}$.    
The inclusion $G_{\ell +1} \subset G_{\ell}$ induces a natural $G$-invariant quotient map $p_{\ell +1} \colon X_{\ell +1} \to X_{\ell}$.
 Introduce the inverse limit 
 \begin{equation} \label{eq-invlimspace}
X_{\infty} \equiv \varprojlim \ \{p_{\ell +1} \colon X_{\ell +1} \to X_{\ell} \mid \ell > 0\} 
\end{equation}
which is a Cantor space with the Tychonoff topology, and basepoint $x_{\infty} = (x_{\ell})$. The actions of $G$ on the factors $X_{\ell}$ induce    a minimal  equicontinuous action on $X_{\infty}$, denoted by  $\Phi_x \colon G \times X_{\infty} \to X_{\infty}$.

 For each $\ell \geq 0$, we have the ``partition coding map'' $\Theta_{\ell} \colon \fX \to X_{\ell}$ which is $G$-equivariant.  The maps $\{\Theta_{\ell}\}$ are compatible with the    quotient maps in \eqref{eq-invlimspace}, and so induce a  limit map $\Theta_x \colon \fX \to X_{\infty}$. The fact that the diameters of the clopen sets $\{U_{\ell}\}$ tend to zero, implies that $\Theta_x$ is a homeomorphism.  This is proved in detail in   \cite[Appendix~A]{DHL2016a}.
Moreover, $\Theta_x(x) =  x_{\infty} \in X_{\infty}$, the basepoint of the inverse limit \eqref{eq-invlimspace}.  
  Let $X_{\infty}$ have a  metric such  that $G$ acts on $X_{\infty}$ by isometries, then let $\dX$ be the metric on $\fX$ induced by the homeomorphism $\whTheta_x$.
The minimal equicontinuous action $(X_{\infty}, G, \Phi_x)$   is called the \emph{odometer model} centered at $x$ for the action $(\fX,G,\Phi)$.

 \subsection{The Ellis group}\label{subsec-ellis}
 
 We next introduce the group chain model for the closure group $\oG$, as the inverse limit of the   finite models for the action $\Phi$ introduced in Section~\ref{sec-models}. 
 For each $\ell \geq 1$, let $C_{\ell}$ denote the largest normal subgroup (the \emph{core}) of the stabilizer group $G_{\ell}$, so
\begin{equation}\label{eq-core}
C_{\ell} ~ = ~ \bigcap_{g \in G} ~ g \ G_{\ell} \ g^{-1} ~ \subset G_{\ell} ~ .
\end{equation}
As   $G_{\ell}$ has finite index in $G$, the same holds for $C_{\ell}$. Observe that for all $\ell \geq 0$,  we have  $C_{\ell +1} \subset C_{\ell}$.
Introduce the quotient group  $Q_{\ell} = G/C_{\ell}$ with identity element $e_{\ell} = e C_{\ell} \in Q_{\ell}$. Then $Q_{\ell}$ acts transitively on the quotient set  $X_{\ell} = G/G_{\ell}$. Let  $\pi_{\ell +1} \colon Q_{\ell +1} \to Q_{\ell}$ be the quotient map induced by the inclusion $C_{\ell +1} \subset C_{\ell}$, which is equivariant for the left $G$-actions. Form    the inverse limit group, 
 \begin{equation} \label{eq-invlimgroup}
\whG_{\infty} ~ \equiv ~ \varprojlim \ \{\pi_{\ell +1} \colon Q_{\ell +1} \to Q_{\ell} \mid \ell \geq 0\} ~ \subset ~ \prod_{\ell \geq 0} ~ Q_{\ell} ~ . 
\end{equation}
For $\ell \geq 0$, let $\whq_{\ell} \colon \whG_{\infty} \to Q_{\ell}$ be the projection map onto the $\ell$-th factor.  Let  $\whC_{\ell}$ denote its kernel. 

Give $\whG_{\infty}$   the relative topology induced by the product (Tychonoff) topology; that is, a basis for the topology of $\whG_{\infty}$ is given by the preimages of points for the maps  $\whq_{\ell}$ for $\ell \geq 1$. For $\ell \geq 1$, 
\begin{equation}\label{eq-nbhdbasis}
  \whC_{\ell} ~ = ~ \whq_{\ell}^{-1}(e_{\ell}) ~ = ~  \{ (g_i') \in \whG_{\infty}   \mid g_i' C_i = e_i \ {\rm for} \ 0 \leq i \leq \ell  \}  
\end{equation}
is a   clopen neighborhood of    $\whe = (e_{\ell}) \in \whG_{\infty}$.
Then for  any $\whg \in \whG_{\infty}$, the collection 
$\{  \whg \ \whC_{\ell} \mid \ell \geq 0\}$   forms a neighborhood basis at $\whg$.

 For each $\ell \geq 0$, let    $\whG_{\ell} \subset \whG_{\infty}$ denote the  clopen subgroup
 \begin{equation}\label{eq-whG}
 \whG_{\ell} \equiv \bigcup_{g \in G_{\ell}/C_{\ell} } ~ g \ \whC_{\ell}  = \bigcup_{g \in G_{\ell}/C_{\ell} } ~   \whC_{\ell} \ g ~ .
 \end{equation}

For $\ell \geq 0$, let $q_{\ell} \colon G \to Q_{\ell}$ be the quotient map. The sequence of maps $\{q_{\ell} \mid \ell \geq 0\}$   induces a   homomorphism $q_{\infty} \colon G \to \whG_{\infty}$  for which $q_{\infty}(e) = \whe$. 
  Let  $\whTheta \colon    \oG \to \whG_{\infty}$ denote the closure of $q_{\infty}$ in the uniform topology on maps. 
We then have the folklore result (see  \cite[Theorem~4.4]{DHL2016a} for a proof):
\begin{thm}\label{thm-quotientspace}
Let    $(\fX,G,\Phi)$ be a Cantor action, and suppose that $\{G_{\ell}\mid \ell \geq 0 \}$ is the group chain associated to an adapted neighborhood basis $\cU$ at $x \in \fX$. 
Then  $\whTheta \colon    \oG \to \whG_{\infty}$  is an isomorphism of topological groups.
\end{thm}

   Recall that for   $\ell \geq 0$,  $D_{\ell} \equiv G_{\ell}/C_{\ell}$ is the isotropy group at the basepoint $x_{\ell} = e_{\ell}   G_{\ell} \in X_{\ell}$  of the action of $Q_{\ell}$ on $X_{\ell}$. Moreover,  by the definition of $\whG_{\infty}$, the quotient maps 
   $\pi_{\ell +1} \colon Q_{\ell +1} \to Q_{\ell}$ in \eqref{eq-invlimgroup} are compatible with this identification.

 \begin{lemma}\label{lem-discrimantcalc}
  The image $\whTheta(\oG_x) \subset \whG_{\infty}$ is the subgroup defined by
  \begin{equation}\label{eq-discriminant}
\cD_x ~ \equiv ~    \varprojlim \ \{\pi_{\ell +1} \colon G_{\ell +1}/C_{\ell +1} \to G_{\ell}/C_{\ell} \mid \ell \geq 0\} \subset \whG_{\infty} ~ .
\end{equation}
\end{lemma}
\proof
For  $h \in \oG_x$, let $h_{\ell} = \whq_{\ell}(\whTheta(h)) \in Q_{\ell}$,   then $\whh = (h_{\ell}) \in \whG_{\infty}$.
The action of $h$ on $X_{\ell}$ fixes the coset $e_{\ell} G_{\ell}$, hence  $h_{\ell} \in G_{\ell}/C_{\ell}$ for each $\ell \geq 1$. Moreover, $\pi_{\ell +1}(h_{\ell+1}) = h_{\ell}$ as the sequence $(h_{\ell}) \in \whG_{\infty}$.   Thus,   $\whh \in \cD_x$ is well-defined in \eqref{eq-discriminant}, and $\whPhi_x(\whh)(x_{\infty}) = x_{\infty}$. It follows that 
  $\whTheta(\oG_x) \subset \cD_x$.

Conversely, given $\whh = (h_{\ell}) \in \cD_x$,  we have    $\Phi(h_{\ell}) \in \Homeo(\fX)$ for $\ell \geq 1$. Then  by   Theorem~\ref{thm-quotientspace}, the   sequence $(\Phi(h_{\ell}))$ of homeomorphisms of $\fX$ converges uniformly to a limit $\ovh \in \oG$. 
Note  that multiplication by $\whq_{\ell}(\whh) \in Q_{\ell}$ fixes the coset $e G_{\ell}$, hence the action of  $ \ovh$ on $X_{\infty}$ leaves the  clopen subset $U_{\ell} \subset X_{\infty}$ invariant for each $\ell \geq 0$. Hence, $\whh \cdot x  = x$ so $\ovh \in \oG_x$. 

Define $\whTheta^{-1} \colon \cD_x \to \oG_x \subset \Homeo(\fX)$ by setting $\whTheta^{-1}(\whh) = \ovh$. By the construction, $\whTheta$ and $\whTheta^{-1}$ are   inverse maps and are each continuous,  thus    $\whTheta \colon \oG_x \to \cD_x$ is a topological isomorphism.
\endproof

The group  $\cD_x$ is  called the \emph{discriminant} of the action, 
in analogy with the interpretation of $\whG_{\infty}$ as  the profinite Galois group associated to the sequence of irregular coverings $Q_{\ell} \to X_{\ell}$. 
 \begin{cor}\label{cor-inducedquotients}
The homeomorphism $\Theta_x \colon \fX \to X_{\infty}$  defined in Section~\ref{subsec-groupchains} agrees with the induced map on quotients,
$\whTheta \colon \oG/\oG_x \to \whG_{\infty}/\cD_x$ .
 \end{cor}

 The abstract formula \eqref{eq-discriminant} is useful for computing the discriminant group, as  in \cite[Sections~6--8]{DHL2016a}. In particular, \cite[Example~8.1]{DHL2016a} shows that the 3-dimensional Heisenberg group admits group chains for which $\cD_x$ is a Cantor group.

  The intersection $K(\cG^x_{\cU}) = \bigcap_{\ell \geq 0} ~ G_{\ell}$ is called the \emph{kernel} of $\cG^x_{\cU}$ and need not be the trivial group.
Let  $NK(\cG) \subset K(\cG)$ denote the largest normal  subgroup of  $K(\cG)$, which is also characterized as the    
  kernel of the homomorphism $q_{\infty} \colon G \to \whG_{\infty}$.

 \subsection{Equivalence of group chains}\label{subsec-chainequiv}  
 The construction of the odometer model for a Cantor action $(\fX,G,\Phi)$ depends upon the choice of a basepoint $x \in \fX$ and an adapted neighborhood basis $\cU$ at $x$, which then yields the group chain $\cG^x_{\cU}$ used for the construction. It is thus important to understand how the   group chain $\cG^x_{\cU}$ depends upon these choices. This leads to introducing notions of equivalence between group chains.
 
 Let     $G$ be a fixed group. A group chain in $G$ is a sequence of nested subgroups
 \begin{equation}
\cG = \{ G = G_0 \supset G_1 \supset G_2 \supset \cdots \}
\end{equation}
where each $G_{\ell +1}$ has finite index in $G_{\ell}$ for $\ell \geq 0$.
The first notion of equivalence between group chains   was   used by Rogers and Tollefson   in their study of equivalence of weak solenoids in \cite{RT1971b}, and corresponds to the standard concept of chain equivalence in terms of interlacing of the chains.
\begin{defn} \label{defn-greq}
Two group chains $\{G_{\ell} \mid \ell \geq 0\}$ and $\{H_{\ell} \mid \ell \geq 0\}$ with $G_0=H_0 = G$ are \emph{equivalent}, if and only if, there is a group chain $\{K_{\ell}\mid \ell \geq 0\}$ in $G$, and infinite subsequences $\{G_{\ell_k}\}_{k \geq 0}$ and $\{H_{j_k}\}_{k \geq 0}$ such that $K_{2k} = G_{\ell_k}$ and $K_{2k+1} = H_{j_k}$ for $k \geq 0$.
\end{defn}
For example, suppose that $\cU$ and $\cU'$ are two choices of adapted neighborhood bases at $x$, then the corresponding   group chains  $\cG^x_{\cU}$ and $\cG^x_{\cU'}$  are equivalent in this sense \cite{FO2002,DHL2016a}.

The second notion of equivalence of group chains in $G$ is a generalization of the above notion, and was   introduced by Fokkink and Oversteegen in \cite{FO2002}, and further developed in the authors' work \cite{DHL2016a}.
\begin{defn} \cite{FO2002}\label{conj-equiv}
Two group chains $\{G_{\ell} \mid \ell \geq 0\}$ and $\{H_{\ell} \mid \ell \geq 0\}$ with $G_0=H_0 = G$  are \emph{conjugate equivalent} if and only if there exists a sequence $(g_{\ell}) \subset G_0$ for which  the compatibility condition   $g_{\ell}G_{\ell} = g_{\ell +1} G_{\ell}$ for all $\ell\geq 0$ is satisfied, and 
such that the group chains $\{g_{\ell} G_{\ell} g_{\ell}^{-1} \mid \ell \geq 0 \}$ and $\{H_{\ell} \mid \ell \geq 0\}$ are equivalent. 
\end{defn}
  
The relation between the equivalences in Definitions~\ref{defn-greq} and \ref{conj-equiv} and their corresponding odometer models is given by the following theorem, which follows from results  in \cite{FO2002}. The following result is stated and proved   in  \cite{DHL2016a}.  
\begin{thm}\label{thm-fundamental}
Let $\{G_{\ell} \mid \ell \geq 0 \}$ and $\{H_{\ell}\mid \ell \geq 0\}$ be group chains   $G_0 = H_0 = G$, and let 
\begin{equation*}
X_\infty   =   \varprojlim    \ \{G/G_{\ell+1} \to G/G_{\ell} \mid \ell \geq 0 \}  \quad , \quad 
Y_\infty   =   \varprojlim  \ \{G/H_{\ell+1} \to G/H_{\ell} \mid \ell \geq 0\}  \ .
\end{equation*}
Then the group chains $\{G_{\ell} \mid \ell \geq 0 \}$ and $\{H_{\ell} \mid \ell \geq 0\}$ are \emph{equivalent} if and only if there exists a homeomorphism $\tau \colon  X_\infty \to Y_\infty$ equivariant with respect to the $G$-actions on $X_\infty$ and $Y_\infty$, and such that $\tau(eG_{\ell}) = (eH_{\ell})$. That is, the conjugating map is basepoint-preserving.

  The group chains $\{G_{\ell} \mid  \ell \geq 0 \}$ and $\{H_{\ell} \mid  \ell \geq 0\}$ are \emph{conjugate equivalent} if and only if there exists a homeomorphism $\tau \colon  X_\infty \to Y_\infty$ equivariant with respect to the $G$-actions on $X_\infty$ and $Y_\infty$.
\end{thm}

 \begin{remark}\label{rmk-models}
 {\rm 
 Given  a Cantor action $(\fX,G,\Phi)$, let  $\cG(\Phi)$ denote the   group chains in $G$ which are \emph{conjugate equivalent} to the  group chain $\cG_{\cU}$ for some choice of   basepoint $x \in \fX$ and adapted neighborhood basis $\cU$ at $x$.
    Theorem~\ref{thm-fundamental} implies that an invariant defined for  group chains in $G$, which is independent of the choice of a representative in $\cG(\Phi)$, is an invariant for the conjugacy class of the action $(\fX,G,\Phi)$.   
 }
 \end{remark}

 \begin{remark}\label{rmk-kernels}
 {\rm 
  Let $\cG = \{ G = G_0 \supset G_1 \supset G_2 \supset \cdots \}$ be a group chain, with associated odometer $(X_{\infty}, G, \Phi_x)$.
Let $NK(\cG) \subset K(\cG)$ denote the largest normal  subgroup of the kernel subgroup. It is immediate that $K(\cG)$ is   an invariant of the equivalence class of $\cG$, but examples show that it need not be invariant under conjugate equivalence. In contrast, it is immediate that $NK(\cG)$ is   an invariant of conjugate equivalence, and is identified with the kernel of the action map $\Phi_x$ on $X_{\infty}$. Thus, $NK(\cG)$ is the trivial group whenever the action is faithful.

Applying these remarks to the group chains in $\cG(\Phi)$ associated  to a Cantor action $(\fX,G,\Phi)$, it follows that for a fixed basepoint $x \in \fX$, there is a well-defined kernel subgroup  $K(\cG_{\cU}^x) \subset G$. On the other hand, given another choice of basepoint $y \in \fX$, its kernel $K(\cG_{\cU}^y)$ need not equal $K(\cG_{\cU}^x)$. In fact, it is possible for   kernel at $x$ to be the trivial group, and that at $y$ be a non-trivial subgroup. However, for free actions, they all agree and are trivial.
  }
 \end{remark}
  \begin{lemma}
A  Cantor action $(\fX,G,\Phi)$ is   free   if and only if, for all $x \in \fX$ and  any adapted neighborhood basis $\cU$ at $x$, the kernel $K(\cG^x_{\cU})$  is the trivial group. 
\end{lemma}

\subsection{Adjoint group chains}\label{subsec-adjoints}
We next consider the adjoint actions for  the discriminant group of an odometer model for a Cantor action $(\fX,G,\Phi)$.
Assume  that we have fixed  an   adapted neighborhood basis   $\cU = \{U_{\ell} \subset \fX  \mid \ell \geq 0\}$  at   $x \in \fX$, and     $\cG^x_{\cU} = \{G_{\ell} \mid \ell \geq 0 \}$ is the associated group chain as in \eqref{eq-groupchain}.

 Recall that there are homeomorphisms $\fX \cong \oG/\oG_x$ and $X_{\infty} \cong \whG_{\infty}/\cD_x$. 
  By Theorem~\ref{thm-quotientspace} there is an isomorphism  $\whTheta \colon    \oG \to \whG_{\infty}$    of topological groups, which maps $\oG_x$ isomorphically to $\cD_x$.
For each $\ell \geq 0$, the clopen subgroup $\whG_{\ell} \subset \whG_{\infty}$ satisfies $\cD_x \subset \whG_{\ell}$, and the quotient $\whU_{\ell} \equiv \whG_{\ell}/\cD_x \subset \whG_{\infty}/\cD_x \cong X_{\infty}$. Then Corollary~\ref{cor-inducedquotients} has a sharper statement, that for each $U_{\ell} \subset \fX$ in the adapted neighborhood basis $\cU$ about $x$, we   have    $\Theta_x \colon U_{\ell} \to  \whU_{\ell}$ is a homeomorphism.

 We introduce the localizations of the   adjoint actions of $\cD_x \subset \whG_{\infty}$ as defined    in \eqref{eq-discriminant}, in terms of the odometer model for the action.  We first give a basic fact, where recall that $\whG_{\ell}$ was defined in \eqref{eq-whG}. 
\begin{lemma}\label{lem-intersections}
$\ds \cD_x = \bigcap_{\ell > 0} ~ \whG_{\ell}$ \ .
\end{lemma}
\proof
For $\whh \in \cD_x$,   by \eqref{eq-discriminant} we have $\whh = (h_i)$ where $h_i \in G_i$ for all $i \geq 1$ and $h_{i+1}^{-1} \ h_i \in G_i$.
Since $G_{i+1} \subset G_i$ for all $i \geq 1$, given any $\ell \geq 1$ this implies that we can choose the sequence $(h_i)$ so that $h_i = h_{\ell} \in G_{\ell}$ for all $i \leq \ell$. 
Thus, if $\whh \in \cD_x$ then  $\whh \in \whG_{\ell}$,  for all $\ell \geq 1$, and so $\cD_x \subset \cap_{\ell > 0} \ \whG_{\ell}$.

Conversely, suppose that $\whh \in \whG_{\ell}$   for all $\ell \geq 0$, and hence $\whh \cdot \whG_{\ell} = \whG_{\ell}$.  
Let $\whU_{\ell} \equiv \whG_{\ell}/\cD_x \subset \whG_{\infty}/\cD_x$. 
Then we have $\Theta_x(U_{\ell}) = \whU_{\ell}$, so that $\whU_{\ell}$ 
  is identified with  the clopen neighborhood $U_{\ell} \subset \fX$ of $x$.  
  Since $\cap_{\ell > 0} ~ U_{\ell} = \{x\}$, it follows that $\whh \cdot x = x$, and so $\whh \in \cD_x$. Thus,    $\cD_x = \cap_{\ell > 0} ~ \whG_{\ell}$.  
  \endproof

Now recall that in Section~\ref{subsec-induced}, 
for adapted sets $V \subset U$, we defined in \eqref{eq-KUV}   the groups $\whH_{U,V}$ and $\whK_{U,V}$.
Their analogs for the odometer model of the action are the groups
  \begin{eqnarray} 
\whH_{X_{\infty}, \whU_{\ell}} &  \equiv &  \{ \whg \in \whG_{\infty} \mid \ \whPhi_x(\whg)(\whU_{\ell}) = \whU_{\ell} \} = \whG_{\ell}  \label{eq-Hinf}\\
 \whK_{X_{\infty}, \whU_{\ell}} & \equiv & \{\whg \in  \whG_{\ell} \mid \ \whPhi_x(\whg) \ {\rm acts ~as ~ the~ identity ~ on }~  \whU_{\ell} \}    \subset  \cD_x  ~ . \label{eq-Kinfty}
\end{eqnarray}
 For   each $\ell \geq 0$,  define     
\begin{equation}\label{eq-isotropykernel}
K_{\ell} \equiv   \whK_{X_{\infty}, \whU_{\ell}}    ~ .
\end{equation}
Then  $\whU_{\ell +1} \subset \whU_{\ell}$ implies that $K_{\ell} \subset K_{\ell +1}$ for all $\ell \geq 0$. Thus,  we obtain an  increasing   chain  
$\cK(\cU) = \{ K_0 \subset K_1 \subset K_2 \subset \cdots \}$ of subgroups of $\cD_x$. This is called the \emph{stabilizer group chain}.

Next,   Lemma~\ref{lem-intersections} implies that for each $\ell \geq 0$,  the adjoint action of $\whh \in \cD_x$  restricts to a map 
  $\Ad_{\ell}(\whh) \colon \whG_{\ell} \to \whG_{\ell}$. Define 
  \begin{equation}\label{eq-restrictedadjoint2}
Z_{\ell} \equiv {\rm ker}\{\Ad_{\ell} \colon \cD_x \to \Homeo(\whG_{\ell}) \}\subset \cD_x ~ .
\end{equation}
Then   $\whG_{\ell+1} \subset \whG_{\ell}$ implies that $Z_{\ell} \subset Z_{\ell +1}$ for all $\ell \geq 0$.
  Thus,  we obtain another increasing   chain $\cZ(\cU) = \{ Z_0 \subset Z_1 \subset Z_2 \subset \cdots \}$ of subgroups of $\cD_x$.    This is called the \emph{centralizer group chain}, due to the following observation.
 For a group $H$ and subgroup $K \subset H$, define the \emph{centralizer} of $K$ in $H$,  
\begin{equation}\label{eq-centralizerdef}
\Cen(H,K) = \{h \in H \mid h \ k = k \ h ~ , ~ {\rm for ~ all}~ k \in K\} \subset H ~ .
\end{equation}
Recall that $q_{\infty} \colon G \to \whG_{\infty}$ has dense image.
\begin{lemma} \label{lem-centralizer}
$Z_{\ell} =   \cD_x  \cap \Cen(\whG_{\infty} ,  \whG_{\ell}) =   \cD_x  \cap \Cen(\whG_{\infty} , q_{\infty}(G_{\ell}))$. 
\end{lemma}
\proof
The equality $Z_{\ell} =    \cD_x \cap \Cen(\whG_{\infty} ,  \whG_{\ell})$ is just a restatement of \eqref{eq-restrictedadjoint2}. 
The second equality follows from the observation that if $\whh \in \whG_{\infty}$ commutes with the image group $q_{\infty}(G_{\ell})$, then it also commutes with every element in its closure $\whG_{\ell}$, hence 
   $\Cen(\whG_{\infty} , q_{\infty}(G_{\ell}))\subset \Cen(\whG_{\infty} ,  \whG_{\ell})$. 
\endproof

  The two group chains $\cZ(\cU)$ and $\cK(\cU)$ are closely related.    
\begin{lemma} \label{lem-centralizer2}
There is an inclusion of group chains $\iota_{\infty} \colon \cZ(\cU) \to \cK(\cU)$. 
\end{lemma}
\proof

By the definition   \eqref{eq-restrictedadjoint2}, for each $\ell \geq 0$ and 
   $\whh \in Z_{\ell}$, the action $\Ad_{\ell}(\whh)$ on $\whG_{\ell}$ is trivial.
Hence by   Lemma~\ref{lem-adjointtrivial}, the induced action   $\whPhi_x(\whh)$ on $\whU_{\ell} \equiv \whG_{\ell}/\cD_x$ is trivial. 
Thus, there is an inclusion map    $\iota_{\ell} \colon Z_{\ell} \subset    K_{\ell}$, which yields the inclusion map   $\iota_{\infty} = (\iota_{\ell})$   on group chains.  
\endproof

Consider $\whh \in Z_{\ell+1}$ such that $\whh \not\in Z_{\ell}$. Then $\whPhi_x(\whh)$ acts trivially on $\whU_{\ell +1}$, and $\Ad_{\ell}(\whh)$ acts non-trivially on $\whG_{\ell}$. However,   it is possible that $\whPhi_x(\whh)$  also acts trivially on $\whU_{\ell}$, while  the adjoint action $\Ad_{\ell}(\whh)$ acts   non-trivially   on the subgroup $\cD_x$ as discussed in Remark~\ref{rmk-fiberaction}.  The point is that comparing the subgroups $Z_{\ell} \subset K_{\ell}$   involves   algebraic properties of the inclusion $\cD_x \subset \whG_{\infty}$.

 Finally,  given  $\whh  \in \cD_x$ which is represented by a sequence   $\whh = (h_{\ell})$ with $h_{\ell} = h \in G$ for all $\ell \geq 0$, then  $h \in \whG_{\ell}$ of all $\ell \geq 0$,  so $h \in K(\cG^x_{\cU})$. 
Conversely, recall that $NK(\cG^x_{\cU})$ is the largest normal subgroup of $K(\cG^x_{\cU})$, and so $NK(\cG^x_{\cU}) \subset C_{\ell}$ for all $\ell \geq 1$, hence each $h \in NK(\cG^x_{\cU})$   acts trivially on $\whG_{\infty}$. Define  the quotient group 
\begin{equation}\label{eq-ratdisc}
\cD_x^{G} \equiv K(\cG^x_{\cU})/NK(\cG^x_{\cU}) \subset \cD_x ~ ,
\end{equation}
then  each $h \in \cD_x^{G}$   is represented by a constant sequence $(h)$ for $h \in G$. One can think of the elements of $\cD_x^{G}$ as the ``$G$-rational points'' of $\cD_x$. Of course, if $K(\cG^x_{\cU})$ is trivial, then there are no such points. 
The examples in Section~\ref{sec-arboreal} illustrate actions with non-trivial $G$-rational points.

\section{Direct limit conjugacy  invariants}\label{sec-limgroups}

In this section, we use the direct limit construction for the group chains $\cZ(\cU)$ and  $\cK(\cU)$ to obtain conjugacy invariants of a Cantor action $(\fX,G,\Phi)$. These are subtle invariants of an action, and each is independent of the other, 
as shown by Theorem~\ref{thm-distinct=nH} and the examples in Section~\ref{sec-arboreal}.
We first briefly recall   the basic formulation of direct limits, and then apply these ideas to obtain   the \emph{stabilizer limit  group} $\Upsilon_s(\Phi)$ and the \emph{centralizer  limit  group} $\Upsilon_c(\Phi)$ for a   Cantor action $(\fX,G,\Phi)$.    
  
\subsection{Direct limits}
We   give the construction and properties of the direct limit in the category of groups. Basic references for  this standard concept are     Eilenberg and Steenrod \cite[Chapter VIII, Section 2]{EilenbergSteenrod1952}, 
and the  text by Munkres \cite[Section 73]{Munkres1984}, which give proofs of the following results.

  \begin{defn}\label{def-directedgroup}
  A directed system of groups $\cG(G_{\lambda}, \phi^{\lambda'}_{\lambda} , \Lambda)$  
  over a directed set $\Lambda$ is collection of groups   $\{G_{\lambda} \mid \lambda \in \Lambda\}$, and for each $\lambda, \lambda' \in \Lambda$ with $\lambda < \lambda'$,  a group homomorphism $\phi^{\lambda'}_{\lambda}  \colon G_{\lambda} \to G_{\lambda'} $. We require that $\phi^{\lambda}_{\lambda} = id_{G_{\lambda}}$, 
  and for   $\lambda < \lambda' < \lambda''$ 
  that $ \phi^{\lambda''}_{\lambda} = \phi^{\lambda''}_{\lambda'} \circ \phi^{\lambda'}_{\lambda}$.
\end{defn}

 \begin{defn}\label{def-directedgroupequivalence}
Let  $\cG(G_{\lambda}, \phi^{\lambda'}_{\lambda} , \Lambda)$  
 be a  directed system of groups  over a directed set $\Lambda$.  Define an equivalence relation on the disjoint union
$ \bigcup_{\lambda \in \Lambda}  \ G_{\lambda}$ where for  $g_{\lambda} \in G_{\lambda}$ and     $g_{\lambda'} \in G_{\lambda'}$, we set $g_{\lambda} \sim g_{\lambda'}$ if for all $\lambda'' \in \Lambda$ with $ \lambda < \lambda''$ and $\lambda' < \lambda''$, then 
$\phi^{\lambda''}_{\lambda}(g_{\lambda}) = \phi^{\lambda''}_{\lambda'}(g_{\lambda'})$.
\end{defn}

   \begin{defn}\label{def-directlimit}
The \emph{direct limit} of a   directed system of groups $\cG(G_{\lambda}, \phi^{\lambda'}_{\lambda} , \Lambda)$   over a directed set $\Lambda$ is defined to be the set of equivalence classes 
\begin{equation}\label{eq-dirlim}
\varinjlim \cG(G_{\lambda}, \phi^{\lambda'}_{\lambda} , \Lambda) ~   \equiv ~ \bigcup_{\lambda \in \Lambda}  \ G_{\lambda}  \bigg/  g_{\lambda} \sim g_{\lambda'} \ .
\end{equation}
The   group structure on $\ds \varinjlim \cG(G_{\lambda}, \phi^{\lambda'}_{\lambda} , \Lambda)$ is inherited from that on the summands.
For  $\lambda \in \Lambda$ and $g_{\lambda} \in G_{\lambda}$ let  $[g_{\lambda}] \in \varinjlim \cG(G_{\lambda}, \phi^{\lambda'}_{\lambda} , \Lambda)$ denote the equivalence class it determines.
\end{defn}

    \begin{defn}\label{def-directedgroupmorphism}
  A map $\Xi$ between directed systems of groups 
  $\cG(G_{\lambda}, \phi^{\lambda'}_{\lambda} , \Lambda)$  
  and $\cG(H_{\omega}, \psi^{\omega'}_{\omega} , \Omega)$
 is an order-preserving map $\xi \colon \Lambda \to \Omega$, and for each $\lambda \in \Lambda$ 
 a group homomorphism $\xi_{\lambda} \colon G_{\lambda} \to H_{\xi(\lambda)}$ such that for $\lambda < \lambda'$ we have
 $$   \xi_{\lambda'} \circ \phi^{\lambda'}_{\lambda} = \phi^{\xi(\lambda')}_{\xi(\lambda)} \circ \xi_{\lambda} \colon G_{\lambda} \to H_{\xi(\lambda')} ~ .$$
\end{defn}

\begin{prop}
 A map $\Xi$ between directed systems of groups 
  $\cG(G_{\lambda}, \phi^{\lambda'}_{\lambda} , \Lambda)$  
  and $\cG(H_{\omega}, \psi^{\omega'}_{\omega} , \Omega)$ induces a homomorphism 
$\ds \underrightarrow{\Xi} \colon \varinjlim \cG(G_{\lambda}, \phi^{\lambda'}_{\lambda} , \Lambda) \longrightarrow \varinjlim \cG(H_{\omega}, \psi^{\omega'}_{\omega} , \Omega)$. 
\end{prop}

We next recall some special cases of directed systems of groups  and maps between directed systems of groups.
We single out two special classes of maps.
   \begin{prop}\label{prop-directedgroupmonomorphism}
  A map $\Xi$ between directed systems of groups 
  $\cG(G_{\lambda}, \phi^{\lambda'}_{\lambda} , \Lambda)$  
  and $\cG(H_{\omega}, \psi^{\omega'}_{\omega} , \Omega)$
 is  a  \emph{monomorphism} if each of the maps  $\xi_{\lambda} \colon G_{\lambda} \to H_{\xi(\lambda)}$ for $\lambda \in \Lambda$ is a  monomorphism of groups. Then the induced map 
 $\ds \underrightarrow{\Xi} \colon \varinjlim \cG(G_{\lambda}, \phi^{\lambda'}_{\lambda} , \Lambda) \longrightarrow \varinjlim \cG(H_{\omega}, \psi^{\omega'}_{\omega} , \Omega)$
 is a group monomorphism.
\end{prop}

\smallskip
   \begin{prop}\label{prop-directedgroupisomorphism}
  A map $\Xi$ between directed systems of groups 
  $\cG(G_{\lambda}, \phi^{\lambda'}_{\lambda} , \Lambda)$  
  and $\cG(H_{\omega}, \psi^{\omega'}_{\omega} , \Omega)$
 is  an \emph{isomorphism} if each of the maps  $\xi_{\lambda} \colon G_{\lambda} \to H_{\xi(\lambda)}$ for $\lambda \in \Lambda$ is an isomorphism of groups.  
  Then the induced map 
 $\ds \underrightarrow{\Xi} \colon \varinjlim \cG(G_{\lambda}, \phi^{\lambda'}_{\lambda} , \Lambda) \longrightarrow \varinjlim \cG(H_{\omega}, \psi^{\omega'}_{\omega} , \Omega)$
 is a group isomorphism.
\end{prop}

A subset $\Lambda_0 \subset \Lambda$ of a directed set is said to be \emph{cofinal} if for each $\lambda \in \Lambda$, there exists $\lambda' \in \Lambda_0$ with $\lambda < \lambda'$. Then we have
\begin{prop}\label{prop-directedgroupcofinal}
Let  $\cG(G_{\lambda}, \phi^{\lambda'}_{\lambda} , \Lambda)$  be a directed systems of groups, and $\Lambda_0 \subset \Lambda$ be a cofinal set. Then the inclusion induces a group isomorphism
$\ds  \varinjlim \cG(G_{\lambda}, \phi^{\lambda'}_{\lambda} , \Lambda_0) \cong   \varinjlim \cG(G_{\lambda}, \phi^{\lambda'}_{\lambda} , \Lambda)$.
\end{prop}

In addition, we introduce two special classes of direct limit systems.

  \begin{defn}\label{def-boundeddirectedgroup}
  A   directed system of groups 
  $\cG(G_{\lambda}, \phi^{\lambda'}_{\lambda} , \Lambda)$   is  said to be   \emph{bounded} if there exists $\lambda_0 \in \Lambda$ such that 
  for all $\lambda_0 < \lambda < \lambda'$ the map $\phi^{\lambda'}_{\lambda}$ is a group isomorphism.  The directed system  is said to be \emph{unbounded} if no such $\lambda_0$ exists.
  \end{defn}

   \begin{defn}\label{def-finitetype}
  A   directed system  
  $\cG(G_{\lambda}, \phi^{\lambda'}_{\lambda} , \Lambda)$   has   \emph{finite type} if each group $G_{\lambda}$ is finite.  
  \end{defn}

 For our purposes, the ordered sets $\Lambda$ and $\Omega$ above will always be assumed to be a subset of the non-negative integers $\mN$ with the natural order.  All maps $\{ \phi^{\lambda'}_{\lambda} \mid \lambda < \lambda' \in \Lambda\}$ in our applications will  be group inclusions,    possibly isomorphisms. For a  morphism $\Xi$ of directed systems of groups, the maps $\{\xi(\lambda) \mid \lambda \in \Lambda\}$ will be be group inclusions (possibly isomorphisms).  
 
 Now assume that $\Lambda \subset \mN$, then     $\ds \varinjlim \cG(G_{\lambda}, \phi^{\lambda'}_{\lambda} , \Lambda)$   inherits a filtration from its definition.

\begin{defn}\label{defn-height}
The \emph{height function} $\lambda \colon \varinjlim \cG(G_{\lambda}, \phi^{\lambda'}_{\lambda} , \Lambda) \to \Lambda$ is defined  for an equivalence class $[g]$ to be    the least $\lambda$ such that $[g] = [g_{\lambda}]$ for some $g_{\lambda} \in G_{\lambda}$. 
\end{defn}

\begin{remark}\label{rmk-order}
{\rm 
 Note that while the height function need not be preserved by a map between directed systems, it does yield a well-defined order, where $[g_1] \leq_{\lambda} [g_2]$ if  $\lambda(g_1) < \lambda(g_2)$.
}
\end{remark}

\begin{remark}\label{rmk-unions}
{\rm 
 Suppose there exists a group $\fG$ such that $G_{\lambda} \subset \fG$ for all $\lambda \in \Lambda$, and the maps $\phi^{\lambda'}_{\lambda}$ are inclusions of subgroups.   Define 
 \begin{equation}\label{eq-limitobject}
G_{\infty} \equiv  \bigcup_{\lambda \in \Lambda} \ G_{\lambda} \ .
\end{equation}
For $g_{\lambda} \in G_{\lambda}$ and $g'_{\lambda'} \in G_{\lambda'}$ then there exists $\lambda'' > \max\{\lambda, \lambda'\}$ so that we can consider $g_{\lambda}, g_{\lambda'} \in G_{\lambda''}$ and then $g_{\lambda} + g'_{\lambda'}$ is defined in $G_{\lambda''}\subset G_{\infty}$. That is, $G_{\infty}$ is a subgroup of $\fG$.

Define a height function   $\lambda_{\infty} \colon G_{\infty} \to \Lambda$ by $\lambda_{\infty}(g) \leq \lambda$ if $g \in \ G_{\lambda}$. Given  the pair $\{G_{\infty}, \lambda_{\infty}\}$, define  $\wtG_{\lambda} \equiv \{g \in G_{\infty} \mid \lambda_{\infty}(g) \leq \lambda\}$, and introduce maps given by inclusions.  Let $\cG(G_{\infty}, \lambda_{\infty})$ denote the resulting direct limit of groups, which is a small abuse of notation. 
 It is then immediate that there is an isomorphism of direct limits $  \cG(G_{\lambda}, \phi^{\lambda'}_{\lambda} , \Lambda) \cong \cG(G_{\infty}, \lambda_{\infty})$. 
    Thus, one can view a direct limit of subgroups of $\fG$ as a subgroup $G_{\infty} \subset \fG$, equipped with a height function $\lambda_{\infty} \colon G_{\infty} \to \Lambda$.
}
\end{remark}

In the case when the directed system of groups 
  $\cG(G_{\lambda}, \phi^{\lambda'}_{\lambda} , \Lambda)$   is     \emph{bounded}, with $\lambda_0$ an upper bound, then the inclusion of the constant sequence 
  $\cG(G_{\lambda_0}, Id_{G_{\lambda_0}} , \{\lambda_0\}) \subset \cG(G_{\lambda}, \phi^{\lambda'}_{\lambda} , \Lambda)$ 
  is an isomorphism, and  we have $G_{\infty} =    G_{\lambda_0}$. In other words, the direct limit of a bounded direct limit collapses.

    \subsection{The   stabilizer and centralizer limit groups}\label{subsec-centralizerlimgroup}
    We now show that  the chain of stabilizer subgroups defined by \eqref{eq-isotropykernel}, and the chain of centralizer subgroups defined in Lemma~\ref{lem-centralizer}, yield direct limit groups whose isomorphism classes     are conjugacy invariants of Cantor actions.

 \begin{defn}\label{def-progroupinvariants}
 Let    $(\fX,G,\Phi)$ be a Cantor action, and $\cG^x_{\cU} = \{G_{\ell}\mid \ell \geq 0 \}$ be the group chain associated to an adapted neighborhood basis $\cU$ at $x \in \fX$.  Let $\{K_{\ell}  \mid \ell \geq 0\}$ be the chain of stabilizer groups associated to $\cG^x_{\cU}$, 
 and let  $\{Z_{\ell}  \mid \ell \geq 0\}$  be the chain of centralizer groups associated to $\cG^x_{\cU}$. Then define:
 \begin{enumerate}
\item The \emph{stabilizer (direct limit) group} is   \quad  $\ds \Upsilon_s^x(\Phi) = \varinjlim \cG(K_{\ell}, \iota^{\ell'}_{\ell} , \mN)$;
\item The \emph{centralizer (direct limit) group} is  \quad  $\ds \Upsilon_c^x(\Phi) = \varinjlim \cG(Z_{\ell}, \iota^{\ell'}_{\ell} , \mN)$; 
\end{enumerate}
where $\iota^{\ell'}_{\ell}$  denotes the inclusion map,  for $0 \leq \ell \leq \ell'$.
 \end{defn}

   We now come to the main result of this section.
 \begin{thm}\label{thm-centraldirectlimgroup}
 Let    $(\fX,G,\Phi)$ be a Cantor action. The direct limit isomorphism classes    of     $\Upsilon_s^x(\Phi)$ and     $\Upsilon_c^x(\Phi)$  are   invariants of the conjugacy class of the action $\Phi$. These isomorphism classes are denoted by   $\Upsilon_s(\Phi)$ and     $\Upsilon_c(\Phi)$, respectively, where we suppress the  basepoint $x$ in the notation.
 \end{thm}
 \proof 
Let  $(\fX,G,\Phi)$ and  $(\fX',G,\Psi)$ be conjugate Cantor actions by a homeomorphism $h \colon \fX' \to \fX$. Let  $\cU = \{U_{\ell} \subset \fX  \mid \ell \geq 0\}$ be an adapted neighborhood basis at $x \in \fX$ for the action $\Phi$, and let  $\cV = \{V_{\ell} \subset \fX' \mid \ell \geq 0\}$ be an adapted neighborhood basis at $z \in \fX'$ for the action $\Psi$. Then $\cU' = \{ U'_{\ell} = h(V_{\ell}) \subset \fX \mid \ell \geq 0\}$ is  an adapted neighborhood basis   at $y = h(z)$ for the action $\Phi$, and the group chain in $G$ associated to $\cV$ and the action $\Psi$ coincides with the group chain in $G$ associated to $\cU'$ and the action $\Phi$. Thus, by Theorem~\ref{thm-fundamental}, we need to show that the direct limits $\Upsilon_s^x(\Phi)$ and $\Upsilon_c^x(\Phi)$, associated  to the group chain $\cG^x_{\cU}$  in $G$, are isomorphic as direct limits to the direct limits $\Upsilon_s^y(\Phi)$ and $\Upsilon_c^y(\Phi)$ associated  to the group chain $\cG^y_{\cU'}$.

 First consider the case where $x=y$, so we are given two  adapted neighborhood bases  at a common basepoint $x$,   $\cU = \{U_{\ell} \subset \fX  \mid \ell \geq 0\}$   and $\cU' = \{U'_{\ell} \subset \fX  \mid \ell \geq 0\}$, with corresponding group chains $\cG^x_{\cU} = \{G_{\ell} \mid \ell \geq 0\}$ and $\cG^x_{\cU'} = \{G'_{\ell}\mid \ell \geq 0\}$.  
 
  As $\cU$ and $\cU'$ are both adapted neighborhood basis at $x$,  there exists   increasing sequences of indices   
 $1 \leq i_1 < i_2 < i _3 < \cdots$ and $1 \leq j_1 < j_2 < j_3 < \cdots$ such that we have a descending sequence of adapted clopen sets at $x$, where $i_0 = j_0 = 0$, 
 $$ \fX = U_0  = U'_0 \supset U_{i_1} \supset U'_{j_1} \supset U_{i_2} \supset U'_{j_2} \supset \cdots \ .$$
 Then by Theorem~\ref{thm-fundamental} and passing to a subsequence, we can   assume without loss of generality that $i_{\ell} = \ell$ and   $j_{\ell} = \ell$, for $\ell \geq 0$.
Introduce a common refinement $\cU'' = \{U''_{\ell} \mid \ell \geq 0\}$ of these chains of clopen sets, where $U''_{2\ell} = U_{\ell}$ and $U''_{2\ell -1} = U'_{\ell}$.

 Let   $\cG^x_{\cU''} = \{G''_{\ell}\mid \ell \geq 0\}$  be the group chain  associated to $\cU''$, then $G''_{2\ell} = G_{\ell}$ and $G''_{2\ell -1} = G'_{\ell}$. Let $X_{\infty}$, $X'_{\infty}$ and $X''_{\infty}$ denote the inverse limit spaces  defined as in \eqref{eq-invlimspace} by the group chains   $\cG^x_{\cU}$,  $\cG^x_{\cU'}$ and  $\cG^x_{\cU''}$, respectively. 
 Let   $\widehat{\cU} = \{\whU_{\ell} \subset X_{\infty}  \mid \ell \geq 0\}$ be the corresponding adapted basis at $X_{\infty}$, and likewise $\widehat{\cU}' = \{\whU'_{\ell} \subset X'_{\infty}  \mid \ell \geq 0\}$ and $\widehat{\cU}'' = \{\whU''_{\ell} \subset X''_{\infty}  \mid \ell \geq 0\}$ for $X'_{\infty}$ and $X''_{\infty}$, respectively.
 
 By  the discussion in Section~\ref{subsec-groupchains},  there are homeomorphisms 
 $\Theta_x \colon \fX \to X_{\infty}$, $\Theta_x' \colon \fX \to X'_{\infty}$ and $\Theta_x'' \colon \fX \to X''_{\infty}$ which intertwine the $G$-actions on these spaces.
Introduce the basepoint preserving homeomorphisms  $\tau = \Theta_x \circ (\Theta_x'')^{-1} \colon X''_{\infty} \to X_{\infty}$ and 
   $\tau' = \Theta_x' \circ (\Theta_x'')^{-1} \colon X''_{\infty} \to X'_{\infty}$. 
 The maps $\tau$ and $\tau'$ have simple descriptions in terms of sequences. 
 Given  $(g''_{i}) \in X''_{\infty}$ then $\tau(g''_{i}) = (g_i) \in  X_{\infty}$ where $g_i = g''_{2i}$, 
 and   $\tau'(g''_{i}) = (g'_i) \in  X'_{\infty}$ where $g'_i = g''_{2i-1}$.

 Let $\whG_{\infty}$, $\whG'_{\infty}$ and $\whG''_{\infty}$ denote the inverse limit groups  defined as in \eqref{eq-invlimgroup}  by the group chains   $\cG^x_{\cU}$,  $\cG^x_{\cU'}$ and  $\cG^x_{\cU''}$, respectively. 
By  Theorem~\ref{thm-quotientspace}, there are topological isomorphisms 
$$\whTheta \colon \oG \to \whG_{\infty} ~ , ~ \whTheta' \colon \oG \to \whG'_{\infty} ~ , ~ \whTheta'' \colon \oG \to \whG''_{\infty} \ .$$ 
There are topological isomorphisms  $\whtau = \whTheta \circ (\whTheta'')^{-1} \colon \whG''_{\infty} \to \whG_{\infty}$ and 
   $\whtau' = \whTheta' \circ (\whTheta'')^{-1} \colon \whG''_{\infty} \to \whG'_{\infty}$. 
 Given a sequence $\whg = (g''_{i}) \in \whG''_{\infty}$ then   $\whtau(g''_{i}) = (g_i = g''_{2i}) \in  \whG_{\infty}$ and 
 $\whtau'(g''_{i}) = (g'_i = g''_{2i-1}) \in  \whG'_{\infty}$.

We first consider the case of the stabilizer group chains. 
Let $\cD_x'' \subset \whG''_{\infty}$ denote the discriminant for the chain $\cG^x_{\cU''}$, then 
   the restriction    $\whtau \colon \cD''_x \to  \cD_x  \subset \whG_{\infty}$  is an isomorphism  by Lemma~\ref{lem-discrimantcalc} and Corollary~\ref{cor-inducedquotients}, and likewise  for $\whtau' \colon \cD''_x \to  \cD'_x \subset \whG'_{\infty}$. 
   
   We show that the stabilizer groups associated to the chains $\cG^x_{\cU''}$ and $\cG^x_{\cU}$ are  isomorphic, with the proof for $\cG^x_{\cU''}$ and $\cG^x_{\cU'}$  being the same.
For each $\ell \geq 0$,   by  \eqref{eq-Kinfty} and \eqref{eq-isotropykernel} we have
\begin{eqnarray*}
K_{\ell} \equiv  \whK_{X_{\infty}, \whU_{\ell}}   & = &   \{\whg \in  \whG_{\ell} \mid \whPhi_x(\whg) \ {\rm acts ~as ~ the~ identity ~ on }~  \whU_{\ell} \} \subset \cD_x.  \\
K''_{\ell} \equiv  \whK_{X''_{\infty}, \whU''_{\ell}}   & = &   \{\whg \in  \whG''_{\ell} \mid \whPhi''_x(\whg) \ {\rm acts ~as ~ the~ identity ~ on }~  \whU''_{\ell} \} \subset \cD''_x\ . 
\end{eqnarray*}
The  isomorphism $\whtau \colon \whG''_{2\ell} \to \whG_{\ell}$ induces a homeomorphism   $\tau_{\ell} \colon \whU''_{2\ell} = \whG''_{2\ell}/\cD''_x  \to \whU_{\ell} = \whG_{\ell}/\cD_x$. 
Thus $\whtau$ induces an isomorphism $\whtau \colon K''_{2\ell} = \whK_{X''_{\infty}, \whU''_{2\ell}} \to \whK_{X_{\infty}, \whU_{\ell}} =  K_{\ell}$.
  
For $0 \leq \ell \leq \ell'$, let $\phi^{\ell'}_{\ell} \colon K_{\ell} \subset K_{\ell'}$ and  $\psi^{\ell'}_{\ell} \colon K''_{\ell} \subset K''_{\ell'}$ denote the inclusion maps.
By Propositions~\ref{prop-directedgroupisomorphism} and \ref{prop-directedgroupcofinal}, we obtain an isomorphism of direct limits 
\begin{equation}\label{eq-isoisotropy}
 \underrightarrow{\whtau} \colon  \varinjlim \cG(K''_{\ell}, \psi^{\ell'}_{\ell} , \mN) \to \varinjlim \cG(K_{\ell}, \phi^{\ell'}_{\ell} , \mN) \ ,
\end{equation}
 as was to be shown.  Hence, the direct limits $\varinjlim \cG(K_{\ell}, \phi^{\ell'}_{\ell} , \mN)$ and $\varinjlim \cG(K'_{\ell}, \psi^{\ell'}_{\ell} , \mN)$  are  isomorphic.

 Next, consider the case of the centralizer group chains, given by
 \begin{equation}\label{eq-centralizerschains}
Z_{\ell} = \Cen(\whG_{\infty} , \iota_{\infty}(G_{\ell})) \cap \cD_x ~ , ~ Z'_{\ell} = \Cen(\whG'_{\infty} , \iota'_{\infty}(G'_{\ell})) \cap \cD'_x ~ , ~ Z''_{\ell} = \Cen(\whG''_{\infty} , \iota_{\infty}(G''_{\ell})) \cap \cD''_x  ~ ,
\end{equation}
which define the group chains $\cZ$, $\cZ'$ and $\cZ''$, respectively.
  
Recall   the topological isomorphisms 
$\whtau \colon \whG''_{\infty} \to \whG_{\infty}$ and $\whtau' \colon \whG''_{\infty} \to \whG'_{\infty}$ are  given by 
    $\whtau(g_{i}) = (g_{2i})$ and  $\whtau'(g_{i}) = (g_{2i-1})$, for   $\whg = (g_{i}) \in \whG''_{\infty}$.
Then $\whtau$ induces   isomorphisms $\whtau_{\ell} \colon Z''_{2\ell} \to Z_{\ell}$ 
and $\whtau'_{\ell} \colon Z''_{2\ell -1} \to Z'_{\ell}$. Thus, as for the case of the stabilizer group chains, we obtain an isomorphism 
\begin{equation}\label{eq-isocentralizer}
 \underrightarrow{\whtau} \colon  \varinjlim \cG(Z''_{\ell}, \psi^{\ell'}_{\ell} , \mN) \to \varinjlim \cG(Z_{\ell}, \phi^{\ell'}_{\ell} , \mN) \ ,
\end{equation}
as was to be shown. Hence, the direct limits $\varinjlim \cG(Z_{\ell}, \phi^{\ell'}_{\ell} , \mN)$ and $\varinjlim \cG(Z'_{\ell}, \phi^{\ell'}_{\ell} , \mN)$  are  isomorphic.

Next,  consider the case where $x$ and $y$ are distinct basepoints, and we are given    adapted neighborhood bases    $\cU = \{U_{\ell} \subset \fX  \mid \ell \geq 0\}$  at $x$  and $\cU' = \{U'_{\ell} \subset \fX  \mid \ell \geq 0\}$ at $y$, with corresponding group chains $\cG^x_{\cU} = \{G_{\ell} \mid \ell \geq 0\}$ and $\cG^y_{\cU'} = \{G'_{\ell}\mid \ell \geq 0\}$.  
Theorem~\ref{thm-fundamental} implies that the chains $\cG^x_{\cU}$ and $\cG^y_{\cU'}$ are conjugate equivalent. That is, 
choose $\whg = (g_i) \in \oG$ such that $\whg \cdot y = x$,  then the collection 
 $\cU'' = \{g_i \cdot U'_{\ell} \subset \fX  \mid \ell \geq 0\}$  is an  adapted neighborhood basis  at $x$, and the associated group chain   
 $\cG^x_{\cU''} = \{G''_{\ell} = g_{\ell}^{-1} \ \whG'_{\ell} \ g_{\ell} \}$ is equivalent to $\cG^x_{\cU}$. 

 A key point is that the isomorphism between $\whG'_{\ell} \subset \whG_{\infty}$ and $\whG''_{\ell} \subset \whG_{\infty}$ is induced by 
  the conjugacy isomorphism $\Ad(g_{\ell}) \colon G_{\ell} \to g_{\ell} G_{\ell} g_{\ell}^{-1}$ which also conjugates the stabilizer groups $K''_{\ell}$ and $K'_{\ell}$, and likewise for the centralizer groups $Z''_{\ell}$ and $Z'_{\ell}$.  Then by Propositions~\ref{prop-directedgroupisomorphism} and \ref{prop-directedgroupcofinal}, we obtain an isomorphism of direct limits.  We are thus reduced to the case shown previously, where we are given direct limit groups obtained from adapted neighborhood systems centered at the same point $x$.
   \endproof
 
 \begin{remark}
 {\rm
The proof of  Theorem~\ref{thm-centraldirectlimgroup} shows that the isomorphism between the direct limits is induced by a group isomorphism 
 $\Ad(\whg) \colon \whG_{\infty} \to \whG_{\infty}$. By  Remark~\ref{rmk-unions}, for $\fG = \whG_{\infty}$,  the adjoint induces the isomorphism of subgroups  $\Ad(\whg) \colon K_{\infty} \cong K'_{\infty}$ and  $\Ad(\whg) \colon Z_{\infty} \cong Z'_{\infty}$. A key point is that this isomorphism preserves the   height filtration on these groups as defined in Remark~\ref{rmk-unions}.
  }
 \end{remark}

 \subsection{Properties of the stabilizer and centralizer  groups} 
The  direct limit groups  $\Upsilon_s(\Phi)$ and     $\Upsilon_c(\Phi)$ are    invariants of Cantor actions whose dynamical implications will be considered in the following sections. We first give  a restatement   of  Lemma~\ref{lem-centralizer2} in the language of direct limits.

 \begin{lemma}\label{lem-inclusionprogroups}
 For $x \in \fX$, there is an inclusion of direct limit groups $\Upsilon_c(\Phi) \subset  \Upsilon_s(\Phi)$.
 \end{lemma}
 Examples show that the inclusion $\Upsilon_c(\Phi) \subset  \Upsilon_s(\Phi)$ can be proper.  One source of this distinction is noted in Remark~\ref{rmk-fiberaction}, that  the adjoint map for $\cD_x$ restricted to $\cD_x$ need not be trivial, while the action on $U_{\ell}$ is trivial.  We next introduce five classes of Cantor actions.
  \begin{defn}\label{def-wildtypes}
A Cantor action     $(\fX,G,\Phi)$ is said to be: 
\begin{enumerate}
\item \emph{stable} if    the stabilizer  group $\Upsilon_s(\Phi)$   is bounded, and  is said to be \emph{wild} otherwise;
\item  \emph{algebraically stable} if the  its centralizer  group $\Upsilon_c(\Phi)$   is bounded, and  is said to be \emph{algebraically wild} otherwise;
\item \emph{wild of finite type}  if    the stabilizer  group $\Upsilon_s(\Phi)$   is unbounded, and represented by a chain of finite groups;
\item \emph{wild of flat type}  if    the stabilizer  group $\Upsilon_s(\Phi)$   is unbounded, and $\Upsilon_c(\Phi) = \Upsilon_s(\Phi)$; 
\item \emph{dynamically wild}  if    the stabilizer  group $\Upsilon_s(\Phi)$   is unbounded, and is not of flat type.
\end{enumerate}
   \end{defn}

 Each of the properties in Definition~\ref{def-wildtypes}   is a \emph{conjugacy} invariant of the action by    Theorem~\ref{thm-centraldirectlimgroup}.
   
   The notion of wild   Cantor actions  was   introduced by the authors in    \cite[Definition~4.6]{HL2018a}, as part of the study of the homeomorphism types of   weak solenoids in their works   \cite{DHL2016c,HL2018a,HL2018b}.  The definition of a wild action in \cite{HL2018a} and   in Definition~\ref{def-wildtypes} coincide. Moreover, the examples of wild actions constructed in \cite[Section~8]{HL2018a} are all of   finite and flat type, while the examples of   actions constructed in Section~\ref{sec-arboreal} below are dynamically wild and not of finite type.
On the other hand,   Corollary~1.7 of \cite{HL2018b} implies that if $G$ is a finitely generated nilpotent group, then every Cantor action of $G$ must be both stable and  algebraically stable.

 \section{Analytic regularity for stable Cantor actions}\label{sec-dynamics}

In this section, we relate the bounded property of  the stabilizer  group $\Upsilon_s(\Phi)$ for a  Cantor action $(\fX, G, \Phi)$ with    ``analytic regularity'' of the action.    We first recall some background context.

  \subsection{Locally quasi-analytic actions}
Haefliger introduced in   \cite{Haefliger1985} the notion of a \emph{quasi-analytic} topological action of a   pseudogroup  on a connected space $X$. 
  {\'A}lvarez L{\'o}pez and Candel in   \cite[Definition~9.4]{ALC2009}, and  later 
   {\'A}lvarez L{\'o}pez and Moreira Galicia in \cite[Definition~2.18]{ALM2016}, adapted the notion of a   quasi-analytic topological action to the more general case  where the action space $X$ need not be connected.
The authors formulated in \cite{HL2018a,HL2018b}   a notion of \emph{locally quasi-analytic} Cantor actions, and showed the relation between this condition  and the stable property for the action.

The quasi-analytic condition for a   Cantor action is a   modification of the notion of a \emph{topologically free}   action for general Cantor actions, which   first appeared in the work of Boyle and Tomiyama \cite{BoyleTomiyama1998} in their study of flip-conjugacy.    Renault   showed in \cite[Section~3]{Renault2008} that an action is topologically free if and only if   the associated action groupoid is \emph{essentially principal}.     We first recall the definition    of a topologically free action and some properties of this definition. Topological freeness and   related ideas are discussed in more detail in  \cite[Section~2]{Li2018}.  
 
The   \emph{isotropy group} at $x \in \fX$  of an  action   $(\fX, G, \Phi)$ is 
$G_x = \{ g \in G \mid g \cdot x = x\}$.
A point $x \in \fX$ is said to have trivial isotropy if   $G_x$ is the trivial group. All points in the orbit of $x$ then also have trivial isotropy, so form a dense subset of $\fX$ for a minimal  action.

Let $\Fix(g) = \{x \in \fX \mid \Phi(g) \cdot x = x \}$, then introduce the \emph{isotropy set}
\begin{equation}\label{eq-isotropy}
 \Iso(\fX,G,\Phi) \equiv \{ x \in \fX \mid \exists ~ g \in G ~ , ~ g \ne id ~, ~ g \cdot x = x    \} = \bigcup_{e \ne g \in G} \ \Fix(g) \ . 
\end{equation}
The action $(\fX,G,\Phi)$ is said to be \emph{topologically free}   if the set $\Iso(\fX,G,\Phi) $ is meager   in $\fX$.  
Note that if $e \ne g \in G$ and $\Phi(g)$ acts trivially on $\fX$, then $\Iso(\fX,G,\Phi)  = \fX$, and thus a topologically free action must be faithful. 
If the group $G$ is abelian, it is an exercise to show that a faithful minimal Cantor action $(\fX,G,\Phi)$ must be topologically free; see for example \cite[Corollary~2.3]{HL2018b}. 

We next restrict attention to Cantor actions, and use the special properties of equicontinuous actions to formulate   local forms  of the topologically free property. 
A Cantor action    $(\fX,G,\Phi)$ is said to be  \emph{locally topologically free} if there exists $\e > 0$ such that for any adapted set $U$ with $\diam(U) < \e$, the action of $G_U$     on $U$ is topologically free. There is another related notion, defined as follows:
  \begin{defn} \cite[Definition~9.4]{ALC2009} \label{def-LQA} A Cantor action       $(\fX,G,\Phi)$  is   \emph{locally quasi-analytic}, or simply   \emph{LQA}, if there exists $\e > 0$ such that for any adapted   set $U \subset \fX$ with $\diam (U) < \e$,  and  for any adapted subset $V \subset U $, and elements $g_1 , g_2 \in G$
 \begin{equation}\label{eq-lqa}
  \text{if the restrictions} ~~ \Phi(g_1)|V = \Phi(g_2)|V, ~ \text{ then}~~ \Phi(g_1)|U = \Phi(g_2)|U. 
\end{equation}
That is, by \cite[Proposition~2.2]{HL2018b},   the action of $H_U = \Phi(G_U)$ on $U$ is topologically free.
\end{defn}
If  \eqref{eq-lqa} holds for $U=\fX$, then the action of $\oG$ is topologically free.

   Examples of equicontinuous Cantor actions   $\fX$ which are locally quasi-analytic, but not quasi-analytic,  are easily constructed, as given in \cite{DHL2016c,HL2018a} for example.
There is also a  generalization of the locally quasi-analytic property for the action of the profinite closure group $\oG$. 
  \begin{defn}  \label{def-LCQA} A Cantor action       $(\fX,G,\Phi)$  is   \emph{locally completely quasi-analytic}, or simply   \emph{LCQA}, if there exists $\e > 0$ such that for any adapted   set $U \subset \fX$ with $\diam (U) < \e$,  and  for any adapted subset $V \subset U $, and elements $\whh_1 , \whh_2 \in \oG$
 \begin{equation}\label{eq-lcqa}
  \text{if the restrictions} ~~ \whh_1 |V = \whh_2 | V, ~ \text{ then}~~ \whh_1 |U = \whh_2|U. 
\end{equation}
Equivalently, set $\whh = \whh_1 \ \whh_2^{-1}$ then we have 
 \begin{equation}\label{eq-lcqa2}
  \text{if the restriction} ~~ \whh |V = id | V, ~ \text{ then}~~ \whh |U = id | U. 
\end{equation}
\end{defn}

Since $\Phi(G) \subset \oG$, the LCQA property implies the LQA property for a Cantor action. 
It is not known however, if there exists a Cantor action    $(\fX,G,\Phi)$  which is LQA but not LCQA. In particular, does there exists a free Cantor action which is not LCQA?  
All examples known to the authors which are not LCQA are also not LQA. 

\subsection{Bounded stabilizer groups}\label{subsec-bddmno}

Recall from Definition~\ref{def-wildtypes} that 
a Cantor action     $(\fX,G,\Phi)$ is  stable  if    the stabilizer  group $\Upsilon_s(\Phi)$   is bounded.
 Here is the main result of this section. 
\begin{thm}\label{thm-boundedLCQA}
Let        $(\fX,G,\Phi)$  be a Cantor action, then $\Phi$ is a locally completely quasi-analytic (LCQA) action if and only if its stabilizer limit group  $\Upsilon_s(\Phi)$ is   bounded.
\end{thm}
 \proof
Let  $(\fX, G, \Phi)$ be a Cantor action, and 
  $\cG^x_{\cU} = \{G_{\ell}\mid \ell \geq 0 \}$ be the group chain associated to an adapted neighborhood basis $\cU$ at $x \in \fX$. 
  Let $\whG_{\infty}$ be defined as in \eqref{eq-invlimgroup}, then by Theorem~\ref{thm-quotientspace}, there is a topological isomorphism 
$\whTheta \colon \oG \to \whG_{\infty}$ which induces an isomorphism $\whTheta \colon  \oG_x \to \cD_x$.

  Let $\cK(\cU) = \{K_{\ell} \subset \cD_x \mid \ell \geq 0\}$ be the increasing chain of stabilizer subgroups defined in \eqref{eq-isotropykernel}.

The map $\whTheta$ induces a  homeomorphism of $G$-spaces
$\Theta_x \colon  \oG/\oG_x \cong  \fX  \to  \whG_{\infty}/\cD_x \cong X_{\infty}$
by Corollary~\ref{cor-inducedquotients}. Let $X_{\infty}$ have a  metric such  that $G$ acts on $X_{\infty}$ by isometries, then let $\dX$ be the metric on $\fX$ induced by the homeomorphism $\Theta_x$.

Recall that for $\ell \geq 0$,  $\Theta_x$ identifies the clopen set  $U_{\ell} \subset \fX$ with $\whU_{\ell} \equiv \whG_{\ell}/\cD_x \subset \whG_{\infty}/\cD_x$.
 
Assume that  $\Phi$ is LCQA, and let $\e > 0$ be as in Definition~\ref{def-LCQA}. Let $\ell_0 \geq 1$ be such that $\diam(U_{\ell}) < \e$ for all $\ell \geq \ell_0$. Thus, for $\ell > \ell_0$ the restricted action of $\whG_{\ell}$ on $\whU_{\ell}$ is quasi-analytic. Given $\whh \in \cD_x$ suppose that 
$\whh$  acts trivially on   $\whU_{\ell}$ then it must act trivially on   $U_{\ell_0}$.
That is,   $K_{\ell} \subset K_{\ell_0}$.   As the converse $K_{\ell_0} \subset K_{\ell}$ always holds, this implies that  $\Upsilon_s^x(\Phi)$ is a bounded direct limit group.

Conversely, assume that  $\Upsilon_s^x(\Phi)$ is a bounded direct limit group.
 Then there exists $\ell_0 > 0$ such that $K_{\ell_0} = K_{\ell}$  for all $\ell \geq \ell_0$.
 Let $\e_0 > 0$ be such that the ball of radius $\e_0$ about $x$ is contained in $U_{\ell_0}$.
 
Set $\e = \e_0/2$. Let $U \subset \fX$ be an open set with with $\diam(U) < \e$, and let $V$ be an open subset with $V \subset U$.
Then by minimality and equicontinuity of the action $\Phi$, there exists $g \in G$ such that $g \cdot x \in V$. 
By the choice of $\e$ we have that $U' = g^{-1} \cdot U \subset U_{\ell_0}$.
In addition, as $\cU$ is a neighborhood basis at $x$, there exists $\ell > 0$ such that $U_{\ell} \subset V' = g^{-1} \cdot V$. 
 
 Let $\whh \in \oG$ such that $\whh$ acts as the identity on $V$, then $\whh' = g^{-1} \ \whh \ \cdot g$ acts as the identity on $V'$, and thus also on  $U_{\ell}$ since $U_{\ell} \subset V'$. In particular, $\whh' \cdot x = x$, so $\whh' \in \oG_x$ hence $\whk = \whTheta(\whh') \in \cD_x$. 
 Moreover, $\whk$ acts as the identity on $\whU_{\ell}$ hence $\whk \in K_{\ell}$.
  We are given that $K_{\ell} = K_{\ell_0}$ so $\whh' \in K_{\ell_0}$ hence acts as the identity on $\whU_{\ell_0}$. 
  Thus, $\whh'$ acts as the identity on $U' \subset U_{\ell_0}$ and so $\whh$ acts as the identity on $U$, as was to be shown.
 \endproof

 \section{Orbit equivalence and rigidity}\label{sec-rigidity} 

Let $(\fX,G,\Phi)$ and $(\fX',G,\Psi)$ be    general Cantor actions.     A homeomorphism $h \colon \fX \to \fX'$ is a \emph{conjugacy} between the two actions  if $\Phi(g) = h^{-1} \circ \Psi(g) \circ h$ for all $g \in G$. A conjugacy map   $h$ thus preserves the   structure of the orbits as $G$-spaces, which is used in constructing conjugacy invariants for an action.
For example, this was used in Theorem~\ref{thm-centraldirectlimgroup} to show that the isomorphism classes of the direct limit groups    $\Upsilon_s(\Phi)$ and $\Upsilon_c(\Phi)$ are invariants of the conjugacy class of a Cantor action  $(\fX, G, \Phi)$.  

An \emph{orbit equivalence} between two actions is a bijective map $h \colon \fX \to \fX'$ which maps orbits of the action $\Phi$ to orbits of the action $\Psi$. In addition, one can impose additional assumptions on the map $h$, such as to assume that it  is a measurable isomorphism   with respect to given  quasi-invariant measures on $\fX$ and $\fX'$.  In the measurable category of actions, this yields the notion of \emph{measurable orbit equivalence}, a topic which has been extensively studied. A celebrated result  by    Connes, Feldman and Weiss in \cite{CFW1981}  showed that for essentially free actions of  amenable groups,   any two such actions are measurably orbit equivalent. For example, if $G$ is a nilpotent group, then any essentially free ergodic action of $G$ is measurably orbit equivalent to an ergodic $\mZ$-action.  Thus, a measurable orbit equivalence need not preserve the   structure on orbits induced by the $G$-action. On the other hand, the   results of Furman in \cite{Furman1999} show that for measure preserving actions of higher rank lattice groups, measure orbit equivalence implies virtual conjugacy of the actions, hence the orbit structures are preserved in this case.  For further discussion of properties of measurable orbit equivalence, see for example  the survey by Gaboriau \cite{Gaboriau2010}. 

For general Cantor actions, \emph{continuous orbit equivalence} is the more natural notion to study.  One assumes  that the orbit equivalence $h$ is a homeomorphism which, in addition,   satisfies a ``locally constant'' property, as stated precisely in Definition~\ref{def-torb1}. This notion 
was introduced by    Boyle in his thesis \cite{Boyle1983}, and see also \cite{GW1995}, and has played a fundamental role in  the classification of general Cantor actions in many subsequent works, as for example in  \cite{BezuglyiMedynets2008,GPS1999,GMPS2010}. The works of Renault \cite{Renault2008} and Li \cite[Theorem~1.2]{Li2018} use the notion of continuous orbit equivalence to classify the  groupoid $C^*$-algebras associated to the action, in the case of  topologically free actions.
The related notion of the topological full group of an action   has   provided a rich source of examples of finitely generated groups with exceptional properties, as discussed for example in \cite{deCornulier2014,Nekrashevych2018}.  

The \emph{rigidity property} for Cantor actions states that topological orbit equivalence implies virtual conjugacy of the actions.  In essence, this is saying that an orbit equivalence which locally preserves the orbit structures, must preserve the orbit structures  on a global scale.  In this way, it can be seen as a property which is analogous to the measure rigidity property studied by  Furman in \cite{Furman1999}.  The work  of Cortez and Medynets   in \cite{CortezMedynets2016}  and Li   in  \cite{Li2018} prove  versions of   topological rigidity.

In this section, we study the consequences of continuous orbit equivalence between Cantor actions, without any assumption that the actions are topologically free.  Our main result is Theorem~\ref{thm-coestable}, which considers the   relation  between the stabilizer  groups  for continuously orbit equivalent actions. 
We deduce  that  the property of being  stable, as well as  algebraically stable, is an invariant of continuous orbit equivalence.

We first recall in Section~\ref{subsec-coe} some notions from the theory of orbit equivalence for Cantor actions. 
In Section~\ref{subsec-rigiditynotions} we discuss various notions of rigidity for Cantor actions. 
 Then   in Section~\ref{subsec-rigidity} the rigidity results from \cite{CortezMedynets2016,Li2018} are described. We then begin  the proof of Theorem~\ref{thm-coestable} in Section~\ref{subsec-rigidity}, and  Section~\ref{subsec-induced} contains some technical results needed to complete the proof. 
 Theorem~\ref{thm-coe=re} gives an extension of Theorem 3.3 by Cortez and Medynets in \cite{CortezMedynets2016}, and  Theorem~1.5   in \cite{HL2018b} by the authors.

 \subsection{Continuous orbit equivalence}\label{subsec-coe}

Let  $(\fX,G,\Phi)$ be a general Cantor action.  Consider the equivalence relation on $\fX$ defined by the action, 
 \begin{equation}\label{eq-ERX}
\cR(\fX, G, \Phi) \equiv \{(x,  g \cdot x)) \mid x \in \fX, g \in G\} \subset \fX \times \fX ~.
\end{equation}
Given   actions $(\fX,G,\Phi)$ and $(\fX',H,\Psi)$, we say they are \emph{orbit equivalent} if there exist a bijection   $h \colon \fX \to \fX'$ which maps 
$\cR(\fX, G, \Phi)$ to $\cR(\fX', H, \Psi)$, and similarly for the inverse map $h^{-1}$.

\begin{defn}\label{def-torb1}
Let $(\fX,G,\Phi)$ and $(\fX',H,\Psi)$ be    general Cantor actions.
A \emph{continuous orbit equivalence}  between the actions  is a homeomorphism    $h \colon \fX \to \fX'$ which is an orbit equivalence, and satisfies the locally constant conditions:      
\begin{enumerate}
\item  for each $x \in \fX$ and $g \in G$, there exists $\alpha(g,x) \in H$ and an open set $x \in U_x \subset \fX$ such that 
$ \Psi(\alpha(g,x))  \circ h | U_x = h \circ \Phi(g) | U_x$; \\
\item for each $y \in \fX'$ and $k \in H$, there exists $\beta(k,y) \in G$ and an open set $y \in V_y \subset \fX'$ such that 
$ \Phi(\beta(k,y))  \circ h | V_y = h \circ \Psi(k) | V_y$.
\end{enumerate}
\end{defn}
Note in particular that these conditions imply that the functions $\alpha \colon G \times \fX \to H$ and $\beta \colon H \times \fX' \to G$ are continuous, as the groups $G$ and $H$ have the discrete topology.

One special class of examples of continuous orbit equivalences, are those for which  the functions $\alpha(g,y)$ and  $\beta(g,x)$  are constant in $y$. Then the   identities \eqref{eq-conjmap1}  and \eqref{eq-conjmap2}  can be considered as   defining ``time-shifts'' along the orbits, in analogy with the case when $G, H = \mZ$.

 Let $(\fX',H,\Psi)$  be   a general Cantor action,   let $h \colon \fX \to \fX'$ be a homeomorphism which implements a continuous orbit equivalence between   $(\fX,G,\Phi)$ and $(\fX',H,\Psi)$, and let  
\begin{equation}\label{eq-conjugateaction2}
\Psi^h \colon H \times \fX \to \fX ~ , ~ {\rm where} ~ \Psi^h(g)(x) = h^{-1}  \circ \Psi(g) \circ h(x) ~, ~ {\rm for} ~ g \in H ~ .
\end{equation}
Then the orbits for the actions $\Psi^h$ and   $\Phi$ are equal, so the identity map is an orbit equivalence between $\cR(\fX, G, \Phi)$ and $\cR(\fX, H, \Psi^h)$. Thus, for the study of orbit equivalence,  it suffices to consider orbit equivalences for which $h$ is the identity map. Then the data of the orbit equivalence is encapsulated in the  continuous, hence locally constant  maps $\alpha \colon G \times \fX  \to H$  and $\beta \colon H \times \fX  \to G$. That is, for $y \in \fX$ and $g \in G$, there exist a clopen set $y \in U_{g,y} \subset \fX$ 
 so that 
 \begin{equation}\label{eq-conjmap1} 
\Psi(\alpha(g,y)) \cdot z      =      \Phi(g) \cdot z   ~ {\rm for} ~   z \in U_{g,y}   ~ ,
\end{equation}
and for $k \in H$, there exists a clopen set $y \in V_{k,y} \subset \fX$ so that 
 \begin{equation}\label{eq-conjmap2}
\Phi(\beta(k,y)) \cdot z      =      \Psi(k) \cdot z   ~ {\rm for}  ~ z \in V_{k,y} ~ .  
\end{equation}
 
Next, recall the notion of a cocycle over an action.
\begin{defn}\label{def-cocycle}
Let $\vp \colon \G \times X \to X$ be an action of the group $\G$ on a space $X$. A map $c \colon \G \times X \to \G'$ is a \emph{cocycle over} the action if for all $\gamma, \gamma' \in \G$ and $x \in X$, we have
\begin{equation}\label{eq-cocycleidentity}
c(\gamma' \cdot \gamma, x) = c(\gamma', \vp(\gamma)(x)) \cdot c(\gamma,x) \ .
\end{equation}
\end{defn}

For an action which is free, or more generally is topologically free, then the identities \eqref{eq-conjmap1}  and \eqref{eq-conjmap2} uniquely determine the values of $\alpha(g,y)$ and $\beta(k,y)$ for all $y \in \fX$ and $k \in G$, and thus the functions $\alpha$ and $\beta$ are cocycles over the actions $\Phi$ and $\Psi$, respectively. However, if the actions are not topologically free, then the identities \eqref{eq-conjmap1}  and \eqref{eq-conjmap2} no longer uniquely specify the values of $\alpha(g,y)$ and $\beta(k,y)$, and so they need not   satisfy   the cocycle identity.

For continuous orbit equivalent topologically free actions whose orbit functions  $\alpha$ and $\beta$ are constant,  the cocycle identity \eqref{eq-cocycleidentity}  implies they define   group homomorphisms. This leads to the notion of 
    \emph{conjugacy} of Cantor actions, as was studied by   Li  in \cite{Li2018},  and also the related notion of \emph{structural conjugacy} by Cortez and Medynets in  \cite{CortezMedynets2016}. We use the following  terminology:   
\begin{defn}\label{def-thetaconjugacy}
The Cantor actions $(\fX, G, \Phi)$ and $(\fX, H, \Psi)$ are said to be  \emph{$\theta$-conjugate} if there exists  a group isomorphism $\theta \colon G \to H$  such that  
\begin{equation}\label{eq-constant}
 \Phi(g)(x) = \Psi^{\theta}(x) \equiv \Psi(\theta(g))(x) \quad , ~ {\rm for ~ all} ~ g \in G ~ ,~  x \in \fX ~ .
\end{equation}
\end{defn}
The identity \eqref{eq-constant} implies that we can choose the function $\alpha(g,x) = \theta(g)$ in Definition~\ref{def-torb1}.  
 
For the case where $G = H = \mZ$, the involution $\theta(n) = -n$ is the only non-trivial isomorphism, and in this case $\theta$-conjugacy is the same as flip-conjugacy as   studied by  Boyle and Tamiyama \cite{BoyleTomiyama1998}.

There is a related notion to Definition~\ref{def-thetaconjugacy},  where we assume that $\fX = \fX'$,  $G = H$,  and      $\theta$ is the identity. Define the \emph{center} of the action:
\begin{equation}\label{eq-aut}
\Aut(\fX, G, \Phi) = \{ h \in \Homeo(\fX) \mid  h(\Phi(g) \cdot x) = \Phi(g) \cdot h(x) ~, ~ {\rm for ~ all} ~ x \in \fX, ~ g \in G\} ~.
\end{equation}
 
Then we have the following relation between the direct limit group $\Upsilon_c^x(\Phi)$ and the center:  
  \begin{lemma}\label{lem-symmetry}
 Let  $(\fX, G, \Phi)$ be a Cantor action, and 
  $\{G_{\ell}\mid \ell \geq 0 \}$ be the group chain associated to an adapted neighborhood basis $\cU$ at $x \in \fX$.
Then for each $\ell \geq 1$,   $Z_{\ell} \subset \Aut(\fX, G_{\ell}, \Phi)$.
 \end{lemma}
 \proof
 Let $\whh \in Z_{\ell}$ which acts on $\whG_{\infty}$ by left multiplication. Then for any $\whg \in G_{\ell}$ we have $\whg \ \whh = \whh \ \whg$ by Lemma~\ref{lem-centralizer}. Thus,  for the isomorphism $\whTheta \colon    \oG \to \whG_{\infty}$ we have $h = \whTheta^{-1}(\whh) \in \Aut(\fX, G_{\ell}, \Phi)$.
 \endproof
 
  \begin{remark}
  {\rm 
 Note that in the context of Lemma~\ref{lem-symmetry},   while the action of $G_{\ell}$ on the clopen set $U_{\ell}$ is minimal, its action on $\fX$ cannot be minimal, as $U_{\ell}$ is adapted hence invariant by the action of $G_{\ell}$.    Lemma~\ref{lem-symmetry} implies that when the action $\Phi$ is restricted to subactions by groups $G_{\ell}$, the groups $Z_{\ell}$ arise as the centers of these restricted actions. That is,  while $\Phi$ may have trivial center, if the direct limit group $\Upsilon_c^x(\Phi)$ is non-trivial, then there exists   restricted actions with non-trivial ``symmetries''.  }
  \end{remark}

 \subsection{Types of rigidity}\label{subsec-rigiditynotions}

Rigidity for a dynamical system can be viewed as asserting that two dynamical systems which are equivalent in some  weaker  sense, are also equivalent, possibly up to some finite indeterminacy, in some  stronger  sense. 
For Cantor systems, we formulate three notions of equivalence for their rigidity, each weaker than the previous one.

The first notion is based on  $\theta$-conjugacy as given in Definition~\ref{def-thetaconjugacy}.
 \begin{defn}\label{def-rigidity}
 A Cantor action $(\fX, G, \Phi)$ is  \emph{rigid} if, given a continuous orbit equivalence $h \colon \fX \to \fX'$ to a Cantor action   $(\fX', G', \Psi)$, then     $\Psi^h$ as defined by  \eqref{eq-conjugateaction2}  is    $\theta$-conjugate to $\Phi$.
\end{defn}
This is the form of the rigidity property  used by Li  \cite{Li2018}.

It is elementary to construct examples of Cantor actions where $G$ is a cross-product of a normal subgroup by a finite quotient group, 
as in Examples~A.2 and A.3   in \cite{HL2018b}, such that   two   actions are continuous orbit equivalent  but are not $\theta$-conjugate actions.  
This is the situation considered by   Cortez and Medynets in \cite{CortezMedynets2016}. They accordingly introduced the weaker notion of 
\emph{structural stability}  and proved their rigidity results for free actions in these terms. Their notion of structural stability coincides with what we call virtual rigidity below.

Let $U \subset \fX$ be adapted for the Cantor action $(\fX, G, \Phi)$. Let $\Phi_U \colon G_U \times U \to U$ denote the restricted action of $G_U$ on $U$. Similarly, for a Cantor action $(\fX', G', \Psi)$ with adapted set $V \subset \fX'$, let $(V, G'_V , \Psi_V)$ denote the restricted action.

\begin{defn}\label{def-Vrigidity}
Let $(\fX, G, \Phi)$  be a Cantor action. Then the   action $(\fX, G, \Phi)$ is  \emph{virtually rigid} if, given a continuous orbit equivalence $h \colon \fX \to \fX'$ to a Cantor action   $(\fX', G', \Psi)$,     there exists  an adapted set $U \subset \fX$ for the action $\Phi$ such that   $V = h(U)\subset \fX'$ is adapted for   the action $\Psi$, and there is  an  isomorphism $\theta_U \colon G_U \to G'_V$ so that the action $(U, G_U, \Phi_U)$ is    $\theta_U$-conjugate to $(V, G'_V, \Psi_V^h)$.
\end{defn}

Note that both Definitions~\ref{def-rigidity} and \ref{def-Vrigidity} are essentially only applicable for topologically free actions on the full space $\fX$. 
The  third  notion we consider is \emph{return equivalence}, which was used in the authors' work \cite{HL2018a} for the study of the homeomorphism types of weak solenoids. For the geometric applications in \cite{DHL2016c,HL2018a},  the holonomy action on a transversal is the fundamental concept. Accordingly, return equivalence for Cantor actions is formulated in terms of the image group $H_U$  for an adapted subset $U \subset \fX$ as defined in \eqref{eq-HUV}.

Let $U \subset \fX$ be adapted for the Cantor action $(\fX, G, \Phi)$. Let $\oPhi_U \colon H_U \times U \to U$ denote the induced action of $H_U = \Phi(G_U) \subset \Homeo(U)$ on $U$. Similarly, for a Cantor action $(\fX', G', \Psi)$ with adapted set $V \subset \fX'$, let $(V, H'_V , \oPsi_V)$ denote the induced action. 

   \begin{defn}\label{def-return}
Two Cantor actions $(\fX, G, \Phi)$ and $(\fX', G', \Psi)$ are  \emph{return equivalent} if there exists an adapted set $U \subset \fX$ for the action $\Phi$, an adapted set $V \subset \fX'$ for the action $\Psi$, and a homeomorphism $h_U \colon U \to V$ which induces a   $\theta_U$-conjugacy between the   action of $H_U$ on $U$ and the action of $H'_V$ on $V$.
\end{defn}
 
While the Definitions~\ref{def-Vrigidity} and \ref{def-return} are similar,   the former requires an isomorphism between the subgroups $G_U \subset G$ and $G'_V \subset G'$, while     the latter only requires  an isomorphism    between their respective  image groups $H_U$ and $H'_V$. 
Thus, Definition~\ref{def-return} is most relevant for the study of Cantor actions
which are not topologically free.

 For topologically free actions, an adaptation of the proof in Theorem
3.3 of \cite{CortezMedynets2016} shows that the fact that the isomorphism $\theta$ is induced
by a continuous orbit equivalence $h$ forces the groups $G_U$ and $G_V'$
to have equal indices in $G$ and $G'$ respectively. In Definition~\ref{def-return} of
return equivalence, the homeomorphism $h_U$ need not be induced by a
homeomorphism of the space $\fX$, so there is no requirement that the
groups $H_U$ and $H_V'$ have  the same index with respect to any
larger groups.
 
\subsection{Rigidity for Cantor actions}\label{subsec-rigidity}

The approach  in the literature to proving     rigidity for Cantor actions is to assume that both actions are free on an invariant dense subset   $\fZ \subset \fX$, then observe that the function   $\alpha$   in Definition~\ref{def-torb1}   satisfies a cocycle identity on $\fZ$, hence   on all of $\fX$ by continuity. The map $\alpha \colon G \times \fX \to H$ is   called the   ``orbit cocycle'' for the orbit equivalence. 

Then either a   cohomological assumption as in Li \cite{Li2018}, or a dynamical assumption as in Cortez and Medynets \cite{CortezMedynets2016},   is used to show that the cocycle $\alpha$ is cohomologous to a constant cocycle,  which implies  that the conditions of Definition~\ref{def-thetaconjugacy} are satisfied. 

The results of Boyle and Tamiyama in \cite{BoyleTomiyama1998} can be interpreted as saying that a minimal Cantor action by $\mZ$ is rigid, as  stated in   \cite[Theorem~3.2]{Li2018}.   As $G$ is abelian and the action is faithful in this case, there is a well-defined orbit cocycle,   and the authors show that   it is cohomologous to a constant.

 Cortez and Medynets showed that \emph{free} equicontinuous Cantor actions are rigid in \cite[Theorem 3.3]{CortezMedynets2016}, and in fact their proof directly extends to show:
  \begin{thm}\label{thm-CM}
  Let $G$, $G'$ be    finitely generated groups, and $(\fX,G,\Phi)$  and $(\fX', G', \Psi)$  be   faithful Cantor actions. Suppose that the actions $\Phi$ and $\Psi$ are topologically free.  Then the actions are continuous orbit equivalent if and only they are virtually rigid in the sense of Definition~\ref{def-Vrigidity}. 
  \end{thm}
    The  proof   of  this result in \cite{CortezMedynets2016} does not explicitly discuss the orbit cocycle. Rather, it   uses the odometer model for an equicontinuous action   to   construct the     isomorphism $\theta_U$ in Definition~\ref{def-Vrigidity} directly from the action.

For $G$ a finitely generated group, 
  Theorems~1.3 and 1.4  in Li \cite{Li2018} give  conditions on a topologically free Cantor action  which imply that the orbit cocycle for the action is cohomologous to a constant cocycle.  The discussion in \cite[Section~4]{Li2018} explains the passage from a constant cocycle to a conjugacy.

The methods in  \cite{CortezMedynets2016,Li2018} were used in the authors' work    \cite[Section~4]{HL2018a} to show that stable actions which are continuous orbit equivalent are \emph{return equivalent}. We  also use these techniques in the proof of the following result, which will be used in the proof of our rigidity result Theorem~\ref{thm-coe=re}.

\begin{thm}\label{thm-coestable}
Let $h \colon \fX \to \fX'$ be  a continuous orbit equivalence between   Cantor actions $(\fX,G,\Phi)$  and $(\fX', G', \Psi)$. If $G$  is finitely generated  and  $(\fX', G', \Psi)$ is stable, then  $(\fX,G,\Phi)$ is stable.
  \end{thm}
   \proof
 Following the discussion in Section~\ref{subsec-coe}, we can assume that $\fX = \fX'$ and $h$ is the identity map. 
Fix a basepoint $x \in \fX$. Then by Theorem~\ref{thm-boundedLCQA}, there exists $V_1 \subset \fX$ adapted to the action $\Psi$ with $x \in V_1$, such that the action of  $H'_1 = \Psi(G'_{V_1}) \subset \Homeo(V_1)$ on $V_1$ is topologically free. Let $\cZ \subset V_1$ be a dense subset which is $H'_1$-invariant and the restriction of the action of $H'_1$ to $\cZ$ is free.

Let $W \subset \fX$ be   adapted set for the action $\Phi$ with $x \in W \subset V_1$. Let $G_{W} \subset G$ be defined as in  \eqref{eq-adapted}.

  Let  $\alpha \colon G \times \fX  \to G'$ be the orbit function  which satisfies the relation  \eqref{eq-conjmap1}, and let 
 $\alpha_W \colon  G_W \times W \to G'$ denote its restriction to $W$.
  
  For each $g \in G_W$ and $y \in W$ we have $\Phi(g)(y) \in W$. Let $h \in G'$ be such that $\Psi(h)(y) = \Phi(g)(y)$, then $W \subset V_1$ implies that $\Psi(h)(V_1) \cap V_1 \ne \emptyset$ hence $\Psi(h) \in H'_1 = \Psi_V(G'_{V_1}) \subset \Homeo(V_1)$. That is,  
  the restriction of $\alpha$ to $G_W \times W$ induces a map  $\whalpha_W = \Psi_{V'} \circ \alpha \colon  G_W \times W \to H_1'$, where
  $\Psi_{V'} \colon G_1' \to H_1'$.

The action of $H'_1$ on $V_1$ is topologically free, so the proof of   \cite[Lemma~2.8]{Li2018} adapts   to yield:

 \begin{lemma}\label{lem-cocycle}
 $\whalpha_W \colon  G_W \times W \to H'_1$ satisfies the cocycle identity \eqref{eq-cocycleidentity}    for the restricted action   $\Phi_W \colon G_W \times W \to W$.
 \end{lemma}
    
 We next show that a properly chosen restriction of the cocycle $\whalpha_W$ is a coboundary. The proof of this fact follows closely the  proof of \cite[Theorem~3.3]{CortezMedynets2016}, with the variation   that we   only assume   the   action of $H'_1$ is topologically free, and do not assume   the   action of $G_W$ on  $W$   is topologically free.

Since $G$ is finitely generated, the same holds for the subgroup $G_W$ of finite index. Choose  a symmetric generating set  $\{g_i \mid  1 \leq i \leq m \}$ for $G_W$.
That is, each $g \in G_W$ can be written as a product 
 $g = g_{i_1} \ g_{i_2} \ \cdots \ g_{i_{\nu}}$ where each $1 \leq i_{\ell} \leq m$.
 The map $\alpha \colon G \times \fX \to G'$ is continuous,  and $G'$ is discrete, hence   $\alpha_W$ is locally constant. 
 Thus,   there exists  $\delta_1 > 0$ so that
\begin{equation}\label{eq-locconstant}
  \alpha_W(g_i, y) = \alpha_W(g_i, z) ~ {\rm    for } ~    1 \leq i \leq m ~ {\rm    and } ~ y,z \in \fX ~ {\rm with} ~ \dX(y,z) < \delta_1 \ .
\end{equation}
 
 Let $W' \subset W$ be an adapted set such that $x \in W'$ and $\diam (g \cdot W') < \delta_1$ for all $g \in G_W$. 
 
 The collection $\cW' \equiv \{g \cdot W' \mid g \in G_W\}$ is a finite clopen partition of $W$, so there exists $0 < \delta_2 < \delta_1$ such that 
     for any $y \in W$   there exists $g \in G_W$ such that $B_{\dX}(y,\delta_2) \subset g \cdot W'$, 
 and $g \cdot W'$ is the unique element of the partition    $\cW'$ which contains it.

By the uniform continuity of the action $\Phi$, there exists  $0 < \delta_3 \leq \delta_2$ such that 
\begin{equation}\label{eq-delta3}
\dX(\Phi(g)(y) ,\Phi(g)(z)) < \delta_2 ~ {\rm for ~ all} ~  y,z \in \fX ~ {\rm with} ~  \dX(y,z) < \delta_3 \ , \  {\rm for ~ all} \  g \in G \ .
\end{equation}

 Choose an adapted set $U_1$  such that $x \in U_1 \subset W'$ and $\diam (g' \cdot U_1) < \delta_3$ for all $g' \in G_{U_1}$.  
 That is, once $W'$ has been chosen,  choose $U_1$ sufficiently small so that all  orbits of points in $U_1$ under the action of $\Phi(G_W)$ have well-defined codings with respect to the partition $\cW'$ of $W$ formed by the translates of $W'$. (The coding method is discussed in detail in \cite[Section~6]{ClarkHurder2013}.)

 Set $G_1 = G_{U_1} \subset G_W$,   let $H_1 = \Phi_{U_1}(G_1) \subset \Homeo(U_1)$, and let $\Phi_1 \colon G_1 \to H_1$ denote the quotient map. We also denote the restricted action of $G_1$ by $\Phi_1 \colon G_1 \times U_1 \to U_1$. Let  $\whalpha_1 \colon G_1 \times U_1 \to H'_1$ 
  denote the cocycle over $\Phi_1$ obtained by restricting $\whalpha_W$. The next result shows that $\whalpha_1$ is a coboundary.
 The  proof below uses a key idea from   the proof of \cite[Theorem~3.3]{CortezMedynets2016}.

\begin{prop}\label{prop-coboundary}
  $\whalpha_1 \colon G_1 \times U_1 \to H'_1$ is   induced by a   homomorphism 
$\whtheta_1 \colon G_1 \to H'_1$. That is, for $g \in G_1$   we have $ \Psi(\whalpha_1(g,x)) = \whtheta_1(g)$, and thus
\begin{equation}\label{eq-conjugatequotient}
\Phi(g)(y)  =     \whtheta_1(g)(y)   \quad {\rm for ~ all} ~  y \in U_1 ~ .
\end{equation}
 \end{prop}  
  \proof
We first show that for $g \in G_1$ the function $\whalpha_1(g, y)$ is constant in $y \in U_1$. 
 
Recall that $\{g_i \mid  1 \leq i \leq m \}$ is a symmetric generating set for $G_W$. 
By the choice of $\delta_1$,   for $\dX(y,z) < \delta_1$ and $1 \leq i \leq m$,   we have   $\whalpha_1(g_i ,y) = \whalpha_1(g_i ,z)$.
Thus  by the choice of $W'$ and \eqref{eq-locconstant},  for any $g' \in G_W$,  the value of $\whalpha_W(g_i,y)$ is constant for $y \in g' \cdot W'$. 

Then by the choice of $\delta_3$ so that \eqref{eq-delta3} holds, and the choice of $U_1$ so that $\diam (g' \cdot U_1) < \delta_3$ for  $g' \in  G_W$, then for   $y,z \in  U_1$ we have $g' \cdot y   \in g' \cdot U_1$ and $g'\cdot z \in g' \cdot U_1$. Thus,  $g' \cdot y $ and  $g'\cdot z$  are contained in the same set of the partition $\cW'$, hence  $\ds  \whalpha_1(g_i , g' \cdot y) = \whalpha_1(g_i ,  g' \cdot z)$ for any $1 \leq i \leq m$.

We now apply this recursively to the elements of $G_1$. Let $g \in G_1$, then   $g = g_{i_1} \cdots g_{i_{\nu}} \in G_1$. 
Set $g' = g_{i_2} \cdots g_{i_{\nu}}$ then for $y,z \in  \cZ \cap U_1$, 
\begin{eqnarray}
 \lefteqn{ \whalpha_1(g , y) = \whalpha_1(g_{i_1} \,  g', y)   =    \whalpha_W(g_{i_1} ,  g' \cdot y) \circ  \whalpha_W(g' , y) } \label{eq-computation}\\
& = &       \whalpha_W(g_{i_1} ,  g' \cdot z) \circ  \whalpha_W(g' , z) 
 =   \whalpha_W(g_{i_1} \,  g', z)
 = \whalpha_1(g , z) \ . \nonumber
\end{eqnarray}
As the values of $\whalpha_W(g , y)  \in H'_1$ and  $\whalpha_W(g , z) \in H'_1$ are defined using the identity  \eqref{eq-conjmap1}, the identity \eqref{eq-computation} holds   for all $y,z \in U_1$ as   the closure   $\overline{\cZ \cap U_1} = U_1$.

Thus, the     calculation \eqref{eq-computation} shows that the value of  $\whalpha_1(g , y) \in H'_1$ 
does not depend on the choice of $y \in U_1$. We are given that $\Phi(g)(y) \in U_1$ for all $g \in G_1$, as $U_1$ is adapted to the action $\Phi$.
Thus, for $g, g' \in G_1$ and $y \in U_1$, we have $z = g' \cdot y \in U_1$ and so 
\begin{equation}\label{eq-constantalpha}
\whalpha_1(g \ g' , y) = \whalpha_1(g   , z) \ \whalpha_1(g'   ,   y) = \whalpha_1(g   , y) \ \whalpha_1(g'   ,   y)  \ .
\end{equation}
Fix $y \in U_1$ and define  $\ds  \whtheta_1 \colon G_1 \to H'_1$ by setting  $\whtheta_1(g) = \whalpha_1(g  , y)$, then \eqref{eq-constantalpha} says that $\whtheta_1$ is a group homomorphism. 
   Then \eqref{eq-conjugatequotient} is just a restatement of the identity \eqref{eq-conjmap1}.
     \endproof

Recall that $H_1 = \Phi_{1}(G_{1}) \subset \Homeo(U_1)$, and $H'_1 = \Psi(G'_{V_1}) \subset \Homeo(V_1)$.

 \begin{cor}\label{cor-inducedhomo}
   $  \whtheta_1$  induces a homomorphism  $\theta_1 \colon H_1 \to H'_1$. 
 \end{cor}
 \proof
 Suppose  that   $g, g' \in G_1$ satisfy $\Phi(g)|U_1 = \Phi(g')|U_1 \in H_1 \subset \Homeo(U_1)$. Then 
  by  the defining identity  \eqref{eq-conjmap1}, for $y \in \cZ \cap U_1$ the action of $H'_1$ is free on the orbit of $y$, hence   
  $\whtheta_1(g) = \whalpha_1(g,y) = \whalpha_1(g',y) = \whtheta_1(g') \in  H'_1$.  Thus,   
 $  \whtheta_1$  induces   a homomorphism  $\theta_1 \colon H_1 \to H'_1$. 
 \endproof

\subsection{The induced map}\label{subsec-induced2}
Recall that we are assuming the hypotheses of Theorem~\ref{thm-coestable}, so that 
  $h \colon \fX \to \fX'$ is  a continuous orbit equivalence between   Cantor actions $(\fX,G,\Phi)$  and $(\fX', G', \Psi)$, that  $G$  is finitely generated,  and  $(\fX', G', \Psi)$ is stable. 
  
  The conclusion of Corollary~\ref{cor-inducedhomo} states that the identity \eqref{eq-conjmap1}, which defines the orbit equivalence between   actions 
 $(\fX,G,\Phi)$ and $(\fX,G',\Psi)$,   induces a   homomorphism  $\theta_1 \colon H_1 \to H'_1$    but does not show that  $\theta_1$ is an isomorphism.
 We   next  show  that the image of $\theta_1$ is ``large enough'' in the profinite topology defined by adapted bases for the actions, so that $\theta_1$ induces a map between well-chosen odometer models for the actions $\Phi$ and $\Psi$. 
   The idea is to recursively construct   adapted neighborhood bases for these actions so that we can recursively apply the   techniques of the proof of Proposition~\ref{prop-coboundary}.  First, set $U_0 = V_0 = \fX$, and let $U_1 \subset V_1$ be as in the proofs of Lemma~\ref{lem-cocycle} and Proposition~\ref{prop-coboundary}. 
Then recursively choose:
\begin{enumerate}
\item $V_{\ell+1}$   adapted to $\Psi$,  with $x \in V_{\ell+1} \subset U_{\ell}$  for $\ell \geq 1$;
\item $U_{\ell}$   adapted to $\Phi$,  with $x \in U_{\ell} \subset V_{\ell}$  for $\ell \geq 2$;
\end{enumerate}
such that  $\cap_{\ell \geq 1} \ U_{\ell}   =   \cap_{\ell \geq 1} \ V_{\ell}   =   \{x\}$.
We   obtain   adapted neighborhood bases  $\cU = \{U_{\ell} \subset \fX  \mid \ell \geq 0\}$ for $(\fX,G,\Phi)$ at $x$, 
  and   $\cV = \{V_{\ell} \subset \fX  \mid \ell \geq 0\}$ for $(\fX, G', \Psi)$ at $x$.  
  
  Let $G_{\ell} \equiv G_{U_{\ell}} \subset G$ for $\ell \geq 0$, then $\cG_{\cU} = \{G_{\ell} \mid \ell \geq 0\}$ is a group chain associated to the action $\Phi$.

Let $G'_{\ell} \equiv G'_{V_{\ell}} \subset G'$  for $\ell \geq 0$,  then $\cG'_{\cV} = \{G'_{\ell} \mid \ell \geq 0\}$ is a group chain associated to the action $\Psi$.  

Set $H_{\ell} \equiv H_{U_{\ell}} = \Phi(G_{\ell})|U_{\ell} \subset \Homeo(U_{\ell})$, 
and $H'_{\ell} \equiv H'_{V_{\ell}} = \Psi(G'_{\ell})|V_{\ell} \subset \Homeo(V_{\ell})$.

  Let  $(U_1, G_1, \Phi_1)$ denote the restricted action $\Phi_1 \colon G_1 \times U_1 \to U_1$   of   $\Phi$. Let $\oPhi_1 \colon H_1 \times U_1 \to U_1$ be the induced action by $H_1 \equiv \Phi_1(G_1) \subset \Homeo(U_1)$.
  
  Let $(V_1, G'_1 , \Psi_1)$ denote the restricted action $\Psi_1 \colon G'_1 \times V_1 \to V_1$   of  $\Psi$.       Let $\oPsi_1 \colon H'_1 \times V_1 \to V_1$ be the induced action by $H'_1 \equiv \Psi_1(G'_1) \subset \Homeo(V_1)$.

We next use the existence of the   map   $\beta$     defined by the   identity   \eqref{eq-conjmap2} to prove the following result, which is used to show that   the map $\theta_1$ induces a map between the odometer models for the induced Cantor actions $(U_1,H_1,\oPhi_1)$ and  $(V_1,H'_1,\oPsi_1)$. Note that it  is not assumed that the action of $H_1$ on $U_1$ is topologically free, 
and as a consequence the choice of the map $\beta$ in \eqref{eq-conjmap2} need not be unique.  Correspondingly, the map $\theta_1$ need not be a monomorphism.

  \begin{lemma}\label{lem-containment}
 For all $\ell \geq 1$,       the subgroup $H'_{\ell} \subset H'_1$ defined as follows,  satisfies
  \begin{equation}\label{eq-filterpreserving}
H'_{\ell+1} \equiv   \Psi_1(G'_{\ell+1})  \subset \whtheta_1(G_{\ell}) = \theta_1 \circ \Phi_1(G_{\ell}) \subset  \Psi_1(G'_{\ell}) \equiv H'_{\ell} \subset H'_1 ~ .
\end{equation}
  \end{lemma}
  \proof
Fix $\ell \geq 1$,     let $k \in G'_{\ell+1}$ and $y \in  V_{\ell+1}$, then by the identity \eqref{eq-conjmap2} there exists $g = \beta(k,y) \in G$ and an open set $V_{k,y} \subset \fX$,  such that   $\Phi(g)(z) = \Psi(k)(z)$ for all $z \in V_{k,y}$.  Note that we are not assuming that the function $\beta$ is a cocycle, just the fact that it exists. 
Then by the definition of $\whalpha_1 \colon G_1 \times U_1 \to H'_1$    as the restriction of  $\whalpha_W$,  
we have that $\whalpha_1(g) = \Psi(k)$ when restricted to $V_{k,y}$.

Recall that the action of  $H'_1 = \Psi_1(G'_{V_1})$ on $V_1$ is topologically free, and $H'_1$ acts freely on the dense $H'_1$-invariant subset $\cZ \subset V_1$. 
The subset $V_{k,y} \cap V_1$   is open in $V_1$ and contained in the closure of $\cZ \cap V_{k,y} \cap V_1$, 
hence the   maps   $\Phi(g)$ and $\Psi(k)$ agree on an open subset of $V_1$, hence agree on $V_1$ since the homomorphism $\theta_1$ does not depend on the choice of $U_1$.    

It is given that $k \in G'_{\ell+1}$ and $y \in  V_{\ell+1} \subset U_{\ell}$,  
so  $\Psi(k)(y) \in V_{\ell+1} \subset U_{\ell}$, 
and thus  $\Phi(g)(U_{\ell}) \cap U_{\ell} \ne \emptyset$, hence   $g \in G_{\ell}$.  
 Thus, by \eqref{eq-conjugatequotient}   we have
 \begin{equation}
\Psi(k)|V_1 = \Phi(g)|V_1 = \whtheta_1(g)|V_1 \in \whtheta_1(G_{\ell}) ~ ,   
\end{equation}
which establishes    the first inclusion in    \eqref{eq-filterpreserving}.

To show the second inclusion in  \eqref{eq-filterpreserving}, let $g \in G_{\ell}$ then  note that $\Phi(g)(U_{\ell}) \cap U_{\ell} \ne \emptyset$.
Then by \eqref{eq-conjugatequotient} we have $\whtheta_1(g) \cdot U_{\ell}  \cap U_{\ell} \ne \emptyset$, 
hence $\whtheta_1(g) \cdot V_{\ell} \cap V_{\ell} \ne \emptyset$  and so $\whtheta_1(g) \in H'_{\ell}$.
 \endproof

Next, we  examine the relation between the 
discriminant groups for the Cantor actions  $(\fX,\oG,\whPhi)$ and   $(U_k,H_k,\oPhi_k)$, for $k \geq 1$, 
where $\cU = \{U_{\ell} \subset \fX  \mid \ell \geq 0\}$ is the adapted neighborhood basis at $x$ chosen above for  $(\fX,G,\Phi)$. 
 Let $\whG_{\infty}$ be the profinite group associated to the group chain $\cG_{\cU}$, 
and let $\whPhi_x \colon  \whG_{\infty} \times X_{\infty} \to X_{\infty}$ be the associated action. 
We  show that the wild property for this action is preserved by restriction.

  By Theorem~\ref{thm-quotientspace}, there is a   topological isomorphism  $\whTheta \colon    \oG \to \whG_{\infty}$, which conjugates the action $\whPhi$  to the action $\whPhi_x$ so we may make all calculations in terms of the action $\whPhi_x \colon \whG_{\infty} \times X_{\infty} \to X_{\infty}$.

Recall that for $\ell \geq 1$, the clopen set  $\whU_{\ell} \equiv \whG_{\ell}/\cD_x \subset \whG_{\infty}/\cD_x \cong X_{\infty}$    is adapted to the action $\whPhi_x$, with  stabilizer    $\whG_{\ell} \subset \whG_{\infty}$ defined by   \eqref{eq-whG}. Recall that $\Phi_{\ell} \colon G_{\ell} \times \whU_{\ell} \to \whU_{\ell}$ denotes the restriction of the action $\Phi_x$ to $\whU_{\ell}$. 
 Now for    $\ell' \geq \ell \geq 0$, define
 \begin{equation}\label{eq-restrictedmonodromies}
H_{\ell, \ell'} \equiv \Phi_{\ell}(G_{\ell'}) \subset \Phi_{\ell}(G_{\ell}) \subset \Homeo(\whU_{\ell}) \ ,
\end{equation}
  and let  $\whH_{\ell,\ell'}$ denote the closure of $H_{\ell, \ell'}$  in the uniform topology on $\Homeo(\whU_{\ell})$. Note that 
$H_{\ell, \ell} = H_{\ell}$ and $\whH_{\ell,\ell} = \whH_{\ell}$ as defined previously. Then the restricted action $\Phi_{\ell}$ induces  a   homomorphism   $\whPhi_{\ell} \colon \whG_{\ell} \to \whH_{\ell}$ which is a surjection, as seen using the   same method of proof as for Lemma~\ref{lem-closures}.

 Recall that  for all $\ell \geq 0$,  $\cD_x \subset   \whG_{\ell}$   so there is a well defined isotropy action $\rho_{\ell} \colon \cD_x \to \Homeo(\whU_{\ell})$. Then by \eqref{eq-isotropykernel} we have $K_{\ell} = \ker \ \{\rho_{\ell} \colon \cD_x \to \Homeo(\whU_{\ell}) \}$.
 Thus,  $K_{\ell}$ is   a normal subgroup of $\cD_x$ and in particular, $K_{\ell}$ is normal in $K_{\ell'}$  for $\ell' > \ell$.
 
For   $k \geq 0$, by   Lemma~\ref{lem-intersections} we have   
$ \cD_x  = \bigcap_{\ell\geq 0} \  \whG_{\ell} = \bigcap_{\ell\geq k} \  \whG_{\ell}$.

 For $k \geq 0$, recall that   $H_{k} = \Phi_{k}(G_{k}) \subset \Homeo(\whU_{k})$, 
 and let $\oPhi_k \colon H_k \times \whU_k \to \whU_k$ 
denote the  quotient action induced   from the restricted action $(\whU_{k}, G_{k}, \Phi_{k})$. Then we have:
 
  \begin{prop}\label{prop-quotientD}
  For $k \geq 0$, 
  the stabilizer direct limit group $\Upsilon_s(\oPhi_{k})$ for   the induced action
$(\whU_{k}, H_{k}, \oPhi_{k})$ is represented by the quotient group chain $\{K_{\ell}/K_{k} \mid \ell \geq k \}$.
   \end{prop}
  \proof
  First, observe that $\cH_{k} \equiv \{H_{k, \ell} \mid \ell \geq k\}$ is an adapted group chain for the action $(\whU_{k}, H_{k}, \oPhi_{k})$ at $x \in \whU_{k}$.  
  Recall that $\whH_{k} \subset \Homeo(\whU_k)$ is the closure of the group $H_k$. 
  
  Let $A_{k} = \ker \ \{\whPhi_{k} \colon \whG_{k} \to \whH_{k}\}$, which is a closed normal subgroup of $\whG_{k}$. 
  Then for $\ell > k$ we have $\whH_{k, \ell} \cong \whG_{\ell}/(\whG_{\ell} \cap A_{k})$, and the 
  discriminant   $\ocD_{k}$ of the action $(\whU_{k}, H_{k}, \oPhi_{k})$ is given by $\ocD_{k}   =   \bigcap_{\ell \geq k} \ \whH_{k, \ell}$.

  Suppose $\whh \in \cD_x$ then the image $\oPhi_{k}(\whh) \in \whH_{k, \ell}$ for all $\ell \geq k$, hence $\whPhi_{k}(\whh) \in \ocD_{k}$. 
   Moreover, as $\whG_{k}$ is sequentially compact, the restricted map $\whPhi_{k} \colon \cD_x \to \ocD_{k}$ is onto, 
  with kernel the closed normal subgroup $K_k \equiv \cD_x \cap A_{k}$. 
Thus,   there is  an   isomorphism   $\whPhi_{k} \colon \cD_x/K_{k}  \cong \ocD_{k}$.

For each $\ell > k$,   let $\orho_{k, \ell} \colon \ocD_{k} \to \whH_{\ell}$ be the restriction of the map $\oPhi_{k} \colon \whH_{k, \ell} \to \whH_{\ell} \subset \Homeo(\whU_{\ell})$, and let  $\oK_{k, \ell} = \ker \ \{\oPhi_{k} \}$. Then the stabilizer  group $\Upsilon_s(\oPhi_{k})$ is represented by the chain $\{\oK_{k, \ell} \mid \ell \geq k\}$.

Let $\whk \in \oK_{k, \ell}$ and  so $\whk$ acts as the identity on $\whU_{\ell}$. Choose $\whh \in \cD_x$ with $\whPhi_{k}(\whh) = \whk$,  then $\whh$ also acts as the identity on $\whU_{\ell}$, hence $\whh \in K_{\ell}$. Moreover, if $\whh \in K_{\ell}$ and $\whPhi_{k}(\whh) = id$ then $\whh  \in K_{\ell} \cap A_{k} = K_{k}$. It follows that $\oK_{k, \ell} \cong K_{\ell}/K_{k}$ for $\ell \geq k$, which yields the claim.
  \endproof

 \begin{cor}\label{cor-wildwild}
For $k \geq 1$,    $(\fX,G,\Phi)$ is  wild if and only if    $(\whU_{k},H_{k},\oPhi_{k})$ is   wild. 
\end{cor}
\proof
 Let   $\ell \geq 0$. Then $(\fX,G,\Phi)$ is a wild action, if and only if   the chain of subgroups $\{K_{\ell} \mid \ell \geq 0\}$ as above is unbounded, if and only if  the chain of quotient groups $\{K_{\ell}/K_{k} \mid \ell \geq k\}$ is   unbounded. Then by the proof  of Proposition~\ref{prop-quotientD} above,  we have $\oK_{k, \ell} \cong K_{\ell}/K_{k}$, so the action is wild if and only if  the    chain $\{\oK_{k, \ell} \mid \ell \geq k\}$ is unbounded, if and only if   the stabilizer group $\Upsilon_s(\oPhi_{k})$ is unbounded.
\endproof

  We can now complete the proof of Theorem~\ref{thm-coestable}.  We are given the Cantor actions $(\fX,G,\Phi)$  and $(\fX, G', \Psi)$, where $G$ and $G'$ are finitely generated groups,  $\Phi$ and $\Psi$  have the same orbits, and satisfy the conditions \eqref{eq-conjmap1} and \eqref{eq-conjmap2}.
We assume that the action   $(\fX,G,\Phi)$ is wild, and derive a contradiction to the assumption that the action $(\fX', G', \Psi)$ is stable.  

By Proposition~\ref{prop-coboundary}, there is a clopen set $U_1 \subset \fX$ adapted to the action $\Phi$, a clopen set $V_1 \subset \fX$ adapted to the action $\Psi$, and a homomorphism   $\theta_1 \colon H_1 \to H'_1$ which satisfies the condition \eqref{eq-conjugatequotient}. 
 
Let   $\cU = \{U_{\ell} \subset \fX  \mid \ell \geq 0\}$ be an  adapted neighborhood basis for $(\fX,G,\Phi)$ at $x$,    and let   $\cV = \{V_{\ell} \subset \fX  \mid \ell \geq 0\}$ be an   adapted neighborhood basis for $(\fX, G', \Psi)$ at $x$,  where $U_1$ and $V_1$ are chosen so that 
Proposition~\ref{prop-coboundary} holds, and the chains $\cU$ and $\cV$ chosen as above so that   the conclusions of Lemma~\ref{lem-containment} hold. We use the notation of the proofs of these results in the following as well.
 
 Let $\cG_{\cU} = \{G_{\ell} \mid \ell \geq 0\}$ be the   group chain associated to $\cU$, and let $\cH_{\cU} = \{H_{1,\ell} \mid \ell \geq 1\}$ be the group chain for the induced action $(U_1, H_1, \oPhi_1)$.

Let  $\whG_{\infty}$ be the profinite group associated to the group chain $\cG_{\cU}$, with discriminant group $\cD_x \subset \whG_{\infty}$.

Let $\whH_{\infty}$ be the profinite group associated to the group chain $\cH_{\cU}$, with discriminant group $\ocD_1 \subset \whH_{\infty}$.

Let $\cK = \{K_{\ell} \mid \ell \geq 0\}$ be the   group chain representing $\Upsilon_s^x(\Phi)$, and let $\ocK \cong \{\oK_{1,\ell} \cong K_{\ell}/K_1 \mid \ell \geq 1\}$ be the   group chain representing $\Upsilon_s^x(\oPhi_1)$. Then by   Corollary~\ref{cor-wildwild}, the chain $\ocK$ is unbounded.

 Next, let  $\cG'_{\cV} = \{G'_{\ell} \mid \ell \geq 0\}$ be the   group chain associated to $\cV$, and let $\cH'_{\cV} = \{H'_{1,\ell} \mid \ell \geq 1\}$ be the group chain for the induced action $(V_1, H'_1, \oPsi_1)$. Note that $H'_{1,1} = H'_1$. 
 
Let  $\whG'_{\infty}$ be the profinite group associated to the group chain $\cG'_{\cV}$, 
with discriminant group $\cD'_x \subset \whG'_{\infty}$.

Let $\whH'_{\infty}$ be the profinite group associated to the group chain $\cH'_{\cV}$, with discriminant group $\ocD'_1 \subset \whH'_{\infty}$.

Let $\cK' = \{K'_{\ell} \mid \ell \geq 0\}$ be the   group chain representing $\Upsilon_s^x(\Psi)$, and let $\ocK' = \{\oK'_{1,\ell} \cong K'_{\ell}/K'_1 \mid \ell \geq 1\}$ be the   group chain representing $\Upsilon_s^x(\oPsi_1)$.

The homomorphism   $\theta_1 \colon H_1 \to H'_1$   satisfies the condition \eqref{eq-filterpreserving},   which implies that 
\begin{equation}\label{eq-nested}
H'_{1,\ell +1} \subset \theta_1(H_{1,\ell}) \subset H'_{1,\ell}
\end{equation}
so induces a map $\whtau_1 \colon \whH_{\cU} \to \whH'_{\cV}$. Then \eqref{eq-nested} and Lemma~\ref{lem-intersections} imply that 
$\whtau_1$ restricts to a surjection $\whtau_1 \colon \ocD_1 \to \ocD'_1$ which satisfies $\whtau_1(\oK_{\ell}) = \oK'_{\ell}$. 
As a result, by Theorem~\ref{thm-fundamental} we obtain maps
\begin{equation}\label{eq-inducedmap2}
  \fX \supset U_1 \cong \whH_{\infty}/\ocD_1 \stackrel{\whtau_*}{\longrightarrow} \whH'_{\infty}/\ocD'_1 \cong V_1 \subset \fX   ~ .
\end{equation}
That is, $\whtau_* \colon U_1 \to V_1$ is the map induced from the homomorphism $\whtau_1$ which is $H_1$-equivariant.

We claim  that if $\whk \in \oK_{\ell +1}$ and $\whk \not\in \oK_{\ell}$, then $\whtau_1(\whk) \not\in \oK'_{\ell-1}$, hence 
$\oK'_{\ell-1} \ne \oK'_{\ell+1}$. 
Let $\whh \in \cD_x$ be chosen so that it has image $\whk \in \ocD_1$. Then $\whk \not\in \oK_{\ell}$ implies that the restricted action of $\Phi_{\ell}(\whh) \colon \whU_{\ell} \to \whU_{\ell}$ is not the identity, while the restricted action $\Phi_{\ell+1}(\whh) \colon \whU_{\ell+1} \to \whU_{\ell+1}$ is the identity.

By the definition of the conjugating map $\theta_1$ in Proposition~\ref{prop-coboundary}, and more precisely the relation \eqref{eq-conjugatequotient},   the map $\theta_1(\Phi_{1}(\whh)) \in \Homeo(\whV_1)$ restricts to $\Phi_{\ell}(\whh)$  on $\whU_{\ell}$, and so is not the identity on $\whU_{\ell}$. Since $\whU_{\ell} \subset \whV_{\ell-1}$ by construction, this implies that $\theta_{1}(\Phi_{2}(\whh))|\whV_{\ell-1} \in \Homeo(\whV_{\ell-1})$ is also not the identity, hence  $\whtau_1(\whk) \not\in \oK'_{\ell-1}$. 

Thus, if the chain $\cK = \{K_{\ell} \mid \ell \geq 0\}$ is unbounded, then the chain $\ocK' = \{\oK'_{\ell}   \mid \ell \geq 1\}$ is unbounded, and hence the action is $(V_1, H'_1, \oPsi_1)$ is wild. This  contradicts the choice of $V_1$, and    completes the proof of Theorem~\ref{thm-coestable}.
 \endproof
 
 \begin{remark}\label{rmk-autoex}
 {\rm
 Here is a simple example that illustrates one of the main ideas in the above   proof of Theorem~\ref{thm-coestable}. 
 Let $\Phi \colon \mZ^n \times \fX \to \fX$ be a faithful action, for $n \geq 2$. So $(\fX, \mZ^n, \Phi)$ is an odometer action of the abelian group $\mZ^n$ and is a free action. Let $\theta \in {\bf SL(n,\mZ)}$ be a non-identity matrix, which acts via conjugation on $\mZ^n$. Let $\Psi = \Phi^{\theta}$ be the   action of $\mZ^n$ on $\fX$ as defined by \eqref{def-thetaconjugacy}, which is again free. Then the identity map is a continuous orbit equivalence between $(\fX, \mZ^n, \Phi)$ and $(\fX , \mZ^n, \Psi)$. The map $\theta_1$ obtained in Proposition~\ref{prop-coboundary} is then just the restriction of $\theta$ to a subgroup of finite index in $\mZ^n$. For a choice of group chains adapted to the actions, the   map $\theta_1$ need not preserve the subgroups in these chains, but the map $\whtau_*$ in \eqref{eq-inducedmap2} is the map on the inverse limit space $\fX$ induced by the action of  $\Theta$. The proof above  is a generalization of this example to the more general (and more complicated)  non-abelian setting.

 }
 \end{remark}

 Finally, we give an application  of Theorem~\ref{thm-coestable} and the techniques used in its proof, to obtain   an extension of Theorem 3.3 by Cortez and Medynets in \cite{CortezMedynets2016}, and  Theorem~1.5 in \cite{HL2018b}.

   \begin{thm}\label{thm-coe=re}
   Let $G$ and $G'$ be   finitely generated groups, and suppose that the Cantor action   $(\fX', G', \Psi)$ is stable. Let $(\fX, G, \Phi)$ be a Cantor action   which is  continuously orbit equivalent to $(\fX', G', \Psi)$,  then the actions are return equivalent. 
   \end{thm}
   \proof
   By Theorem~\ref{thm-coestable}, the action $(\fX, G, \Phi)$  is stable. Thus, in the proof of Theorem~\ref{thm-coestable}, we can choose the clopen set $U_1 \subset \fX$ so that the restricted action $\oPhi_1 \colon H_1 \times U_1 \to U_1$ is topologically free.
   
  Then in the proof of Lemma~\ref{lem-containment},  the choice of 
   $g = \beta(k,y) \in G$   such that   $\Phi(g)(z) = \Psi(k)(z)$ is unique on $U_1$, so the map $\theta_1$ is injective. That is, the 
   homomorphism   $\theta_1 \colon H_1 \to H'_{V_1}$ is an isomorphism onto its image. Set $H''_1 = \theta_1(H_1) \subset H'_{V_1}$.
   
   Then \eqref{eq-nested} implies that the group chain $\cH'_{\cV} = \{H'_{1,\ell} \mid \ell \geq 1\}$   for the induced action $(V_1, H'_1, \oPsi_1)$,  is equivalent  in the sense of 
   Definition~\ref{defn-greq}, to the group chain $\cH'' = \{H''_{\ell} \mid \ell \geq 1\}$ where $H''_1 = H'_{V_1}$, and 
   $H''_{\ell} = \theta_1(H_{1,\ell})$ for $\ell \geq 1$. Let $W'_1 = \whtau_*(U_1) \subset V'_1$ be the clopen subset corresponding to the truncated chain   $\{H'_{1,\ell} \mid \ell \geq 1\}$, where $\whtau_*$ is the map in \eqref{eq-inducedmap2}.  Let $h_1 \colon U_1 \to W'_1$ denote  the map $\whtau_*$ onto its image, which then  conjugates the induced action   $(U_1, H_1, \oPhi_1)$ with the Cantor action 
     $(W'_1, H''_{\ell}, \oPsi_1)$ which is induced from $(\fX', G', \Psi)$. Thus, the actions $\Phi$ and $\Psi$ are return equivalent. 
   \endproof

 \section{Wild actions and non-Hausdorff elements}\label{sec-hausdorff} 
  
  In this section, we consider Cantor actions $(\fX,G,\Phi)$ for which   the stabilizer direct limit group  $\Upsilon_s(\Phi)$ is  unbounded, and derive some of their dynamical properties. Of special interest will be the existence of \emph{non-Hausdorff elements} for the action, which are defined in Section~\ref{subsec-groupoid}. 
  
  The notion of a non-Hausdorff element is local in $\fX$, and our first result   is a generalization  of  Proposition~3.1 by Renault in \cite{Renault2008}.   Proposition~\ref{prop-hausdorff2} below implies that   a Cantor action with a non-Hausdorff element has unbounded stabilizer group $\Upsilon_s(\Phi)$.  Then in Theorem~\ref{thm-distinct=nH}, we show that a Cantor action with a non-Hausdorff element must be dynamically wild.  
  The converse implication, which is to give criteria for the existence of non-Hausdorff elements, is a much more subtle problem.    

  \subsection{The germinal groupoid}\label{subsec-groupoid}
The reduced  $C^*$-algebra $C_r^*(\fX,G,\Phi)$ associated to a Cantor action  $(\fX,G,\Phi)$ is an invariant of the continuous orbit equivalence class of the action. The study of the $K$-theory of  $C_r^*(\fX,G,\Phi)$ offers another approach to the classification of Cantor actions, as used for example in the work \cite{GPS2019}. This $C^*$-algebra can also be constructed using the  germinal groupoid $\cG(\fX, G, \Phi)$ associated to the action,    as discussed for example by Renault in \cite{Renault1980}.  
Then one can ask how the    properties of the $C^*$-algebra  are related to the dynamical properties of the action, as discussed by Renault in  \cite{Renault2006}. 

In the work \cite{Renault2008}, Renault assumes that the Cantor action is topologically free, and thus the germinal groupoid $\cG(\fX, G, \Phi)$ is a Hausdorff topological space, in order to avoid technical difficulties that arise otherwise. In fact, as remarked in Corollary~\ref{cor-Hgroupoid}, there are also wild actions for which $\cG(\fX, G, \Phi)$ is still Hausdorff. Thus,   the case where $\cG(\fX, G, \Phi)$  has non-Hausdorff topology may be considered to be  exceptional, and the fact that the topology is non-Hausdorff    has implications for the algebraic structure of   $C_r^*(\fX,G,\Phi)$, as discussed   in  \cite{BCFS2014,Exel2011}. 

 Recall first the definition of the germinal groupoid $\cG(\fX, G, \Phi)$.
For $g_1, g_2 \in G$, we say that $\Phi(g_1)$ and $\Phi(g_2)$ are \emph{germinally equivalent} at $x \in \fX$ if $\Phi(g_1) (x) = \Phi(g_2)(x)$, and there exists an open neighborhood $x \in U \subset \fX$ such that the restrictions agree, $\Phi(g_1)|U = \Phi(g_2)|U$. We then write $\Phi(g_1) \sim_x \Phi(g_2)$. For $g \in G$, denote the equivalence class of  $\Phi(g)$ at $x$ by $[g]_x$.
 The collection of germs $\cG(\fX, G, \Phi) = \{ [g]_x \mid g \in G ~ , ~ x \in \fX\}$ is given the sheaf topology, and   forms an \emph{\'etale groupoid} modeled on $\fX$. 
 We recall the following formulation of the Hausdorff property that was given by  Winkelnkemper:
 \begin{prop}\cite[Proposition 2.1]{Winkelnkemper1983}\label{prop-hausdorff}
The germinal groupoid $\cG(\fX, G, \Phi)$ is Hausdorff at $[g]_x$  if and only if, for all $[g']_x \in \cG(\fX, G, \Phi)$ with $g \cdot x = g' \cdot x$, if there exists a sequence   $\{x_i\} \subset \fX$ which converges to $x$ such that $[g]_{x_i} = [g']_{x_i}$ for all $i$, then $[g]_{x} = [g']_{x}$. 
\end{prop}

  Winkelnkemper showed in \cite[Proposition~2.3]{Winkelnkemper1983} that for a smooth foliation $\F$ of a connected manifold $M$ for which the associated holonomy pseudogroup $\cGF$ is generated by real analytic maps, then $\cGF$ is a Hausdorff space. For Cantor actions, an analogous result holds for the LQA property. 

\begin{prop}\label{prop-hausdorff2}
If an action $(\fX,G,\Phi)$   is   locally quasi-analytic, then  $\cG(\fX, G, \Phi)$ is    Hausdorff.     
 \end{prop}
 \proof
Assume that $\cG(\fX, G, \Phi)$ is  not  Hausdorff.  Then there exists $g \in G$ and $x \in \fX$ such that   $\cG(\fX, G, \Phi)$ is non-Hausdorff at $[g]_x$. By Proposition~\ref{prop-hausdorff},  there exists $[g']_x \in \cG(\fX, G, \Phi)$ with $g \cdot x = g' \cdot x$, and a sequence   $\{x_i\} \subset \fX$ which converges to $x$ such that $[g]_{x_i} = [g']_{x_i}$ for all $i$, but $[g]_x \ne [g']_x$. Let  $g'' = g^{-1} g' \in G$,  then $g'' \cdot x = x$ and $[g'']_{x_i} = [id]_{x_i}$ for all $i$, but $[g'']_{x} \ne [id]_{x}$. 

The fact that the germ  $[g'']_{x} \ne [id]_{x}$ means that the action of $\Phi(g'')$ is not the identity in any open neighborhood of $x$.
On the other hand, the germinal equalities $[g]_{x_i} = [g']_{x_i}$ for $i \geq 1$ imply there  exists a sequence of open sets $x_i \in W_i \subset \fX$ for which the restriction of $\Phi(g'')$  to $W_i$ is the identity. Hence, there does not exists $\e >0$ such that  $\Phi(g)|U$ is   quasi-analytic  for any open neighborhood $x \in U$ with $\diam \ (U) < \e$. Thus, the action $(\fX,G,\Phi)$   is not  locally quasi-analytic.
\endproof

 \subsection{Non-Hausdorff elements in the closure}\label{subsec-nHclosure}

Recall that by  Theorem~\ref{thm-boundedLCQA}, a Cantor action   $(\fX,G,\Phi)$  is a locally completely quasi-analytic (LCQA) action if and only if the stabilizer group  $\Upsilon_s(\Phi)$ is   bounded.  Based on Proposition~\ref{prop-hausdorff}, we introduce the following notion:

\begin{defn}\label{def-nHelement}
 Let $(\fX,G,\Phi)$  be a  Cantor action, and let $\oG \subset \Homeo(\fX)$ denote the closure of the action. Then $\whg \in \oG$ is a \emph{non-Hausdorff element}  at $x \in \fX$  if: 
 \begin{enumerate}
\item $\whg \cdot x = x$; and there exists
\item a sequence   $\{x_i \mid i \geq 1\} \subset \fX$   converging to $x$; and 
\item clopen subsets $x_i \in W_i \subset  \fX$;
\end{enumerate}
such that for any clopen subset $x \in U \subset \fX$ the restriction of $\whg$ to $U$ is not the identity, while  for all $i \geq 1$, we have $\whg \cdot W_i = W_i$ and the restriction of $\whg$   to $W_i$ is the identity.
 \end{defn}

We then have the following consequence of Proposition~\ref{prop-hausdorff2}, which gives a connection between the unbounded property for $\Upsilon_s(\Phi)$  and the dynamics of the action. 
 \begin{cor}\label{cor-nH=unbounded}
  Let $(\fX,G,\Phi)$  be a  Cantor action. Suppose that $\oG$ contains a non-Hausdorff element, then the    stabilizer limit group   $\Upsilon_s(\Phi)$ is unbounded. That is, the action is wild.
    \end{cor}
 \proof
Assume that  $\oG$ contains a non-Hausdorff element,   then by Proposition \ref{prop-hausdorff2} applied to the action of $\oG$,   the action $(\fX,G,\Phi)$   is  not locally completely quasi-analytic. Then by Theorem \ref{thm-boundedLCQA} $\Upsilon_s(\Phi)$ is unbounded. 
 \endproof

 \subsection{Non-Hausdorff dynamics}\label{subsec-ucg}
 
We next give in  Theorem~\ref{thm-distinct=nH}  a sharper version of the conclusion of Corollary~\ref{cor-nH=unbounded}. As an application, 
  Corollary~\ref{cor-producttype}     implies that  the examples of wild actions constructed in  \cite[Section~9]{HL2018a}   cannot have non-Hausdorff elements in $\oG$. The arboreal constructions in Section~\ref{sec-arboreal}  yield examples of Cantor actions with non-Hausdorff elements in $\Phi(G)$.

    The difficulty in constructing non-Hausdorff elements can be seen from the following considerations.   
    The hypothesis that  $\Upsilon_s(\Phi)$ is unbounded implies that for $x \in \fX$, there exists a sequence of adapted clopen neighborhoods of   $x$ and elements in $\oG_x$ which do not act as the identity on this open neighborhood, but do  act as the identity  on some smaller clopen neighborhoods. On the other hand, the non-Hausdorff property asserts there is some fixed element $\whh \in \oG_x$ so that for any adapted neighborhood  with $x \in U$, the restricted action of $\Phi(\whh)$ is not the identity on $U$, but the action is the identity on some clopen subset $V \subset U$.  The distinction is that the non-Hausdorff condition in Proposition~\ref{prop-hausdorff} is a statement about the local action of a fixed element $\whh$, while the wild hypothesis is a statement about the behavior of a sequence of elements in $\oG$.      We use the interplay of these two notions in the proof of Theorem~\ref{thm-distinct=nH}.
  
Recall from Definition~\ref{def-finitetype} that the limit group $\Upsilon_s(\Phi)$ for a Cantor action $(\fX,G,\Phi)$ has finite type if it is represented by an increasing  chain $\cK = \{K_{\ell} \mid \ell \geq 0\}$ where each $K_{\ell}$ is a finite group.  Recall from Definition \ref{def-wildtypes} that a Cantor action $(\fX,G,\Phi)$ is dynamically wild if the stabilizer group $\Upsilon_s(\Phi)$ is unbounded, and there is a proper inclusion $\Upsilon_c(\Phi) \subset \Upsilon_s(\Phi)$.

 \begin{thm}\label{thm-distinct=nH}
  Let $(\fX,G,\Phi)$  be a  Cantor action with a non-Hausdorff element   $\whh \in \oG_x$,  then 
    $\Upsilon_s(\Phi)$  does not have finite type, and 
  the action is dynamically wild. 
  \end{thm}
 
\proof
 Let $x \in \fX$, with   $\whg \in \oG_x$ a non-Hausdorff element at $x$. Then by Definition~\ref{def-nHelement},  for any adapted clopen set $U$ with  $x \in U$,  the restriction of $\whg$ to $U$ is not the identity. In addition, 
   there exists a sequence   $\{x_i \mid i \geq 1\} \subset \fX$   of distinct points converging to $x$, with  $\whg \cdot x_i = x_i$ for all $i \geq 1$, and   clopen subsets $x_i \in W_i \subset  \fX$ such that   $\whg \cdot W_i = W_i$   and the restriction of $\whg$   to $W_i$ is the identity. 
   It follows   that $x \not\in W_i$ for all $i \geq 1$.
 
 Let $\cU = \{U_{\ell} \subset \fX  \mid \ell \geq 0\}$ be an adapted neighborhood basis for the action $\Phi$ at $x$.

Let   $\cG^x_{\cU} = \{G_{\ell}\mid i \geq 0 \}$ be the group chain associated to   $\cU$. 
 Let $\{K_{\ell}  \mid \ell \geq 0\}$ be its chain of stabilizer groups, and       $\{Z_{\ell}  \mid \ell \geq 0\}$  be its chain of centralizer groups as in 
 Definition~\ref{def-progroupinvariants}.

Recall from  Corollary~\ref{cor-inducedquotients} that the map $\whTheta$ in Theorem~\ref{thm-quotientspace} 
induces a  homeomorphism of $G$-spaces
$\Theta_x \colon  \oG/\oG_x \cong  \fX  \to  \whG_{\infty}/\cD_x \cong X_{\infty}$.  

  For each $i \geq 1$, choose $\whg_i \in \whG_{\infty}$   such that $x_i = \whg_i \cdot x$. 
  
As $x_i \in W_i$ and $x \not\in W_i$, there exists $\ell_i > i$ such that $U_{\ell_i} \cap W_i = \emptyset$ and $\whg_i \cdot U_{\ell_i} \subset W_i$, hence   $\whg_i \cdot U_{\ell_i} \cap U_{\ell_i} = \emptyset$. Moreover, as $x_i$ limits to $x$ there exists $j > i$ such that $x_j \in U_{\ell_i}$ and $x_j \not\in W_i$.

As $\{U_{\ell} \mid \ell \geq 0\}$ is a neighborhood system about $x$, by passing to   subsequences chosen recursively, 
we can assume that for $\ell \geq 1$ we have: 
\begin{equation}\label{eq-indexing}
 U_{\ell} \cap W_{\ell} = \emptyset    ~ , ~  x_{\ell} = \whg_{\ell} \cdot x        ~ , ~ \whg_{\ell} \cdot U_{\ell} \subset W_{\ell}    ~ , ~  x_{\ell +1} \not\in W_{\ell} ~ .
\end{equation}

Set $\whh_{\ell} =  \whg_{\ell}^{-1} \  \whg \ \whg_{\ell}$.
Then $\whg \cdot x_{\ell} = x_{\ell}$ implies that $(\whg_{\ell}^{-1}  \ \whg \ \whg_{\ell}) \cdot x = x$,  hence  $\whh_{\ell} \cdot x = x$ so $\whh_{\ell} \in \cD_x$.

Since $\whg$ acts as the identity on $W_{\ell}$ and $\whg_{\ell} \cdot U_{\ell} \subset W_{\ell}$, we obtain  that $\whh_{\ell} \in K_{\ell}$.

Set  $y_{\ell} = \whg_{\ell}^{-1} \cdot x$ and observe that $\whh_{\ell} \cdot y_{\ell} = y_{\ell}$ and that $y_{\ell} \not\in U_{\ell+1}$ as $\whg_{\ell} \cdot U_{\ell+1} \cap U_{\ell+1} = \emptyset$.

We claim   this implies that  each $K_{\ell}$ is not a finite group, and hence $\Upsilon_s(\Phi)$  does not have finite type.

 For $\ell \geq 0$, recall  from  \eqref{eq-whG} that $\whG_{\ell} = \whC_{\ell} \ \cD_x$ is a clopen subset of $\whG_{\infty}$. 
Then  $\Theta_x$ identifies the clopen set  $U_{\ell} \subset \fX$ with $\whG_{\ell}/\cD_x \subset \whG_{\infty}$. 
   Thus,    $\whC_{\ell}$ acts transitively on   the quotient space 
$$\whU_{\ell} = \whG_{\ell}/\cD_x = (\whC_{\ell} \ \cD_x)/ \cD_x = \whC_{\ell}/(\whC_{\ell} \cap \cD_x) \ .$$
   It is given that $\whC_{\ell} \subset \whG_{\ell}$ is a normal subgroup of $\whG_{\infty}$, 
  hence for any $g \in G$ we have $g \ \whC_{\ell} \ g^{-1} = \whC_{\ell}$, thus $\whC_{\ell}$ also acts transitively on the clopen set $\whU_{\ell}^g \equiv g \cdot \whU_{\ell}$.

For each $\whk \in \whC_{\ell}$ define  the conjugate element $\whh_{\ell}^{\whk} = \whk \ \whh_{\ell} \ \whk^{-1}$. 
Since $\whh_{\ell}$ acts as the identity on $U_{\ell}$, and   $\whk  \cdot U_{\ell} = U_{\ell}$ we have  $\whh_{\ell}^{\whk} \in K_{\ell}$.

Finally, recall that  $\whg$ is not the identity on any clopen neighborhood of $x$, and    $\whh_{\ell} =  \whg_{\ell}^{-1} \  \whg \ \whg_{\ell}$ satisfies  $\whh_{\ell} \cdot y_{\ell} = y_{\ell}$ and is not the identity map in any open neighborhood of $y_{\ell}$. That is, the map $\whh_{\ell}$ has a non-trivial germ  at $y_{\ell}$.
 Thus, each  conjugate map $\whh_{\ell}^{\whk}$ has a fixed point at $\why = \whk \cdot y_{\ell} \in  U_{\ell}^{\whg_{\ell}^{-1}}$, and acts non-trivially on any neighborhood of $\why$.  
  Let $\cB = \{\whh_{\ell}^{\whk} \mid \whk \in  \whC_{\ell} \} \subset K_{\ell}$ be the collection of all such maps.

For $\phi \in \Homeo(\whU_{\ell}^{\whg_{\ell}^{-1}})$ let $\Fix(\phi) = \{\why \in \whU_{\ell}^{\whg_{\ell}^{-1}} \mid \phi(\why) = \why \}$, and let ${\rm Germ}(\phi , \whU_{\ell}^{\whg_{\ell}^{-1}})$ denote the subspace of $\Fix(\phi)$ such that $\phi$ has non-trivial germ at $\why$. If $\why \in W$ is a clopen subset of ${\rm Germ}(\phi , \whU_{\ell}^{\whg_{\ell}^{-1}})$, then $W \subset \Fix(\phi)$, and then the germ of $\phi$ is trivial at $\why$. This is a contradiction, so the set ${\rm Germ}(\phi , \whU_{\ell}^{\whg_{\ell}^{-1}})$ has no interior.

Now, suppose that  the collection $\cB$ is countable.  Then the set of points   $\why \in \whU_{\ell}^{\whg_{\ell}^{-1}}$  for which there exists some $\phi \in \cB$ for which $\why$ is a fixed point with non-trivial holonomy is a countable union of subsets of   $\whU_{\ell}^{\whg_{\ell}^{-1}}$ without interior. 
Since $\whU_{\ell}^{\whg_{\ell}^{-1}}$ is a Cantor space hence is Baire, by the Baire Category Theorem, this union has no interior. This contradicts the previous observation that every point of $\whU_{\ell}^{\whg_{\ell}^{-1}}$ is  a fixed point for some $\phi \in \cB$ with non-trivial germ. Thus, $\cB$ must be an uncountable collection. In particular, $K_{\ell}$ is  an uncountable group for   all $\ell \geq 0$.

Note that the above Baire argument is based on the same ideas as in the proofs of Theorem~1 by  Epstein, Millet and Tischler in \cite{EMT1977}, and Theorem~3.6 by Renault in \cite{Renault2008}.

Next, we show that the action  $(\fX,G,\Phi)$ is dynamically wild. Suppose not, then we have $\Upsilon_c(\Phi) = \Upsilon_s(\Phi)$ which implies that there exists an increasing  subsequence $\{ \ell_i \mid i \geq 1\}$ such that the inclusion $Z_{\ell_i} \subset K_{\ell_i}$ is an isomorphism for all $i \geq 1$. Then for  each $i \geq 1$ we have  the collection of maps  $\{\whh_{\ell_i}^{\whk} \mid \whk \in  \whC_{\ell_i} \} \subset Z_{\ell_i}$. For simplicity, set $\ell_i = \ell$. 

Recall that by construction for $\ell \geq 1$ the maps $h_{\ell}$ have non-trivial germinal holonomy at $y_{\ell}$. Moreover, since $\whg_{\ell}^{-1}$ conjugates $y_\ell$ to $x$, $\whh_{\ell}$ is non-Hausdorff at $y_{\ell}$. Therefore, there is a clopen subset $W' \subset \whU_{\ell}^{\whg_{\ell}^{-1}}$ such that $\whh_{\ell}$ is the identity on $W'$.

Since $\whC_\ell$ preserves and acts transitively on $\whU_{\ell}^{\whg_{\ell}^{-1}}$, there is $\whk \in \whC_\ell$ such that $\whk \cdot y_\ell \in W'$. Then the element $\whh_{\ell}^{\whk}$ fixes an open neighborhood of $y_\ell$. But this contradicts the fact that $\whh_{\ell}^{\whk} = \whk \whh_\ell \whk^{-1} = \whh_\ell$ is not the identity map on any clopen neighborhood of the  fixed point $y_{\ell}$. Thus,  $Z_{\ell} = K_{\ell}$ is impossible.
  \endproof

\begin{cor}\label{cor-producttype}
  Let $(\fX,G,\Phi)$  be a  Cantor action with discriminant group $\cD_x$ at $x \in \fX$.
Suppose that for each $\whg \in \cD_x$, the intersection $\cD_x \cap \{\whh^{-1} \ \whg \ \whh \mid \whh \in \whG_{\infty} \}$ is finite,  then the germinal groupoid of the action is Hausdorff.
\end{cor}
\proof
Suppose that $\whg \in \cD_x$ is a non-Hausdorff element. Then the proof of Theorem~\ref{thm-distinct=nH} shows that the  the set of   conjugacy classes of $\whg$ in $\cD_x$ is an infinite set.
\endproof
 
   It is a basic question to find a converse  to the conclusion of Theorem~\ref{thm-distinct=nH}. That is,  to find sufficient conditions so that  a Cantor action which is dynamically wild must have  a non-Hausdorff element. 

In the next Section~\ref{sec-arboreal}, examples of wild Cantor actions  are constructed using the   ``automata'' method, which is a well-known technique in Geometric Group Theory, and defines a homeomorphism  using a  recursive definition along the branches of a tree.  It would be very interesting to understand if the use of automata  is the only approach    to constructing actions with non-Hausdorff elements, or whether there are  possibly alternative general methods for their construction.

   \section{Examples of tree automorphisms}\label{sec-arboreal}

 In this section,  we present examples   constructed  as actions on the boundary of an infinite binary tree to   illustrate some of the properties of wild Cantor actions.     The boundary of a tree $T$ is identified with the set of all infinite paths in the tree, and can also be viewed   as a collection of   infinite sequences of $0$'s and $1$'s. It is a Cantor set by a standard argument. To construct our examples, we   use   an approach well-known in Geometric Group Theory, of defining a homeomorphism of the Cantor set at the boundary recursively, by specifying how it acts at each level of the tree. For example, the   Grigorchuk group and the Basilica group are usually defined this way; see \cite{Grigorchuk1984,Nekrashevych2005} and other works. 
 
 \subsection{Actions on trees}\label{subsec-arboreal}
We start by explaining the notation and the method for the recursive construction of a homeomorphism of the boundary of a tree. 
 
Let $T$ be a binary tree, that is, $T$ consists of the vertex set $V = \bigsqcup_{n \geq 0} V_n$ and of the set of edges $E$ with the following properties: For $n \geq 0$
  \begin{enumerate}
  \item   ${\rm card}(V_n) = 2^n$.
 \item For   each $v \in V_n$, there are exactly two vertices in $V_{n+1}$ joined to $v$ by edges.
 \item For each  $v \in V_n$, with $n > 0$, there is exactly one vertex in $V_{n-1}$ joined to $v$ by an edge.
  \end{enumerate}

We denote by $\cP(T)$ the set of all infinite connected paths in $T$.

It is sometimes convenient to label the vertices in $V$ by $0$'s and $1$'s as follows. The single vertex in $V_0$ is not labelled; the two vertices in $V_1$ are labelled one by $0$, and another one by $1$. Now suppose the vertices in $V_n$ are labelled by words of length $n$ consisting of $0$'s and $1$'s. Since there are $2^n$ such distinct words, we can assign to each vertex in $V_n$ a unique word, and to every word of length $n$ we can assign a vertex. Proceed to label the vertices in $V_{n+1}$ as follows. Let $v \in V_n$ be labelled by a word $s_1s_2\cdots s_n$. Since $T$ is binary, there are two vertices, $w_1$ and $w_2$ in $V_{n+1}$ which are joined to $v$ by edges. Label $w_1$ by $s_1s_2 \cdots s_n 0$ and $w_2$ by $s_1 s_2 \cdots s_n 1$.

It follows that every infinite sequence $s_1s_2 \cdots$, where $s_i \in \{0,1\}$ for   $ i \geq 1$, corresponds to a unique path $s$ in the tree $T$. Namely, $s$ passes through  the vertex labelled by $s_1$ in $V_1$, by $s_1s_2$ in $V_2$ and, inductively, $s$ passes through  the vertex labelled by $s_1s_2 \cdots s_n$ in $V_n$ for $n \geq 1$.

Next, let $w= w_1w_2 \cdots w_n$ be a word of length $n$ in $0$'s and $1$'s. We denote by $wT$  the subtree of $T$ which contains all infinite paths which start with $w$, that is, the path space of $wT$ is given by
  $$\cP(wT) = \{s_1s_2\cdots s_m \cdots \in \cP(T) \mid s_i = w_i \textrm{ for }1 \leq i \leq n\}.$$
Clearly every path in $\cP(wT)$ contains the vertex in $V_n$ labelled by $w$. The set $\cP(wT)$ is a clopen subset of $\cP(T)$. There is an obvious homeomorphism 
   $$\psi_w: \cP(wT) \to \cP(T): w_1 w_2 \cdots w_n s_{n+1} s_{n+2} \cdots  \mapsto s_{n+1} s_{n+2} \cdots \ ,$$
which induces a homeomorphism between $wT$ and $T$ which preserves the paths. Also, we can write $T =  v_0 \cup 0T \cup 1T$,  where $v_0$ is the unique vertex in $V_0$.
  
We now explain the recursive definition of homeomorphisms of $\cP(T)$. Let $\sigma$ be the non-trivial permutation of a set of two elements, that is, if this set is $\{0,1\}$, then $\sigma$ interchanges $0$ and $1$. Denote by $\Aut(T)$ the automorphism group of the tree $T$, that is, an element $h \in  \Aut(T)$ is a homeomorphism of $\cP(T)$ which for all $n \geq 1$ restricts to a permutation of $V_n$. Let $h,g \in \Aut(T)$, then define  an  automorphism  $(h,g)$ of $T$ by declaring that it restricts to the trivial permutation of $V_1$, and for every $s = s_1 s_2 \cdots \subset 0T$ one has
  \begin{equation}\label{eq-leftmap}
  (h,g)(s) = \psi_0^{-1} \circ h \circ \psi_0(s) \ , 
  \end{equation}
and for every $s \subset 1T$ one has
  \begin{equation}\label{eq-rightmap}
  (h,g)(s) = \psi_1^{-1} \circ g \circ \psi_1(s) \ . 
  \end{equation}
That is, the automorphism $h$ is applied to the branch of the tree through the first vertex $\{0\}$ of $V_1$, and $g$ is applied to the branch of the tree through the second vertex $\{1\}$ of $V_1$. We compose the maps on the left.

We now give an example of a recursive definition of a homeomorphism of $T$. Define $a$ such that 
   \begin{align}\label{eq-odom}a = (a,1)\sigma, \end{align}
 where $1$ denotes the identity in $\Aut(T)$. Then for every $s = s_1 s_2 \cdots $ the map $\sigma$ maps $s_1$ to $(s_1+1) \mod 2$ and fixes the rest of the sequence, so we can write, with a slight abuse of notation,
  $$\sigma (s) = \sigma(s_1) s_2 s_3 \cdots \ .$$
 Next, we apply $(a,1)$ to $\sigma(s)$. If $\sigma(s_1) = 1$, then by the rule \eqref{eq-rightmap}, the element  $a$ acts as the identity map on $\sigma(s_1)T = 1 T$. 
 If $\sigma(s_1) = 0$, then by \eqref{eq-leftmap} 
   $$a(s) = \sigma(s_1) a(s_2 s_3 \cdots) = 0 a(s_2 s_3 \cdots),$$ 
 so we need to compute $a(s_2 s_3 \cdots)$.  We compute that $\sigma (s_2) = (s_2 +1) \mod 2$, and then we have to apply $a$ to $\sigma(s_1) \sigma(s_2) s_3 \cdots = 0 \sigma(s_2) s_3 \cdots$, that is, the image of $s$ is computed recursively. For example,
   \begin{align*} a(0001110^\infty) = 1001110^\infty, \quad a(11001^\infty) = 00101^\infty, \quad a(1^\infty) = 0^\infty, \end{align*}
 where $1^\infty$ denotes an infinite sequence of $1$'s, and $0^\infty$ denotes an infinite sequence of $0$'s. 
 
The element $a$, defined by \eqref{eq-odom}, generates a group isomorphic to the integers $\mZ$, which acts freely and transitively on every vertex set $V_n$, $n \geq 1$. Such an arboreal action is  often called an \emph{odometer}, a \emph{standard odometer} or an \emph{adding machine} in the literature \cite{Nekrashevych2005,Pink2013}. In this paper, we use the term ``odometer'' in a more general sense; that is, an odometer is an equicontinuous minimal action of \emph{any} group $G$ on a Cantor set $\fX$. To distinguish between the two notions, we use the term \emph{cyclic odometer} for the action generated by a single element $a$ as described above.
 
 \subsection{Rigidity of Cantor actions} \label{subsec-oecocycles}
 
 In this section we illustrate the discussion in Section \ref{subsec-rigiditynotions} of various definitions of rigidity. We present an example of orbit equivalent actions, where one of the actions is free, and another one is stable but not topologically free. This example motivates the introduction of the notion of \emph{return equivalence} in Section \ref{subsec-rigiditynotions}.

Let $T$ be a binary tree with the Cantor set boundary $\cP(T)$, and set $\fX = \cP(T)$. Let $a_1 = (a_1,1)\sigma$ be the cyclic odometer defined in \eqref{eq-odom}. Let $a_2 = (a_1,1)$, that is, $a_2$ acts as a cyclic odometer on the clopen set $\cP(0T)$, and as the identity map on $\cP(1T)$.  We also note that $a_2|\cP(0T) = a_1^2|\cP(0T)$.
 
 Define $H = \langle a_1 \rangle$  and $G = \langle a_1,a_2 \rangle $, so we have the actions $(\fX,G)$ and $(\fX,H)$, where we omit $\Phi$ and $\Psi$ from the notation since $G$ and $H$ are already defined as subgroups of $Homeo(\fX)$. Both actions are minimal, since the orbit of any path $(v_n)_{n \geq 1} \in \cP(T)$ under the powers of $a_1$ is dense in $\fX$. It is easy to see that the group $G$ is non-commutative.
 
The action $(\fX,G)$ is not topologically free, as $a_2$ acts as the identity on the clopen set $\cP(1T)$. On the other hand,   $(\fX,H)$ is a free action of a cyclic odometer.

 We now show the actions $(\fX,G)$ and $(\fX, H)$   are continuously orbit equivalent. Let $h: \fX \to \fX$ be the identity map, then we define  the maps  $\alpha: \fX \times G \to H$ and $\beta: \fX \times H \to G$ (as in Section \ref{subsec-coe})   as follows.
 
  Since $a_2$ acts on $0T$ as $a_1^2$ and on $1T$ as the identity map, for any $x \in \fX$ the orbits of the action of $G$ and $H$ coincide as sets. Since the action of $H$ is free, given $x \in \fX$ and $g \in G$ there exists a unique element $h \in H$ such that $g \cdot x = h \cdot x$. We define $\alpha(g,x) = h$. We show that $\alpha$ is constant when restricted to $\cP(0T)$, and to $\cP(1T)$.
 
 So let $\cU = \{\cP(0T), \cP(1T)\}$ be a partition of $\cP(T)$. Any element $g \in G$ can be written as a finite word $g = a_1^{k_n}a_{2}^{k_{n-1}} \cdots a_1^{k_{2}}a_2^{k_1}$, with $k_j \in \mZ$ for $1 \leq j \leq n$. The actions of any power of $a_2$ and of even powers of $a_1$ preserve $\cP(wT)$, and the actions of odd powers of $a_1$ interchange $\cP(0T)$ and $\cP(1T)$. It follows that $x$ and $y$ are in the same set of the partition $\cU$ if and only if $a_2^{k_1} \cdot x$ and $a_2^{k_1}(y)$ are in the same set of $\cU$. Implementing induction on the number of elements in the word decomposition of $g$, we obtain that  $x$ and $y$ are in the same set of $\cU$ if and only if $g \cdot x$ and $g \cdot y$ are in the same set of $\cU$.

 Define the map $\beta: H \times \fX \to G$ by $\beta(a_1^k,x) = a_1^k$. The map $\beta$ is constant on $\fX$, so is continuous. Thus for $h = id$, $\alpha$ and $\beta$ yield  a continuous orbit equivalence between the actions $(\fX,G)$ and $(\fX,H)$.  By Theorem \ref{thm-main3} the action $(\fX,G)$ is stable.

  We also note that the map $\alpha_x: G \to H$ obtained by specializing $\alpha_x = \alpha(-,x)$ for  $x \in \cP(T)$ is not injective. For example, for $x \in \cP(T)$ we have
  $$\alpha_x^{-1}(a_1^{2})  \supset \{a_1^{2},a_2\} \ .$$ 
 Note that this example can be considered as an explicit form of the construction in    \cite[Example~A.4]{HL2018b}.
 
 \subsection{Dynamically wild actions}\label{subsec-adjointactions}
 
 Let $T$ be a binary tree with the Cantor set boundary $\cP(T)$, we set $\fX = \cP(T)$. Let $r \geq 3$, $2 \leq s \leq r$, and consider a group $G$ generated by the homeomorphisms
  \begin{equation}\label{eq-generatorss1} 
  a_1 = \sigma, a_{s+1} = (a_s,a_r), a_i = (a_{i-1},1) \textrm{ for } 3 \leq i \leq s, \, s+1 \leq i \leq r \ . 
  \end{equation}
 
Groups generated by the recursive rules \eqref{eq-generatorss1} arise as iterated monodromy groups of quadratic post-critically finite polynomials \cite{Nekrashevych2005,BN2008}. They have been extensively studied, for example, in \cite{BN2008}, and also in other works. For example, in \cite{Pink2013} it was shown that the closure of the action of an iterated monodromy group, associated to a post-critically finite quadratic polynomial with strictly pre-periodic post-critical orbit of length $r$, where the periodic part has length $r-1$, is conjugate to the closure of the action of the group $G$ generated by \eqref{eq-generatorss1}.

Let $G$ be a group generated by \eqref{eq-generatorss1}, and let $\oG$ be the closure of the action. It was shown in \cite[Theorem 1.3]{Lukina2018b} (see also \cite{Lukina2018a}) that $G$ contains non-Hausdorff elements, namely the generators $a_{s+1}, \ldots, a_r$ are non-Hausdorff. Then by   Theorem~\ref{thm-distinct=nH}, the action of the group $G$ is dynamically wild; that is, there is a proper inclusion of the direct limit groups  $\Upsilon_c(\Phi) \subset \Upsilon_s(\Phi)$.
 
\subsection{Infinite stabilizer group}\label{subsec-unboundedexample}

 Theorem \ref{thm-distinct=nH} shows that  if a group $G$, or more generally its closure $\oG$ has a non-Hausdorff element, then the groups $K_0 \subset K_1 \subset K_2 \subset \cdots \cD_x$ in the stabilizer chain $\cK(\cU)$ must be infinite. In specific examples it may be difficult to compute the groups $K_\ell$, $\ell \geq 1$, explicitly. We now give an example of an action for which this can be done.

 Let $T$ be a binary tree with the Cantor set boundary $\cP(T)$, we set $\fX = \cP(T)$, and let $G$ be generated by \eqref{eq-generatorss1} with $r = 2s$ and $s \geq 2$.

By definition $a_2,\ldots, a_s, a_{s+2}, \ldots, a_{2s}$ are the identity on the clopen set $\cP(1T)$ of all sequences starting with digit $1$. Since $a_{s+1} = (a_s,a_{2s})$, and $a_{2s}$ is the identity on $\cP(1T)$, $a_{s+1}$ is the identity on the set $\cP(11T)$ of all sequences starting with a finite word $11$. 

Choose a path $x \in \cP(11T)$, that is, $x = 11x_3x_4 \ldots$, and for $\ell \geq 1$ let $U_\ell$ be a clopen neighborhood
  $$U_\ell = \{s_1s_2 s_3\ldots \mid x_i = s_i \textrm{ for } 1 \leq i \leq \ell \}.$$
Then every $U_\ell \subset \cP(11T)$, and so the generators $a_{2}, \ldots, a_{2s}$ act trivially on $U_\ell$, in other words,
  $$\langle a_2,\ldots, a_s, a_{s+1}, \ldots, a_{2s}\rangle \subset K_\ell.$$
Consider the compositions $a_ia_j \in K_\ell$, with $i < j$. By \cite[Proposition 3.1.9]{Pink2013} if $j = i+s$, then $a_ia_j$ has infinite order. In particular,  $a_sa_{2s}$ has infinite order. It follows that $G \cap K_\ell$ is infinite, and so $K_\ell$ is a Cantor group.



\begin{thebibliography}{10}

 

  \bibitem{ALC2009}
{J.~{\'A}lvarez L{\'o}pez and A.~Candel},
\newblock {\it Equicontinuous foliated spaces},
\newblock {\bf Math. Z.}, 263:725--774, 2009.

\bibitem{ALM2016}
{J.~{\'A}lvarez L{\'o}pez and M.~Moreira Galicia},
\newblock {\it Topological {M}olino's theory},
\newblock {\bf Pacific. J. Math.}, 280:257--314, 2016.

 
     
\bibitem{Auslander1988}
{J.~Auslander},
\newblock {\bf Minimal flows and their extensions},
\newblock {North-Holland Mathematics Studies}, Vol. 153, {North-Holland Publishing Co., Amsterdam}, 1988.


\bibitem{BN2008}
{L.~Bartholdi and V.~Nekrashevych}
\newblock{\it Iterated monodromy groups of quadratic polynomials, I},
\newblock{\bf Groups Geom. Dyn.}, 2:309-336, 2008.

 \bibitem{BezuglyiMedynets2008}
{S.~Bezuglyi and K.~Medynets},
\newblock {\it Full groups, flip conjugacy, and orbit equivalence of {C}antor minimal systems},
\newblock {\bf Colloq. Math.}, 110:409--429, 2008.

\bibitem{Boyle1983}
{M.~Boyle},
\newblock {\bf Topological orbit equivalence and factor maps in symbolic dynamics},
\newblock {Ph.D. Thesis}, University of Washington, 1983.     

 \bibitem{BoyleTomiyama1998}
{M.~Boyle and J.~Tomiyama},
\newblock {\it Bounded topological orbit equivalence and {$C^*$}-algebras},
\newblock {\bf J. Math. Soc. Japan}, 50:317--329, 1998.     

 \bibitem{BCFS2014}
{J.~Brown, L.~Clark, L.~Orloff and C.~Farthing},
\newblock {\it Simplicity of algebras associated to \'etale groupoids},
\newblock {\bf Semigroup Forum}, 88:433--452, 2014.     

    
\bibitem{ClarkHurder2013}
{A.~Clark and S.~Hurder},
\newblock {\it Homogeneous matchbox manifolds},
\newblock {\bf Trans. A.M.S.}, 365:3151--3191, 2013.

\bibitem{CHL2018a}
{A.~Clark, S.~Hurder and O.~Lukina},
\newblock {\it Classifying matchbox manifolds},
\newblock {\bf Geom \& Top},  23:1-38, 2019;  {arXiv:1311.0226}.


 \bibitem{CFW1981}
{A.~Connes, J.~Feldman and B.~Weiss},
\newblock {\it An amenable equivalence relation is generated by a single transformation},
\newblock {\bf Ergodic Theory Dynamical Systems}, 1:431--450, 1981.


 \bibitem{CortezPetite2008}
{M.-I.~Cortez and S.~Petite},
\newblock {\it $G$-odometers and their almost one-to-one extensions},
\newblock {\bf J. London Math. Soc.}, 78(2):1--20, 2008. 

      
\bibitem{CortezMedynets2016}
{M.I.~Cortez and K.~Medynets},
\newblock {\it Orbit equivalence rigidity of equicontinuous systems},
\newblock {\bf J. Lond. Math. Soc. (2)}, 94:545--556, 2016.

\bibitem{CortezPetite2018}
{M.I.~Cortez and S.~Petite},
\newblock {\it On the centralizers of minimal aperiodic actions on the Cantor set},
\newblock {\it preprint}; {arXiv:1807.04654}.


 \bibitem{deCornulier2014}
{Y.~de~Cornulier},
\newblock {\it Groupes pleins-topologiques (d'apr\`es {M}atui, {J}uschenko, {M}onod, {$\ldots$})},
\newblock {\bf Ast\'erisque}, 361:183--223, 2014.

   
 \bibitem{Do2005}
{T.~Downarowicz}, 
\newblock   {\it Survey of odometers and {T}oeplitz flows}, 
\newblock in {\bf Algebraic and topological dynamics},
\newblock  {Contemp. Math.}, Vol. 385: 7--37,  {Amer. Math. Soc.}, {Providence, RI}, {2005}. 
  
 
  
\bibitem{DHL2016a}
{J.~Dyer, S.~Hurder and O.~Lukina},
\newblock {\it The discriminant invariant of Cantor group actions},
\newblock {\bf Topology Appl.}, 208: 64--92, 2016.


\bibitem{DHL2016b}
{J.~Dyer, S.~Hurder and O.~Lukina},
\newblock {\it Growth and homogeneity of matchbox manifolds},
\newblock {\bf Indag. Math.}, 28:145--169, 2017.

\bibitem{DHL2016c}
{J.~Dyer, S.~Hurder and O.~Lukina},
\newblock {\it Molino theory for matchbox manifolds},
\newblock {\bf Pacific J. Math.}, 289:91-151, 2017.

 
\bibitem{EilenbergSteenrod1952}
{S.~Eilenberg and N.~Steenrod},
\newblock {\bf Foundations of algebraic topology},
\newblock {Princeton University Press, Princeton, New Jersey}, 1952.

     
\bibitem{Ellis1960}
{R.~Ellis},
\newblock {\it A semigroup associated with a transformation group},
\newblock {\bf Trans. Amer. Math. Soc.}, 94:272--281, 1969.


 \bibitem{EllisGottschalk1960}
{R.~Ellis and W.H.~Gottschalk},
\newblock {\it Homomorphisms of transformation groups},
\newblock {\bf Trans. Amer. Math. Soc.}, 94:258--271, 1969.

  \bibitem{EMT1977}
{D.B.A.}~Epstein, {K.C.}~Millet, and {D.}~Tischler,
\newblock {\it Leaves without holonomy},
\newblock {\bf Jour. London Math. Soc.}, 16:548--552, 1977.

   
\bibitem{Exel2011}
{R.~Exel},
\newblock {\it Non-{H}ausdorff \'etale groupoids},
\newblock {\bf Proc. Amer. Math. Soc.}, 139:897--907, 2011.
     
\bibitem{FO2002}
{R.~Fokkink and L.~Oversteegen},
\newblock {\it Homogeneous weak solenoids},
\newblock {\bf Trans.  Amer. Math. Soc.}, 354(9):3743--3755, 2002.

\bibitem{Furman1999}
{A.~Furman},
\newblock {\it Orbit equivalence rigidity},
\newblock {\bf Ann. of Math. (2)}, 150:1083--1108, 1999.


\bibitem{Gaboriau2010}
{D.~Gaboriau},
\newblock {\it Orbit equivalence and measured group theory},
\newblock {\bf Proceedings of the {I}nternational {C}ongress of {M}athematicians. {V}olume {III}}, 
\newblock {Hindustan Book Agency, New Delhi}, 2010, pages 1501--1527.


       

\bibitem{GPS1999}
{T.~Giordano, I.~Putman and C.~Skau},
\newblock {\it Full groups of {C}antor minimal systems},
\newblock {\bf Israel J. Math.}, 111:285--320, 1999.

     
\bibitem{GMPS2010}
{T.~Giordano, H.~Matui, I.~Putnam, and C.~Skau},
\newblock{\it Orbit equivalence for {C}antor minimal {$\mathbb Z^d$}-systems},
\newblock{\bf Invent. Math.}, 179: 119--158, 2010.


\bibitem{GPS2019}
{T.~Giordano, I.~Putman and C.~Skau},
\newblock {\it {$\mathbb Z^d$}-odometers and cohomology},
\newblock {\bf Groups Geom. Dyn.}, 13:909--938, 2019.  
  
           
\bibitem{GW1995}
{E.~Glasner and B.~Weiss},
\newblock {\it Weak orbit equivalence of {C}antor minimal systems},
\newblock {\bf Internat. J. Math.}, 6:559--579, 1995.

  \bibitem{Glasner2007}
{E.~Glasner},
\newblock {\it Enveloping semigroups in topological dynamics},
\newblock {\bf Topology Appl.}, 154:2344--2363, 2007.

      
           
\bibitem{Grigorchuk1984}
{R.I.~Grigorchuk},
\newblock {\it Degrees of growth of finitely generated groups and the theory of invariant means},
\newblock  {\bf Izv. Akad. Nauk SSSR Ser. Mat.},  48:939--985, 1984.


 \bibitem{Haefliger1985}
{A.~Haefliger},
\newblock {\it Pseudogroups of local isometries}, in Differential Geometry (Santiago de Compostela, 1984), edited by L.A. Cordero,
\newblock {\bf Res. Notes in Math.}, 131:174--197, Boston, 1985.

        
\bibitem{HL2018a}
{S.~Hurder and O.~Lukina},
\newblock {\it Wild solenoids},
\newblock {\bf Transactions A.M.S.}, 371:4493-4533, 2019.
 
\bibitem{HL2018b}
{S.~Hurder and O.~Lukina},
\newblock {\it Orbit equivalence and   classification of weak solenoids},
\newblock {\bf Indiana Univ. Math. J.}, to appear; {arXiv:1803.02098}.
 

 
   
 \bibitem{KN1969}
{S.~Kobayashi and K.~Nomizu},
\newblock {\bf Foundations of differential geometry. {V}ol. {II}},
\newblock {Interscience Tracts in Pure and Applied Mathematics, No. 15 Vol. II },
 {Interscience Publishers John Wiley \& Sons, Inc., New York-London-Sydney}, 1969.
  
 \bibitem{Li2018}
{X.~Li},
\newblock {\it Continuous orbit equivalence rigidity},
\newblock {\bf Ergodic Theory Dynam. Systems}, 38:1543--1563, 2018.

 \bibitem{Lubotzky1993}
{A.~Lubotzky},
\newblock {\it Torsion in profinite completions of torsion-free groups},
\newblock {\bf Quart. J. Math. Oxford Ser. (2)} 44:327--332, 1993.

  \bibitem{Lukina2018a}
{O.~Lukina},
\newblock {\it Arboreal Cantor actions},
\newblock {\bf J. Lond. Math. Soc. (2)}, 99:678--706, 2019.
 
\bibitem{Lukina2018b}
{O.~Lukina},
\newblock {\it Galois groups and Cantor actions},
\newblock {\it submitted}; {arXiv:1809.08475}.    

        
      
      
 \bibitem{McCord1965}
{C.~McCord},
\newblock {\it Inverse limit sequences with covering maps},
\newblock {\bf Trans. A.M.S.}, 114:197--209, 1965.

\bibitem{Munkres1984}
{J.~Munkres},
\newblock {\bf Elements of algebraic topology},
\newblock {Addison-Wesley Publishing Company, Menlo Park, CA}, 1984. 
       

      
\bibitem{Nekrashevych2005}
{V.~Nekrashevych},
\newblock {\bf Self-similar groups},
\newblock {Mathematical Surveys and Monographs}, Vol. 117, 
\newblock {American Mathematical Society, Providence, RI}, 2005.  
 
\bibitem{Nekrashevych2018}
{V.~Nekrashevych},
\newblock {\it Palindromic subshifts and simple periodic groups of intermediate growth},
\newblock {\bf Ann. of Math. (2)}, 187:667--719, 2018.

 
\bibitem{Pink2013}
{R.~Pink},
\newblock {\it Profinite iterated monodromy groups arising from quadratic polynomials}; 
\newblock {\it preprint}; arXiv: 1307.5678.

    
\bibitem{Renault1980}
{J.~Renault},
\newblock {\it A groupoid approach to {$C^*$}-algebras},
\newblock {\bf Lecture Notes in Math.},  vol. 793, 1980.


        
\bibitem{Renault2006}
{J.~Renault},
\newblock {\it Transverse properties of dynamical systems},
 \newblock in {\bf Representation theory, dynamical systems, and asymptotic combinatorics}, 
 \newblock {Amer. Math. Soc. Transl. Ser. 2}, Vol. 217, {Amer. Math. Soc., Providence, RI}, 2006, pages 185--199. 
      
\bibitem{Renault2008}
{J.~Renault},
\newblock {\it Cartan subalgebras in {$C^*$}-algebras},
 \newblock {\bf Irish Math. Soc. Bull.}, 61:29--63, 2008.
 
\bibitem{RT1971b} 
{J.T.~Rogers, Jr. and J.L.~Tollefson}, 
\newblock {\it Homogeneous inverse limit spaces with non-regular covering maps as bonding maps},
\newblock {\bf Proc. Amer. Math. Soc.} {\bf 29}: 417--420, 1971.

       
 \bibitem{Schori1966}
{R.~Schori},
\newblock {\it Inverse limits and homogeneity},
\newblock {\bf Trans. A.M.S.}, 124:533--539, 1966.
       

 \bibitem{Winkelnkemper1983}
{E.~Winkelnkemper}, 
\newblock {\it The graph of a foliation}, 
\newblock {\bf Ann. Global Ann. Geo.}, 1:51--75, 1983.
      
   
\end{thebibliography}
\end{document}